\documentclass[EJP, preprint]{ejpecp} % replace ECP by EJP if needed.
\SHORTTITLE{Subgraph count of dynamic ER random graphs}

\TITLE{Functional Central Limit Theorem\\ for the simultaneous subgraph count of \\ dynamic Erd\H{o}s--R\'enyi random graphs}
    
\AUTHORS{%
  Rajat Subhra Hazra\footnote{Mathematical Institute, Leiden University, P.O. Box 9512,
2300 RA Leiden
    \EMAIL{r.s.hazra@math.leidenuniv.nl}}
  \and %% remove this line and below if single author
  Nikolai Kriukov\footnote{Korteweg-de Vries Institute for Mathematics, University of Amsterdam, Science Park 904, 1098 XH Amsterdam, The Netherlands \BEMAIL{n.kriukov@uva.nl}}
  \and
  Michel Mandjes \footnote{Mathematical Institute, Leiden University, P.O. Box 9512,
2300 RA Leiden, Korteweg-de Vries Institute for Mathematics, University of Amsterdam, Amsterdam\BEMAIL{m.r.h.mandjes@math.leidenuniv.nl}}}
 % Separate items with ;

\AMSSUBJ{05C80; 60F05} % Edit. Separate items with ;
\SUBMITTED{} % Edit.
\ACCEPTED{} % Edit.

\VOLUME{0}
\YEAR{2023}
\PAPERNUM{0}
\DOI{10.1214/YY-TN}

\ABSTRACT{In this paper we consider a dynamic {Erd\H{o}s--R\'{e}nyi} random graph with independent identically distributed edge processes. Our aim is to describe the joint evolution of the entries of a subgraph count vector. The main result of this paper is a functional central limit theorem: we establish, under an appropriate centering and scaling, the joint functional convergence of the vector of subgraph counts to a specific multidimensional Gaussian process. The result holds under mild assumptions on the edge processes, most notably a Lipschitz-type condition.}

\usepackage{amsmath}
\usepackage{bbm}
\usepackage{amssymb}
\usepackage{mathtools}
\usepackage{stmaryrd}

\numberwithin{equation}{section}

\usepackage{color}
\usepackage[dvipsnames]{xcolor}

\definecolor{c20}{rgb}{0.,0.7,0.}
\definecolor{c30}{rgb}{0.,0.,1.}
\definecolor{c40}{rgb}{1,0.1,0.7}
\definecolor{c50}{rgb}{1,0,0}
\definecolor{c60}{rgb}{1,0.9,0.1}

\newcommand{\abs}[1]{\left\lvert #1 \right\rvert}

\newcommand{\sprod}[1]{\langle#1\rangle}

\newcommand{\E}[1]{\mathbb{E}\left\{ #1\right\}}

\newcommand{\pk}[1]{\mathbb{P} \left\{ #1 \right \} }

\newcommand{\R}{\mathbb{R}}

\newcommand{\N}{\mathbb{N}}
\newcommand{\Z}{\mathbb{Z}}

\newcommand{\BQN}{\begin{eqnarray}}
\newcommand{\EQN}{\end{eqnarray}}
\newcommand{\BQNY}{\begin{eqnarray*}}
\newcommand{\EQNY}{\end{eqnarray*}}

\newcommand{\BS}{\begin{sat}}
\newcommand{\ES}{\end{sat}}
\newcommand{\BT}{\begin{theo}}
\newcommand{\ET}{\end{theo}}
\newcommand{\BK}{\begin{korr}}
\newcommand{\EK}{\end{korr}}

\newcommand{\BD}{\begin{de}}
\newcommand{\ED}{\end{de}}

\newcommand{\BIT}{\begin{itemize}}
\newcommand{\EIT}{\end{itemize}}
\newcommand{\BDI}{\begin{description}}
\newcommand{\EDI}{\end{description}}

\newcommand{\BRM}{\begin{remarks}}
\newcommand{\ERM}{\end{remarks}}

\newcommand{\BEL}{\begin{lem}}
\newcommand{\EEL}{\end{lem}}

\newtheorem{theo}{Theorem}[section]
\newtheorem{sat}[theo]{Proposition}
\newtheorem{de}[theo]{Definition}
\newtheorem{lem}[theo]{Lemma}

\newtheorem{ass}[theo]{Assumption}
\newtheorem{korr}[theo]{Corollary}
\newtheorem{remarks}[theo]{Remarks}
\newtheorem{prop}[theo]{Proposition}

\newcommand{\nelem}[1]{{Lemma \ref{#1}}}
\newcommand{\neprop}[1]{{Proposition \ref{#1}}}
\newcommand{\netheo}[1]{{Theorem \ref{#1}}}
\newcommand{\nekorr}[1]{{Corollary \ref{#1}}}

\newcommand{\COM}[1]{}

\def\td{\text{\rm d}}

\newcommand{\QED}{\hfill $\Box$}

\newcommand{\kb}[1]{\boldsymbol{#1}}
\newcommand{\vk}[1]{\kb{#1}}

\begin{document}

%%%%%%%%%%%%%%%%%%%%%%%%%%%%%%%%%%%%%%%%%%%%%%%%%%%%%%%%%%%%%%%%%%%
%%                                                               %%
%% No need for \maketitle.                                       %%
%%                                                               %%
%%%%%%%%%%%%%%%%%%%%%%%%%%%%%%%%%%%%%%%%%%%%%%%%%%%%%%%%%%%%%%%%%%%

%%%%%%%%%%%%%%%%%%%%%%%%%%%%%%%%%%%%%%%%%%%%%%%%%%%%%%%%%%%%%%%%%%%
%%                                                               %%
%% Please replace what follows by the body of your article       %%
%% (up to the bibliography):                                     %%
%%                                                               %%
%%%%%%%%%%%%%%%%%%%%%%%%%%%%%%%%%%%%%%%%%%%%%%%%%%%%%%%%%%%%%%%%%%%

\section{Introduction}\label{section:introduction}

Over the past decades, inspired by applications in a broad range of disciplines (such as biology, physics, economics, and computer science), various types of random graph models have been introduced and analyzed.  
Random graph models are specified via an underlying probabilistic mechanism. 
The {Erd\H{o}s--R\'enyi} random graph \cite{erdos1959random, gilbert1959random}, in which there are $N$ vertices and in which every possible vertex pair is connected independently with a given probability $p\in(0,1)$, is arguably the most well-known.

Most literature on random graphs focuses on {\it static} models, which represent a single instance of the random graph without accounting for them changing over time. However, real-world networks typically {\it do} evolve over time, motivating the study of {\it dynamic} random graphs, particularly dynamic {Erd\H{o}s--R\'enyi} random graphs (e.g., \cite{zhang2017random, mandjes2019dynamic, holme2012temporal, holme2015modern, braunsteins2023sample}).  In these models, edges evolve independently over time, alternating between being present and absent. Essentially all properties that have been studied for static random graphs can be analyzed for their dynamic counterparts as well. For instance, in the context of dynamic random graphs there is now a notion of a {\it giant component process}, and a notion of a {\it spectrum process}. In the present paper we concentrate on the {\it subgraph counts process}, in that our aim is to describe the joint evolution of the entries of a subgraph count vector. We thus capture, for instance, the joint distribution of the number of edges and the number of triangles over time. 

\medskip

Subgraph counts play a pivotal role in network analysis, as they encode crucial structural information about the geometry underlying the graph (see e.g.\ \cite{frank1986markov, kendall1940method, moran1947method, holland1971transitivity, holland1976local, harary1979matrix, frank1980transitivity, frank1982cluster, andrade:bhattacharya:2025}). 
There is a mature body of work that probabilistically analyzes subgraph counts in the static case.
We proceed by providing a brief account of the existing literature in this area, with a focus on the central limit regime. It should be borne in mind, though, that this branch of research is relatively mature, rendering our overview inherently non-exhaustive. 

 The first results showing asymptotic normality for subgraph counts in a static {Erd\H{o}s--R\'enyi} random graph appeared already in \cite{erdos1960evolution}: an appropriately centered and normalized version of the subgraph count converges, as $N\to\infty$, to a zero-mean Gaussian random variable. These results were further generalized to the regime where the edge probability $p=p_N$ tends to zero as the number of vertices $N$ grows large; see e.g.\ \cite{nowicki1988subgraph} that covers specific types of subgraphs. Finally, in \cite{rucinski1988small} necessary and sufficient conditions on the decay of $p_N$ were established for obtaining asymptotic normality for any subgraph count.

A separate stream in the literature fully concentrates on the case of a constant edge probability. In this context we mention \cite{nowicki1989asymptotic}, which proves the joint asymptotic normality for different subgraph counts. More recently, the limiting distributions for a single subgraph count in the graphon setting were addressed in~\cite{hladky2021limit,bhattacharya2023fluctuations}. 

\medskip

Importantly, there seem to be no asymptotic normality results covering the setting of random graphs that evolve over time, even in the simplest case in which the edge processes do not depend on the number of vertices $N$. The main objective of this paper is to settle this, by establishing a joint functional central limit theorem for dynamic {Erd\H{o}s--R\'enyi} random graphs. Concretely, we identify conditions under which an appropriately centered and normalized version of a vector of subgraph counts has, at the process level, a Gaussian limit. Moreover, we explicitly characterize the limiting Gaussian process. 
{The type of convergence considered is in distribution, in $\mathbb{D}([0,T],\R^{m})$, which is the set of all c\`adl\`ag functions $f(\cdot)\colon [0,T]\to \R^m$, equipped with $J_1$ distance, i.e., for any functions $f(\cdot),\,g(\cdot)\in\mathbb{D}([0,T],\R^{m})$
\begin{align*}
    d\bigl(f(\cdot),\,g(\cdot)\bigr) = \inf_{\lambda\in\Lambda}\max\left(\sup_{t\in[0,T]}\abs{\lambda(t) - t},\,\sup_{t\in[0,T]}\| f(t) - g(\lambda(t))\|\right),
\end{align*}
for some finite, deterministic horizon $T$, where $\Lambda$ is the class of strictly increasing, continuous mappings of $[0,T]$ onto itself. We refer to e.g.\ \cite[Section 12]{billingsley2013convergence} for a more detailed account}.

In our setup, we do not focus on one specific type of edge dynamics, in that we  aim to work with a specification of the dynamic {Erd\H{o}s--R\'enyi} random graph that is as general as possible. We follow the `traditional' path in proving process-level convergence: we first establish finite-dimensional convergence, followed by a tightness step that yields the desired functional convergence. In order to accommodate this tightness, we have to impose some assumptions on the edge dynamics. More concretely, in our framework the probabilities of the edge switching once or twice (in an interval of given length, that is) are assumed to obey a Lipschitz-type condition. At the methodological level, some of the steps of our analysis follow their counterparts in \cite{hazra2024principle}, where the principal eigenvalue process of the graph's adjacency matrix is analyzed in the regime in which the edge processes do not depend on $N$.

\medskip

The main result of this paper, \netheo{main}, establishes that the limiting behavior (as the size of our graph grows large) of the subgraph count process is Gaussian, and strictly correlated with the subgraph count for the {least common subgraph} of our graph. This {least common subgraph} is defined as a subgraph with at least one edge and the asymptotically smallest expected number in our dynamic {Erd\H{o}s--R\'{e}nyi} random graph. For instance, if the edge process does not depend on the size of the random graph (as before denoted by $N$), the expected subgraph count asymptotically scales as $N^{v}$, where $v$ is the size of the subgraph under consideration. In this scenario, the {least common subgraph} for any graph with at least one edge is a single edge, and thus every subgraph count process is strictly correlated with the edge count. This implies that all possible subgraph counts are strictly correlated with one another. We provide a detailed exploration of this phenomenon in \nekorr{main_corr}. 

{It is noted that in the regime where the `edge probabilities' decay as $N$ grows to infinity, the speed of this convergence can influence which subgraph becomes the {least common one}. {Here the notion of `edge probability' refers to the probability that a given edge exists; observe that these probabilities are functions of time.} In the setting of decaying edge probabilities (as $N\to\infty$), two additional assumptions are in place, alongside the Lipschitz-type condition. In the first place, we assume that the edge probabilities vanish, as $N\to\infty$, uniformly over the time interval $[0,T]$ considered. This is for convenience, because without this assumption we would have to split the time interval into subintervals that are to be dealt with separately.  Second, we focus solely on those subgraphs for which the probability of non-appearance in the dynamic random graph tends to zero as $N\to\infty$. This assumption, which links the asymptotic behavior of the edge probabilities to the structure of the subgraph, is adopted from \cite{rucinski1988small} and is discussed in greater detail in Section~\ref{section:graph_notations}.

We also like to highlight a specific consequence of our results: \netheo{main} reveals that if two different graphs have distinct {least common} subgraphs, their respective subgraph count processes are asymptotically independent.}

\medskip

\noindent {\it Notation ---} Throughout the paper we consistently use the following notations. For any {natural number} $n$, we write $\sprod{n}$ for the set $\{1,\ldots,n\}$. All vectors are written in bold and, vice versa, all bold symbols refer to vectors. All entries of a vector share the same letter, but do not use bold font and have a subscript which represents the index of the entry. For example if $\vk b\in\R^d$, then $\vk b = (b_1,\ldots,b_d)^{\top}$. If a vector has a subscript itself, for each element of this vector we first write a subscript of the vector, and after the subscript of the element. For example, if $\vk b_i\in\R^d$, then $\vk b_i = (b_{i,1},\ldots,b_{i,d})^{\top}$. For any finite set $A$ we denote the cardinality of $A$ as $\abs{A}$.

\medskip

\noindent{\it Outline ---} This paper is organized as follows. In Section \ref{section:model} we introduce the crucial properties of the edge process. 
In particular, this section explicitly states the condition that we impose,  along with some natural processes that meet them.
{In Section \ref{section:graph_notations} we introduce a series of key objects which we heavily rely on when proving the main result. Section \ref{section:main_results} presents our main result: we state the multivariate functional central limit theorem for a (vector-valued) subgraph count process.} Sections \ref{section:proof_of_main}--\ref{section:proof_of_tightness} are dedicated to the various steps of the proof of \netheo{main}: in  Section \ref{section:proof_of_main} we show three crucial lemmas to build up the proof of \netheo{main}, and Sections \ref{section:proof_of_covariance}--\ref{section:proof_of_tightness} cover the proofs of these three lemmas, respectively.  Appendix~\ref{section:appendix} contains proofs of various auxiliary results.

\section{Model}\label{section:model} 
We start by constructing our dynamic {Erd\H{o}s--R\'{e}nyi} random graph process, in the sequel denoted by $G_N(t)$, with $N\in\N$ vertices over the time interval $t\in[0,T]$ for some given time horizon $T>0$. Each pair of vertices $u,v\in\sprod{N}$ can be connected by an edge (which we associate with the state `1') or not (which we associate with the state `0'). Throughout this paper, the resulting edge process of the vertex-pair $(u,v)$ is denoted by $(a_{N,u,v}(t))_{t\in[0,T]}$. We restrict ourselves to what could be called {\it homogeneous simple graphs}, i.e., across all $u,v\in\sprod{N}$ with $u<v$ and for any fixed $N\in\N$, the processes $(a_{N,u,v}(t))_{t\in[0,T]}$ are independent, identically distributed, and $G_N(t)$ does not have any self-loops and double edges. This in particular entails that (i)~$a_{N,u,u}(t) = 0$ for all $N\in\N$, and $t\in[0,T]$, and (ii)~$a_{N,u,v}(t)= a_{N,v,u}(t)\in\{0,1\}$ for all $N\in\N$, $u,v\in\sprod{N}$, $u\not= v$ and $t\in[0,T]$.

\medskip

In the sequel expectations and covariances pertaining to the processes $(a_{N,u,v}(t))_{t\in[0,T]}$ play a key role. Throughout this paper, we use the following notations:
\begin{align*}
    p_N(t) &= \E{a_{N,u,v}(t)} = \pk{a_{N,u,v}(t)  =1},\\
    \kappa_N(s,t) &= \operatorname{Cov}\bigl(a_{N,u,v}(s), a_{N,u,v}(t)\bigr).
\end{align*}
% Throughout the paper we would like to consider the classical case $p_N(t)\to p_{\lim}(t)\in [0,1)$ as $N\to\infty$ However, as now we have a set of such sequences, indexed by $t\in[0,T]$, we would like all them to have the same asymptotic behaviour as $N\to\infty$. 
The papers on asymptotic normality for static {Erd\H{o}s--R\'{e}nyi} random graphs, such as the ones cited in the introduction, work with some assumption on the behavior of the edge probability in the regime that the number of vertices $N$ grows large. 
In our dynamic counterpart we have to impose a similar assumption.
\begin{ass}\label{ass:ass1}
There exists a control sequence $\varrho_N>0$ such that either $\lim_{N\to\infty}\varrho_N = 0$ or $\varrho_N\equiv 1$, and, for any given $s,t\in[0,T]$,
\begin{align} 
    \lim_{N\to\infty}\frac{p_N(t)}{\varrho_N}&=p_{\star}(t)\in[0,\infty),\label{p_asymptotics}\\
    \lim_{N\to\infty}\frac{\kappa_N(s,t)}{\varrho_N}& = \kappa_{\star}(s,t)\in[0,\infty)\label{kappa_asymptotics}
\end{align} 
for some limiting functions $p_{\star}(t)$, $\kappa_{\star}(s,t)$. In addition, the convergences in \eqref{p_asymptotics} and \eqref{kappa_asymptotics} are uniform in $s,t\in[0,T]$. 
\end{ass}

\begin{remark}{\em 
    {Under Assumption \ref{ass:ass1}, the sequence $\varrho_N$ is defined up to a constant factor, which we determine as follows. For subsequent calculations, it is important to control the limits of the sequences $N^v\varrho_N^e$ as $N\to\infty$, where $v$ and $e$ are elements of $\Z$. We initially fix the constant arbitrarily. If $\lim_{N\to\infty}N^v\varrho_N^e\in\{0,\infty\}$ for all pairs $v,e\in\Z$, the arbitrarily chosen constant remains unchanged. However, if for some $v^{\star},e^{\star}\in\Z$, the limit $\lim_{N\to\infty}N^{v^{\star}}\varrho_N^{e^{\star}}\in(0,\infty)$, we set $\varrho_N = N^{-v^{\star}/e^{\star}}$. In this case, for any $v,e\in\Z$,} \begin{align*}{
        \lim_{N\to\infty}N^v\varrho_N^e = \lim_{N\to\infty}N^{v - e\frac{v^{\star}}{e^{\star}}} = \begin{cases} 0 \quad &\text{if } e/v >e^\star/{v^{\star}}, \\
        1 \quad &\text{if } e/v = e^\star/{v^{\star}}, \\
        \infty \quad &\text{if } e/v < e^\star/{v^{\star}}. \end{cases} }
    \end{align*} 
    {Thus, in all cases, the sequence $\varrho_N$ can be chosen so that for any pair of constants $v,e\in\mathbb{Z}$, it satisfies the following property:} 
    \begin{align} 
        \lim_{N\to\infty}N^{v}\varrho_N^{e} \in\{0,1,\infty\}.\label{rho_choice}
    \end{align}
    {Technically, to check \eqref{rho_choice} we need to have $v, e\in \mathbb{N}$, but at various places in our proofs we work with powers that potentially have a negative sign; this explains why in \eqref{rho_choice} we chose to include $v, e\in \mathbb{Z}$.}
    \hfill ${\Diamond}$}
\end{remark}

\COM{
\begin{remark}\footnote{{N: Here the upper bound $p_+$ is important, but however we do not need the lower bound $p_-$}}{\em 
    As the convergence in \eqref{p_asymptotics} is uniform and the limiting function $p_{\star}(t)$ is positive and continuous on $[0,T]$, we can find positive constants $p_-, p_+$ such that for large enough $N$ and all $t\in[0,T]$, the following lower and upper bound apply:
    \begin{align}
        0<p_-\leqslant \frac{p_N(t)}{\varrho_N}\leqslant p_+<\infty
    \end{align}} \hfill ${\Diamond}$
\end{remark}
}

\begin{remark}{\em 
    We can always choose the same control sequence for both expectation and covariance function, due to the fact that
    $\kappa_N(t,t) = p_N(t) - p_N^2(t).$}
    \hfill ${\Diamond}$
\end{remark}

It is our objective to work with a setup that is as general as possible. To this end, we chose to not uniquely specify the stochastic dynamics of the edge processes. Our results require the following mild additional assumption on the processes $(a_{N,u,v}(t))_{t\in[0,T]}$. 

\begin{ass}\label{ass:ass2}
 {There exists a sequence of processes $(a_{N}(t))_{t\in[0,T]}$ and $N\in\N$ such that, for every $u,v\in\sprod{N}$ with $u\not= v$ and every $N\in\N$, the process $a_{N,u,v}(\cdot)$ has the same distribution as the generic process $a_N(\cdot)$. Additionally,} there exists a positive constant ${\mathfrak{C}}>1$ such that, for all  $0\leqslant r\leqslant  s\leqslant t\leqslant T$ and $N\in\N$,
\begin{align}
    &\pk{a_{N}(s)\not= a_{N}(r)}\leqslant {\mathfrak{C}}\,\varrho_N(s-r),\label{claim_1}\\
    &\pk{a_{N}(s)\not= a_{N}(r),\, a_{N}(t)\not= a_{N}(s)}\leqslant {\mathfrak{C}}\,\varrho_N(t-r)^2.\label{claim_2}
\end{align}
\end{ass}

\begin{remark}{\em 
    Using \eqref{claim_1} we can obtain that, for any $n\in\N$, the function $p_N(t)$ is necessarily Lipschitz-continuous with Lipschitz-constant ${\mathfrak{C}}$. Indeed,
    \begin{align}
        \abs{\,p_N(s)-p_N(r)}\leqslant {\mathfrak{C}}\,\varrho_N(s-r)\label{p_inc_bound}
    \end{align}
for all  $0\leqslant r\leqslant  s\leqslant T$. {It implies that the function $p_N(\cdot)/\varrho_N$ is continuous for any $N\in \N$. Hence, due to \eqref{p_asymptotics}, there exists constant $p_{+}\in(1,\infty)$ such that
    \begin{align}
        {\sup_{N\in\N}\sup_{t\in[0,T]}}\frac{p_N(t)}{\varrho_N}&\leqslant p_+,\label{p_star_bound}
    \end{align}}} \hfill    ${\Diamond}$
\end{remark}

    One of the most natural class of processes which satisfies Assumptions \ref{ass:ass1} and \ref{ass:ass2} is the class of  \emph{alternating  processes} {$a_N(\cdot)$, where $a_N(t)\in\{0,1\}$ for all $t\in[0,T]$}. An alternating  process is initiated by sampling the Bernoulli random variable $a_N(0)$, and then it changes between on and off at random points in time. We let $\xi_{N,x,i}$, for $x\in\{0,1\}$ and $i,N\in\N$, denote the duration of the $i$-th `sojourn time'   given that process starts at the state $x$, where such a sojourn time is either an on-time or an {off-time}, depending on the parity of $x+i$; it is typically assumed that all the on- and off-times are independent. If these random variables are in addition from a continuous distribution, then under a couple of mild restrictions such a process satisfies \eqref{claim_1} and \eqref{claim_2} by virtue of the following proposition.
    % \footnote{\tt \CLM{What we wrote down was not -renewal-, because we did not assume the on-times (off-times, resp.) to be identically distributed. That's why I left out "renewal". OK? In Example 2.7 is "renewal" OK, I suppose.}}

    \begin{prop}\label{example}
    Let {$a_N(\cdot)$} be the alternating process described above, which satisfies the following properties:
    \begin{itemize}
        \item[$\circ$]{There exists a sequence $\varrho_N$ such that the limit \begin{align*}p_{\star}(t) = \lim_{N\to\infty}\frac{p_N(t)}{\varrho_N}\end{align*} exists uniformly in $t$.}
        \item[$\circ$]{The random variables $\xi_{N,x,i}$, with $x\in\{0,1\}$ and $i\in{\mathbb N}$, are from a continuous distribution. in addition, their densities $\rho_{N,x,i}(t)$ are uniformly bounded by the constant $P$ for any $N,i\in\N$ and $x\in\{0,1\}$.}
        \item[$\circ$]{ For any $N\in\N$, $\rho_{N,0,1}<P\varrho_N$.}
    \end{itemize}
    Then $a_N(\cdot)$ satisfies Assumption $\ref{ass:ass2}$ for ${\mathfrak{C}} = (p_{\star}(0)+2)Pe^{PT}$.
    \end{prop}

    The proof of this proposition is presented in Appendix \ref{section:appendix}.

\begin{example}\label{ex:process}\em 
    As a particular example of the alternating processes that satisfy the assumptions of \neprop{example}, we can consider a stationary alternating renewal process with exponential on- and off-times. This concretely means that $\xi_{N,x,i}$ has an exponential distribution with parameter $\lambda_{N, {\rm off}}$ if $x+i$ is even and with parameter $\lambda_{N, {\rm on}}$ if $x+i$ is odd, and $p_{N}(0) = \lambda_{N, {\rm off}} / (\lambda_{N, {\rm on}}+\lambda_{N, {\rm off}})$. In this case we have for any $t\in[0,T]$
    \begin{align*}
        p_N(t) = \frac{\lambda_{N, {\rm off}}}{\lambda_{N, {\rm on}}+\lambda_{N, {\rm off}}}.
    \end{align*}
    When considering the specific parametrization $\lambda_{N, {\rm off}}= N^{-\alpha}\lambda_{{\rm off}}$ and $\lambda_{N, {\rm on}}\equiv \lambda_{{\rm on}}$ for some parameters $\lambda_{\rm on},\lambda_{\rm off},\alpha>0$, we obtain that the process {$a_N(\cdot)$} satisfies {Assumptions} \ref{ass:ass1} and \ref{ass:ass2} with $\varrho_N = N^{-\alpha}$,      \begin{align*}
        \kappa_{\star}(s,t) = \frac{\lambda_{\rm off}}{\lambda_{\rm on}}e^{-\lambda_{\rm on}\abs{s-t}},
    \end{align*}
    and $p^{\star}(t) \equiv \lambda_{\rm off}/\lambda_{\rm on}$. \hfill$\clubsuit$
\end{example}

\section{Graph notations}\label{section:graph_notations}

In this section we introduce various graph-related notations. {Due to the dynamic nature of our graph process, we need to introduce some new concepts that were not required in the static setting.}
\begin{itemize}
    \item[$\circ$] {Throughout the paper we consider only non-empty graphs whose vertex set is a subset of $\N$. We denote the set of such graphs by ${\tt G}_{\infty}$. For any graph $H\in{\tt G}_{\infty}$ we denote by ${\tt V}(H)$ (${\tt E}(H)$, respectively) the set of vertices (edges, respectively) of the graph $H$. By definition of the set ${\tt G}_{\infty}$, ${\tt V}(H)\subset \N$ and ${\tt E}(H)\subset \{(u,v)\in\N^2\colon u<v\}$. In particular, for any $t\in[0,T]$ the random graph $G_N(t)$ can be considered as a random element of ${\tt G}_{\infty}$ with $N$ vertices.}
    \item[{$\circ$}] {For} two graphs {$H,{H^{\star}}\in{\tt G}_{\infty}$} we write ${H}\sim {H^{\star}}$ if they are {isomorphic}, and we write ${H}\subset {H^{\star}}$ if $H$ is a subgraph of ${H^{\star}}${, i.e., ${\tt V}(H)\subseteq{\tt V}(H^{\star})$ and ${\tt E}(H)\subseteq{\tt E}(H^{\star})$}.  Denote {for any $N\in\N$} by $K_N$ a complete graph with vertices labeled $1,\ldots,N${, i.e., ${\tt V}(K_N) = \sprod{N}$ and ${\tt E}(K_N) = \{(u,v)\in\N^2\colon u<v\}$.} 
    \item[$\circ$]
{For any graph {$H\in{\tt G}_{\infty}$} we denote by $\mathcal{V}(H)$ ($\mathcal{E}(H)$, respectively) the number of vertices (edges, respectively) in $H${, i.e., $\mathcal{V}(H) = \abs{{\tt V}(H)}$ and $\mathcal{E}(H) = \abs{{\tt E}(H)}$}; we also use the notation $\mathcal{A}(H)$ to denote the number of automorphisms of $H$.}    
 \item[$\circ$]
 {For any pair  $g_1,g_2 \in {\tt G}_{\infty}$ we use the notation $g_1\cup g_2$ for the graph with ${\tt V}(g_1\cup g_2) = {\tt V}(g_1)\cup{\tt V}(g_1)$ and ${\tt E}(g_1\cup g_2) = {\tt E}(g_1)\cup {\tt E}(g_2)$. The graphs $g_1\cap g_2$ and $g_1\setminus g_2$ can be defined in a similar way.}
\item[$\circ$]
The main object that is studied in this paper, concerns the $m$-dimensional subgraph count process $(X_{N,1}(\cdot),\ldots,X_{N,m}(\cdot))$, whose $i$-th entry represents the number of subgraphs of $G_N(\cdot)$ which are isomorphic to some fixed graphs {$H_i\in{\tt G_{\infty}}$}. Given any sequence of subgraphs $H_1,\ldots,H_m{\in{\tt G}_{\infty}}$, {we can} interpret $X_{N,i}(t)$ for $i\in\sprod{m}$ as the number of {graphs} $g\, {\in{\tt G}_{\infty}}$ such that $g\sim {H}_i$ and $g{\,{\subset}\,} G_N(t)$. {Using that $G_N(t)\subset K_N$  for any fixed $N\in N$ almost surely, we can consider only such graphs $g\,\in{\tt G}_{\infty}$ that ${\tt V}(g)\subseteq\sprod{N}$.} As such, it is convenient {for any $H\in{\tt G}_{\infty}$} to define the set of all possible {graphs} $g\,{\in{\tt G}_{\infty
}}$ {with ${\tt V}(g)\subseteq\sprod{N}$} such that they are isomorphic to some fixed graph $H$: 
\begin{align*}
    {\tt G}_N({H}) = {\{g\in{\tt G}_{\infty}\colon {\tt V}(g)\subseteq\sprod{N},\, g\sim {H}\}}.
\end{align*}
The cardinality of the set {${\tt G}_N(H)$ is}
\begin{align}
    \abs{{\tt G}_N({H})} = \frac{N!}{(N-\mathcal{V}({H}))!\,\mathcal{A}({H})}.\label{G_cardinality}
\end{align}
{For example, if $H$ is a triangle, then $\abs{{\tt G}_N(H)} = \binom{N}{3}$.} Define also for any ${g\,{\in{\tt G}_{\infty}}}$ the indicator function of this {graph} $g$ being a subgraph {of $G_N(t)$:}
\begin{align}
    \mathbb{I}_N(g,t) = \mathbb{I}\{g{\,{\subset}\,} G_N(t)\} = \prod_{(u,v)\in{\tt E}(g)}a_{N,u,v}(t).\label{I_N_def}
\end{align}

\item[$\circ$]
{Throughout the paper we need to control the way how different subgraph patterns can be obtained when the graphs are intersecting. To formalize this concept, denote the quotient set of all graphs with at least one edge by ${\tt G}_{\sim} = \{g\in{\tt G}_{\infty}\colon \mathcal{E}(g)\geqslant 1\}/\sim$, where $\sim$ was introduced above. We refer to any element of ${\tt G}_{\sim}$ as a \emph{subgraph pattern}. We denote by $\pi:{\tt G}_{\infty}\to{\tt G}_{\sim}$ the classical projection on the quotient set. Then, define for any pair of graphs $H,\,{H}^{\star}\in{\tt G}_{\infty}$ a set of all common subgraph patterns with at least one edge:
\begin{align*}
    {\tt CS}({H,H^{\star}}) = \{g\in{\tt G}_{\sim}\colon \exists g_1,g_2\in{\tt G}_{\infty} \colon g_1\subset {H},\,g_2\subset {H^{\star}},\, \pi(g_1)=\pi(g_2)= g\}.
\end{align*}
A graph pattern $g\in {\tt CS}({H,H^{\star}})$ represents an unlabeled subgraph structure common to both $H$ and $H^\star$.} 

{
    For example, let $H$ be a triangle with vertex set $\{1,2,3\}$ and $H^{\star}$  an edge with vertices $\{4,5\}$. Then, even though graphs $H$ and $H^{\star}$ do not have any common subgraph, they still have a common subgraph pattern, which is the edge pattern.
}

For brevity we write ${\tt CS}({H}) = {\tt CS}({H},{H})$, which basically stands for the set of all {subgraph patterns} of graph $H$. 
{Note that functions $\mathcal{V}(\cdot)$, $\mathcal{E}(\cdot)$ and $\mathcal{A}(\cdot)$ are invariant with respect to graph isomorphism, so we can as well define them for the subgraph patterns.}

{In addition, for any graph $H\in{\tt G}_{\infty}$ and any graph pattern $g\in{\tt G}_{\sim}$ we define by $\mathcal{S}(H,g)$ the total number of different subgraphs of $H$ which are isomorphic to $g$, i.e.,
   \begin{align}
        \mathcal{S}(H,g) = \sum_{g_0\,\subset\, {H}}\mathbb{I}\{\pi(g_0)= g\}.\label{S_cal_def}
    \end{align}
}
\item[{$\circ$}]
For any graph $H{\in{\tt G}_{\infty}}$ (or graph pattern {$H\in{\tt G}_{\sim}$)} with at least one edge we define {for any $N\in\N$}
%\KB{RSH: check in the bracket if this is $H\in{\tt G}_{\sim}$}
\begin{align}
    \mathcal{F}_N({H}) = N^{-\mathcal{V}({H})}\varrho_N^{-\mathcal{E}({H})}.\label{F_claim_0}
\end{align}

{In this paper we mostly consider such graphs $H$ that satisfy the following assumption.
\begin{ass}\label{ass:ass3}
    {For any $g\subset H$ such that $\mathcal{V}(g)\geqslant 1$,
    \begin{align}
        \lim_{N\to\infty}{\mathcal{F}_N(g)} = 0.\label{graph_assumption}
    \end{align}}
\end{ass}
\COM{Equivalently, Assumption \ref{ass:ass3} tells us that for any $g\subset H$ such that $\mathcal{V}(g)\geqslant 1$
\begin{align*}
    \lim_{N\to\infty}\mathcal{F}_N(g) = 0.
\end{align*}}

\begin{remark}\label{ass3_remark}{\em
{Consider the following scenario where Assumption \ref{ass:ass3} is  not satisfied: suppose that for some subgraph $g \subset H$,
\begin{align}
    \lim_{N\to\infty} N^{-\mathcal{V}(g)}\varrho_N^{-\mathcal{E}(g)} = \infty. \label{non_ass3}
\end{align}
{Denote by $X_N(t)$ the number of subgraphs of $G_N(t)$ which are isomorphic to $H$, and $X_N^\star(t)$ the number of subgraphs of $G_N(t)$ which are isomorphic to $g$.}

Then,
\begin{align*}
    \pk{X_N(t) \geqslant 1} \leqslant \pk{X_N^\star(t) \geqslant 1} \leqslant \E{X_N^\star} = \frac{N! \, p_N^{\mathcal{E}(g)}(t)}{(N-\mathcal{V}(g))! \mathcal{A}(g)} \leqslant \frac{N^{\mathcal{V}(g)} p_N^{\mathcal{E}(g)}(t)}{\mathcal{A}(g)}.
\end{align*}
From \eqref{non_ass3}, it follows that
\begin{align*}
    \lim_{N\to\infty}\pk{X_N(t) = 0} \geqslant 1 - \lim_{N\to\infty} \frac{1}{\mathcal{A}(g) N^{-\mathcal{V}(g)} p_N^{-\mathcal{E}(g)}(t)} = 1,
\end{align*}
which implies that, if \eqref{non_ass3} holds, it is impossible to obtain a non-degenerate limit of the sequence $X_N(t)$.

In the borderline case where, for some subgraph $g \subset H$,
\begin{align*}
    \lim_{N\to\infty} N^{-\mathcal{V}(g)}\varrho_N^{-\mathcal{E}(g)} = C \in (0, \infty),
\end{align*}
it is known from \cite{rucinski1988small} that even in the static case, the limiting subgraph count is non-Gaussian. {We exclude this case from our paper as it requires a radically different approach}.

The conclusion of the above considerations is that Assumption \ref{ass:ass3} is strictly necessary to ensure a non-degenerate subgraph count as $N \to \infty$.
}
\hfill
${\Diamond}$}

\end{remark}

The asymptotic behavior of functions $\mathcal{F}_N(g)$ may be different for different subgraphs $g$. It is therefore useful to identify, {for any pair of graphs $H,H^{\star}\in{\tt G}_{\infty}$, common  subgraph patterns} for which the decay of $\mathcal{F}_N(g)$ is the slowest. {More, concretely, we consider {common subgraph patterns} $g$ for which $\mathcal{F}_N(g)$ 
remains asymptotically at least as large as the decay $\mathcal{F}_N(g^\star)$ of any other {common subgraph pattern} $g^\star$}:
\begin{align}
    {{\tt OCS}({H},{H^{\star}}) = \left\{g\in{\tt CS}({H},{H^{\star})}\colon \forall g^{\star}\in{\tt CS}({H},{H^{\star}})  :\lim_{N\to\infty}\frac{\mathcal{F}_N(g^{\star})}{\mathcal{F}_N(g)}\in[0,\infty)
    \right\}.\label{OCS_def}}
\end{align}
In particular, {as \eqref{rho_choice} tells us that for any $g,g^{\star}\in {\tt CS}(H,H^{\star})$}
\begin{align*}
    {\lim_{N\to\infty}\frac{\mathcal{F}_N(g^{\star})}{\mathcal{F}_N(g)} = \lim_{N\to\infty}N^{\mathcal{V}(g) - \mathcal{V}(g^{\star})}\varrho_N^{\mathcal{E}(g) - \mathcal{E}(g^{\star})}\in\{0,1,\infty\} ,}   
\end{align*}
{due to \eqref{OCS_def}} for any $g\in{\tt OCS}({H},{H^{\star}})$ and any $g^{\star}\in{\tt CS}({H},{H^{\star}})$,
\begin{align}
    \lim_{N\to\infty}\frac{\mathcal{F}_N(g^{\star})}{\mathcal{F}_N(g)} = \begin{cases} 1,\quad \text{iff} \quad g^{\star}\in{\tt OCS}({H},{H^{\star}}),\\
        0,\quad \text{iff} \quad g^{\star}\not\in{\tt OCS}({H},{H^{\star}}).
    \end{cases}\label{opt_property}
\end{align}
Hence, for any graphs {${H},{H^{\star}}\in{\tt G}_{\infty}$}, {if we define $\mathcal{F}^{\rm opt}_N({H},{H^{\star}}) = \mathcal{F}_N(g^{\star})$ {for some \mbox{$g^{\star}\in{\tt OCS}(H,H^{\star})$}},} {then due to \eqref{opt_property} for any} {$g\in{\tt CS}({H},{H^{\star}})$,}
\begin{align}
    \lim_{N\to\infty}\frac{\mathcal{F}_N(g)}{\mathcal{F}^{\rm opt}_N({H},{H^{\star}})} = \mathcal{F}_{\infty}(g,{H},{H^{\star}}) =\begin{cases} 1,\quad \text{iff} \quad g\in{\tt OCS}({H},{H^{\star}}),\\
        0,\quad \text{iff} \quad g\not\in{\tt OCS}({H},{H^{\star}}).\label{F_cal_opt_def}
    \end{cases}
\end{align}
Again, for simplicity we write $\mathcal{F}^{\rm opt}_N({H}):=\mathcal{F}^{\rm opt}_N({H},{H})$. 
We call two graphs {\mbox{$H,\,H^{\star}\in{\tt G}_{\infty}$}} {\it equioptimal} if
\begin{align} 
    \lim_{N\to\infty} \frac{\mathcal{F}_N^{\rm opt} ({H}, {H^{\star}})}{\mathcal{F}_N^{\rm opt} ({H})}=\lim_{N\to\infty} \frac{\mathcal{F}_N^{\rm opt} ({H}, {H^{\star}})}{\mathcal{F}_N^{\rm opt} ({H^{\star}} )}=1.\label{equioptimal}
\end{align} 

\begin{remark} {\em
    Again, by virtue to \eqref{rho_choice}, the limit in \eqref{equioptimal} cannot be positive, finite and different from $1$.}
    \hfill ${\Diamond}$
\end{remark} 

\item[$\circ$]
We write ${H}\cong {H^{\star}}$ if graphs {$H,\,{H^{\star}}\in{\tt G}_{\infty}$} are equioptimal. A combination of \eqref{equioptimal} and \eqref{F_cal_opt_def} tells us that if ${H}\cong {H^{\star}}$ then
\begin{align}
    {\tt OCS}({H},{H^{\star}}) = {\tt OCS}({H})\cap {\tt OCS}({H^{\star}}),\label{GCS_cong}
\end{align}
and if ${H}\not\cong {H^{\star}}$ then
\begin{align}
    {\tt OCS}({H})\cap {\tt OCS}({H^{\star}})=\varnothing.\label{GCS_not_cong}
\end{align}

\item[$\circ$]
Finally, for any graphs ${H},{H^{\star}}{\in{\tt G}_{\infty}}$ and $g\in{\tt CS}({H},{H^{\star}})$,  denoting ${\mathcal{N}(H,H^{\star},g)}=\mathcal{V}({H})+\mathcal{V}({H^{\star}}) - \mathcal{V}(g)$, we define
\begin{align}
    \mathcal{C}({H},{H^{\star}},g):=\sum_{\substack{g_1\in{\tt G}_{{\mathcal{N}(H,H^{\star},g)}}({H})\\ g_2\in{\tt G}_{{\mathcal{N}(H,H^{\star},g)}}({H^{\star}}) \\ {\pi(g_1\cap g_2)= g}}}\frac{1}{{\mathcal{N}(H,H^{\star},g)}!}{.}\label{C_def}
\end{align}
Then, for any $N\geqslant {\mathcal{N}(H,H^{\star},g)}$ we can see that the number of {pairs of graphs $g_1,g_2\in{\tt G}_{\infty}$ such that $g_1\sim H,\,g_2\sim {H^{\star}}$, ${\tt V}(g_1),{\tt V}(g_2)\subset \sprod{N}$ and $\pi(g_1\cap g_2)= g$} equals
\begin{align*}
    \sum_{\substack{g_1\in{\tt G}_{N}({H})\\ g_2\in{\tt G}_{N}({H^{\star}}) \\ {\pi(g_1\cap g_2)\sim g}}}1 = \frac{N!}{(N-{\mathcal{N}(H,H^{\star},g)})!}\mathcal{C}({H},{H}^{\star},g).
\end{align*}
The following lemma shows a different representation of the constant $\mathcal{C}$ defined in \eqref{C_def} in case of equioptimal graphs under {\eqref{graph_assumption}}.
\begin{lem}\label{lem:constant_formula} 
    Let $H,\,H^{\star}{\in{\tt G}_{\infty}}$ be two graphs which satisfy {Assumption \ref{ass:ass3}}. If $H\cong H^{\star}$ and $g\in {\tt OCS}(H, H^{\star})$, we have the following representation for the constant $\mathcal{C}(H,H^{\star},g)$:
    \begin{align*} 
        \mathcal{C}(H, H^{\star}, g) =\left(\frac{\sqrt{\mathcal{A}(g)}\,\mathcal{S}(H,g)}{\mathcal{A}(H)}\right)\left(\frac{\sqrt{\mathcal{A}(g)}\,\mathcal{S}(H^{\star},g)}{\mathcal{A}(H^{\star})}\right),
    \end{align*}
    where {the constants $\mathcal{S}(H,g)$ and $\mathcal{S}(H^{\star},g)$ are defined in \eqref{S_cal_def}.}
\end{lem}

The proof of this lemma is provided in the appendix.}
\end{itemize}

\section{Main result}\label{section:main_results} 

In this section we state our functional central limit theorem for a vector of  subgraph count processes. {To this end,} we fix, for a given $m\in{\mathbb N}$,  (deterministic) graphs ${H}_1,\ldots,{H}_m{\in{\tt G}_{\infty}}$, which satisfy Assumption \ref{ass:ass3}; think for instance of $m=2$ and ${H}_1$ corresponding to edges and ${H}_2$ to triangles. 
Recall that from \eqref{graph_assumption} we know that for any $i\in\sprod{m}$ and any $g\subset {H}_i$ with at least one vertex
\begin{align}
    \lim_{N\to\infty}\mathcal{F}_N(g) = 0.\label{F_cal_zero_limit}
\end{align}

The key object in our analysis is the $m$-dimensional process
\begin{align*}
    \bigl(\vk X_N(t)\bigr)_{t\in[0,T]} = \big((X_{N,1}(t)),\ldots, (X_{N,m}(t))\big)_{{t\in[0,T]}}^{\top},
\end{align*}
where, for $i\in \sprod{m}$, $X_{N,i}(t)$ records the number of {subgraphs of $G_N(t)$ which are isomorphic to ${H}_i$}. {Using \eqref{I_N_def} the marginals of the process {${\boldsymbol X}_N(t)$} can be alternatively defined via
\begin{align}
    X_{N,i}(t)=\sum_{g\in{\tt G}_N({H}_i)}\mathbb{I}\{g{\,{\subset}\,} G_N(t)\}=\sum_{g\in{\tt G}_N({H}_i)}\prod_{(u,v)\in{\tt E}(g)}a_{N,u,v}(t).\label{X_def_with_g}
\end{align}}
% The aim of this paper is to find a limiting process 
% \begin{align*} 
%     \vk X(t)=\lim_{N\to\infty}\frac{\vk X_N(t)-\E{\vk X_N(t)}}{\sqrt{\Var\bigl(\vk X_N(t)\bigr)}}.
% \end{align*}

For each {graph pattern $g\in{\tt G}_{\sim}$} we define a centered Gaussian process $(X^{\star}_{g}(t))_{t\in[0,T]}$, characterized via the covariance function
\begin{align*}
    \operatorname{Cov}(X^{\star}_{g}(t),X^{\star}_g(s)) = \big(\kappa_{\star}(s,t)\big)^{\mathcal{E}(g)};
\end{align*}
in addition, all processes $(X^{\star}_{g}(t))_{t\in[0,T]}$ are independent (across the {graph patterns} $g$). 
We can now state our main result. 
% It uses some notation that will be introduced later: the functions $\mathcal{F}_{N}^{\rm opt}(\cdot)$, $\mathcal{S}(\cdot,\cdot)$ and the set $\mathfrak{OCS}(\cdot,\cdot)$ are discussed in Section \ref{section:graph_notations}.

\begin{theo}\label{main}
    {Let $G_N(t)$ for $t\in[0,T]$ be a dynamic {Erd\H{o}s--R\'{e}nyi} random graph with $N$ vertices which satisfies Assumptions \ref{ass:ass1} and \ref{ass:ass2}, and let $H_1,\ldots, H_{m}{\in{\tt G}_{\infty}}$ be {graphs} which satisfy Assumption \ref{ass:ass3}. Then, as $N\to\infty$, the process $\vk X_N(\cdot)$ defined in \eqref{X_def_with_g} satisfies the following limiting theorem:}
    \begin{align*}
        \left(\frac{X_{N,1}(t) - \E{X_{N,1}(t)}}{\sqrt{\mathcal{F}^{\rm opt}_N({H}_1)}/\mathcal{F}_N({H}_1)},\ldots,\frac{X_{N,m}(t) - \E{X_{N,m}(t)}}{\sqrt{\mathcal{F}^{\rm opt}_N({H}_m)}/\mathcal{F}_N({H}_m)}\right)^{\top}\to \vk X(t),
    \end{align*}
    in distribution in $\mathbb{D}([0,T],\R^{m})$, where {$\vk X(\cdot)$} is centered vector-valued Gaussian process with components
    \begin{align*}
        X_i(t) = \sum_{g\in{\tt OCS}({H}_i)}\frac{\sqrt{\mathcal{A}(g)}\,\mathcal{S}({H}_i,g)}{\mathcal{A}({H}_i)}\bigl(p_{\star}(t)\bigr)^{\mathcal{E}({H}_i) - \mathcal{E}(g)}X^{\star}_{g}(t).
    \end{align*}
\end{theo}

{

\begin{remark}{\em
    The covariance matrix function $\Sigma(s,t)$ of the process $\vk X(\cdot)$ defined in \netheo{main} can be evaluated more explicitly. A standard computation reveals that it can be represented, in terms of the function $p_\star(\cdot)$, as 
    \begin{align*}
        \Sigma_{i,j}(s,t) = \mathbb{I}\{{H}_i\cong {H}_j\}\sum_{g\in{\tt OCS}({H}_i,{H}_j)}\frac{\mathcal{A}(g)\,\mathcal{S}({H}_i,g)\,\mathcal{S}({H}_j,g)}{\mathcal{A}({H}_i)\,\mathcal{A}({H}_j)}\frac{\bigl(p_{\star}(s)\bigr)^{\mathcal{E}({H}_i)}\bigl(p_{\star}(t)\bigr)^{\mathcal{E}({H}_j)}}{\bigl(p_{\star}(s)\,p_{\star}(t)\bigr)^{\mathcal{E}(g)}}\kappa^{\mathcal{E}(g)}_{\star}(s,t)
    \end{align*}
    for $s,t\geqslant 0$.}
   \hfill ${\Diamond}$
\end{remark}

{The following corollary covers the case when the distribution of the edge processes of $G_N({\cdot})$ is independent of the size of the graph.}

\begin{korr}\label{main_corr}
     Let $G_N(t)$ for $t\in[0,T]$ be a dynamic {Erd\H{o}s--R\'{e}nyi} random graph with $N$ vertices, {which satisfies  Assumption \ref{ass:ass2} (for $\varrho_N = 1$), and assume that the edge process $a_{N}(\cdot)$ does not depend on $N$.} {Let  $H_1,\ldots, H_{m}\in{\tt G}_{\infty}$} be graphs with at least one edge. Then, the process $\vk X_N(\cdot)$ defined in \eqref{X_def_with_g} satisfies the following limiting theorem: as $N\to\infty$,
    \begin{align*}
        \left(\frac{X_{N,1}(t) - \E{X_{N,1}(t)}}{N^{\mathcal{V}(H_1)-1}},\ldots,\frac{X_{N,m}(t) - \E{X_{N,m}(t)}}{N^{\mathcal{V}(H_m)-1}}\right)^{\top}\to \vk X(t)
    \end{align*}
    in distribution in $\mathbb{D}([0,T],\R^{m})$, where {$\vk X(\cdot)$} is a centered vector-valued Gaussian process with components
    \begin{align*}
        X_i(t) =\frac{\sqrt{2}\,\mathcal{E}(H_i)}{\mathcal{A}({H}_i)}\bigl(p(t)\bigr)^{\mathcal{E}({H}_i) - 1}X^{\star}(t),
    \end{align*}
    and where $p(t) = p_N(t)$, $\kappa(s,t) = \kappa_N(s,t)$ (which are independent of $N$ due to our assumption), and {$X^{\star}(\cdot)$} is a centered Gaussian process characterized via its covariance function
    \begin{align*}
        \operatorname{Cov}(X^{\star}(s),X^{\star}(t)) = \kappa(s,t).
    \end{align*}
\end{korr}

\begin{example}\label{ex:subgraphs} \em To get an impression of a typical application of Theorem \ref{main}, consider the case of $m=2$, with the subgraphs corresponding to wedges and triangles.
Consider the edge processes defined in Example~\ref{ex:process} with $\alpha\in(0,1)$. Then one can check that both wedge graph (denoted by ${H}_1$) and triangle graph (denoted by ${H}_2$) satisfy \eqref{graph_assumption}, and
\begin{align*}
    &\mathcal{F}_N({H}_1) = N^{2\alpha - 3},\qquad \mathcal{F}_N({H}_2) = N^{3\alpha - 3},\\
    &\mathcal{F}^{\rm opt}_N({H}_1) = N^{\alpha - 2},\qquad \mathcal{F}^{\rm opt}_N({H}_2) = N^{\max(\alpha - 2, 3\alpha - 3)}.
\end{align*}
{The next step is to describe the sets ${\tt OCS(H_1)}$ and ${\tt OCS(H_2)}$. The first one is easy to find: ${\tt OCS(H_1)} = \{g_1\}$, where ${g_1\in {\tt G}_{\sim}}$ is {the pattern of} an edge. However, for the second one we need to distinguish three cases:}
\begin{itemize}
    \item[$\circ$]{If $\alpha\in(0,1/2)$ then ${\tt OCS}({H}_2) = \{g_1\}$},
    \item[$\circ$]{If $\alpha=1/2$ then ${\tt OCS}({H}_2) = \{g_1,g_2\}$},
    \item[$\circ$]{If $\alpha\in(1/2,1)$ then ${\tt OCS}({H}_2) = \{g_2\}$},
\end{itemize}
where ${g_2\in{\tt G}_{\sim}}$ is {the pattern of} a triangle. It is readily verified  that
\begin{align*}
    \mathcal{A}({H}_1) = 2,\qquad \mathcal{A}({H}_2) = 6,\qquad \mathcal{A}(g_1) = 2,\qquad \mathcal{A}(g_2) = 6,\\
    \mathcal{S}({H}_1,g_1) = 2,\qquad \mathcal{S}({H}_2,g_1) = 3,\qquad \mathcal{S}({H}_2,g_2) = 1.
\end{align*}
Applying \netheo{main} we obtain that, as $N\to\infty$,
\begin{align*}
    \left(\frac{X_{N,1}(t) - \E{X_{N,1}(t)}}{N^{ 2 - 3\alpha/2}} , \frac{X_{N,2}(t) - \E{X_{N,2}(t)}}{N^{\max(2-5\alpha/2,\, 3/2 - 3\alpha/2)}}\right)^{\top}\to (X_1(t),X_2(t))^{\top},
\end{align*}
where
\begin{align*}
    X_1(t) &= \frac{2\sqrt{2}}{2}\frac{\lambda_{\rm off}}{\lambda_{\rm on}}X_1^{\star}(t),\\
    X_2(t) &= \mathbb{I}\{\alpha\leqslant 1/2\}\frac{3\sqrt{2}}{6}\left(\frac{\lambda_{\rm off}}{\lambda_{\rm on}}\right)^2X_1^{\star}(t) + \mathbb{I}\{\alpha\geqslant 1/2\}\frac{\sqrt{6}}{6}X^{\star}_2(t);
\end{align*}
the coefficients in these expressions can be simplified in the obvious manner, and $X_1^{\star}(\cdot)$ and $X_2^{\star}(\cdot)$ are two independent centered Gaussian processes, characterized via the covariance functions
\begin{align*}
    \operatorname{Cov}(X_1^{\star}(s),X_{1}^{\star}(t)) &= \frac{\lambda_{\rm off}}{\lambda_{\rm on}}e^{-\lambda_{\rm on}\abs{s-t}},\\
    \operatorname{Cov}(X_2^{\star}(s),X_{2}^{\star}(t)) &= \left(\frac{\lambda_{\rm off}}{\lambda_{\rm on}}\right)^{3}e^{-3\lambda_{\rm on}\abs{s-t}},
\end{align*}
respectively.
In particular, for $\alpha<1/2$ we have that processes {$X_1(\cdot)$} and {$X_2(\cdot)$} are perfectly correlated, whereas for $\alpha>1/2$ they are independent. {In case $\alpha=1/2$, the processes {$X_1(\cdot)$} and {$X_2(\cdot)$} are positively but not perfectly correlated.} \hfill$\clubsuit$
\end{example}

}

\section{Proof of \netheo{main}}\label{section:proof_of_main}
The main objective of this section is to establish \netheo{main}. The following objects play a crucial role: for $N\in\N$ and $i\in\sprod{m}$, we define the processes
\begin{align}
    X^{\star}_{N,i}(t) = \frac{X_{N,i}(t) - \E{X_{N,i}(t)}}{\sqrt{\mathcal{F}^{\rm opt}_N({H}_i)}/\mathcal{F}({H}_i)}.\label{X_star_def}
\end{align}
The proof of \netheo{main} is split up in a number of steps. In the first place, using that the stochastic processes {$\vk X_N^{\star}(\cdot)$} are centered, we show the convergence of its covariance function to the covariance function of the process {$\vk X(\cdot)$} featuring in \netheo{main}. We state this result in the following lemma.  

\begin{lem}\label{lem:exp_and_cov}
    For any $N\in\N$, $i\in\sprod{m}$ and $t\in[0,T]$,
    \begin{align}
        \E{X_{N,i}(t)} = \frac{N!\,p_N^{\mathcal{E}({H}_i)}(t)}{(N-\mathcal{V}({H}_i))!\mathcal{A}({H}_i)}{,}\label{E_formula}
    \end{align}
    and for any $t,s\in[0,T]$, $i,j\in\sprod{m}$ 
    \begin{align}
        &\lim_{N\to\infty}\frac{\mathcal{F}_N({H}_i)\mathcal{F}_N({H}_j)}{\mathcal{F}_N^{\rm opt}({H}_i,{H}_j)}\operatorname{Cov}\bigl(X_{N,i}(t),\,X_{N,j}(s)\bigr) \notag\\
        &\qquad\qquad\qquad= \bigl(p_{\star}(s)\bigr)^{\mathcal{E}({H}_i)}\bigl(p_{\star}(t)\bigr)^{\mathcal{E}({H}_j)} \sum_{g\in{\tt OCS}({H}_i,{H}_j)}\mathcal{C}({H}_i,{H}_j,g)\left(\frac{\kappa_{\star}(s,t)}{p_{\star}(t)\,p_{\star}(s)}\right)^{\mathcal{E}(g)}.\label{Cov_formula}
    \end{align}
\end{lem}

More precisely, we use the following corollary of this lemma, due to \eqref{equioptimal}.

\begin{korr}\label{korr:exp_and_cov}
    For any $i,j\in\sprod{m}$ and any $t,s\in[0,T]$
    \begin{align*}
        \lim_{N\to\infty}\operatorname{Cov}(X_{N,i}^{\star}(t),&X_{N,j}^{\star}(s))\\&=\mathbb{I}\{{H}_i\cong {H}_j\}\bigl(p_{\star}(s)\bigr)^{\mathcal{E}({H}_i)}\bigl(p_{\star}(t)\bigr)^{\mathcal{E}({H}_j)} \sum_{g\in{\tt OCS}({H}_i,{H}_j)}\mathcal{C}({H}_i,{H}_j,g)\left(\frac{\kappa_{\star}(s,t)}{p_{\star}(t)\,p_{\star}(s)}\right)^{\mathcal{E}(g)}.
    \end{align*}
\end{korr}

We postpone the proof of \nelem{lem:exp_and_cov} to Section \ref{section:proof_of_covariance}. 
To link the results of \nelem{lem:exp_and_cov} with the covariance function of the process {$\vk X(\cdot)$}, we rely on the following proposition, proven in the appendix.
\begin{prop}\label{prop:X_vk_cov}
    The covariance matrix function of the process $\vk X(t)$ satisfies the following representation for any $i,j\in\sprod{m}$ and any $t,s\in[0,T]$:
    \begin{align*}
        \operatorname{Cov}(X_i(t),X_j(s)) = \mathbb{I}\{{H}_i\cong {H}_j\}\bigl(p_{\star}(s)\bigr)^{\mathcal{E}({H}_i)}\bigl(p_{\star}(t)\bigr)^{\mathcal{E}({H}_j)} \sum_{g\in{\tt OCS}({H}_i,{H}_j)}\mathcal{C}({H}_i,{H}_j,g)\left(\frac{\kappa_{\star}(s,t)}{p_{\star}(t)\,p_{\star}(s)}\right)^{\mathcal{E}(g)}.
    \end{align*}
\end{prop}

The next step of proving \netheo{main} is to show the convergence of finite-dimensional distributions, which can be formulated as follows.

\begin{lem}\label{lem:asymptotic_normality}
    For any $n\in\N$ and $t_1,\ldots,t_n\in[0,T]$ random vector
    \begin{align*}
        \bigl(\vk X^{\star}_N(t_1),\ldots, \vk X^{\star}_N(t_n)\bigr)
    \end{align*}
    converges in distributions to a centered Gaussian random vector with covariance block matrix $\vk \Sigma$ which elements are given via
    \begin{align*}
        \bigl(\vk \Sigma_{k,l}\bigr)_{i,j} =\mathbb{I}\{{H}_i\cong {H}_j\} \bigl(p_{\star}(t_{k})\bigr)^{\mathcal{E}({H}_i)}\bigl(p_{\star}(t_{l})\bigr)^{\mathcal{E}({H}_j)}\sum_{g\in{\tt OCS}({H}_i,{H}_j)}\mathcal{C}({H}_i,{H}_j,g)\left(\frac{\kappa_{\star}(t_k,t_l)}{p_{\star}(t_k)\,p_{\star}(t_l)}\right)^{\mathcal{E}(g)}
    \end{align*}
    for $k,l\in\sprod{n}$ and $i,j\in \sprod{m}$.
\end{lem}
Note that each element of the block matrix $\vk \Sigma$ is a matrix itself, as it is the covariance of two Gaussian vectors.
Again, the proof of this lemma is postponed to a separate section, namely Section~\ref{section:proof_of_normality}.

Finally, to justify the claim of \netheo{main}, we apply the analogue of \cite[Theorem 13.5]{billingsley2013convergence} for vector valued process, in particular, the condition given in \cite[Eqn.\ (13.14)]{billingsley2013convergence}. According to this condition, {to show that the sequence {$\vk X_N(\cdot)$} is tight in $\mathbb{D}([0,T],\R^{m})$}, it is sufficient to find a continuous function $F(t)$ {for $t\in[0,T]$} such that for any {$0\leqslant r\leqslant s\leqslant t\leqslant T$,} 
\begin{align*}
    \E{\abs{\vk X^{\star}_N(t)-\vk X^{\star}_N(s)}^2\abs{\vk X^{\star}_N(s)-\vk  X^{\star}_N(r)}^2}\leqslant \abs{F(t)-F(r)}^2.
\end{align*}

The following lemma provides us with the required function.
    
\begin{lem}\label{lem:tightness} For any $0\leqslant r<s<t\leqslant T$, we have
    \begin{align}
        \E{\abs{\vk X^{\star}_N(t)-\vk X^{\star}_N(s)}^2\abs{\vk X^{\star}_N(s)-\vk 
        X^{\star}_N(r)}^2}\leqslant \abs{F(t)-F(r)}^2{,}\label{tightness_claim}
    \end{align}

    where 
    \begin{align*}
        F(t) = t\sum_{i,j=1}^{m}\biggl(\sqrt{51\,F_{i,j}^{\star}C_{i,j}^{\star}}\biggr),
    \end{align*}

    and for any $i,j\in\sprod{m}$
    \begin{align}
        F_{i,j}^{\star} &=4{\sum_{V=\max(\mathcal{V}(H_i),\mathcal{V}(H_j))}^{2\mathcal{V}(H_i) + 2\mathcal{V}(H_j) - 2}}\left(\frac{V!}{(V-\mathcal{V}({H}_i))!\mathcal{A}({H}_i)}\right)^2\left(\frac{V!}{(V-\mathcal{V}({H}_j))!\mathcal{A}({H}_j)}\right)^2,\label{frak_f_def}\\
        C_{i,j}^{\star} &= 16\,C_0\max(1,T^2)\,{\mathfrak{C}^4}\,\mathcal{E}^2({H}_i)\,\mathcal{E}^2({H}_j)\,p_+^{2\mathcal{E}({H}_i)+2\mathcal{E}({H}_j)-2},\label{frak_c_def}\\
        C_0 &= 2\,{\mathfrak{C}^2} \left(\max_{i\in\sprod{m}}\mathcal{E}({H}_i)\right)^2\,{p_+^{4\max_{i\in\sprod{m}}\mathcal{E}({H}_i)}}.\label{C_0_def}
    \end{align}
\end{lem}
The proof of this lemma being provided in Section \ref{section:proof_of_tightness}. {We} obtain the claimed convergence by combining \cite[Theorem 13.5]{billingsley2013convergence} with {Lemmas~\ref {lem:asymptotic_normality} and \ref{lem:tightness} and \neprop{prop:X_vk_cov}.}

\section{Proof of \nelem{lem:exp_and_cov}}\label{section:proof_of_covariance}

We start by establishing the claim \eqref{E_formula}. Taking expectation of the both sides of \eqref{X_def_with_g}, and using that all $a_{N,v,u}(t)$ are independent for all $(v,u)\in{\tt E}(g)$ and all $g\in{\tt G}_N({H_i})$, we obtain
\begin{align*}
    \E{X_{N,i}(t)} &= \sum_{g\in{\tt G}_{N}({H_i})}\E{\mathbb{I}\{g{\,{\subset}\,} G_N(t)\}}\notag\\
               &= \sum_{g\in{\tt G}_{N}({H_i})}\pk{\prod_{(u,v)\in{\tt E}(g)}a_{{N,u,v}}(t)=1}\notag= \sum_{g\in{\tt G}_{N}({H_i})}\prod_{(u,v)\in{\tt E}(g)}\pk{a_{N,u,v}(t)=1}\notag\\
               &= \sum_{g\in{\tt G}_{N}({H_i})}p_N^{\mathcal{E}({H_i})}(t)= \frac{N!p_N^{{\mathcal{E}}({H_i})}(t)}{(N-\mathcal{V}({H_i}))!\mathcal{A}({H_i})},
\end{align*}
where in the last equality we use \eqref{G_cardinality} for ${H=H_i}$. Hence, \eqref{E_formula} follows. The covariance \eqref{Cov_formula} can be calculated by again applying~\eqref{X_def_with_g}:
\begin{align}
    \operatorname{Cov}\bigl(X_{N,i}(s), X_{N,j}(t)\bigr) &= \operatorname{Cov}\left(\sum_{g\in{\tt G}_N({H_i})}\mathbb{I}\{g{\,{\subset}\,} G_N(s)\}, \sum_{g\in{\tt G}_N({H_j})}\mathbb{I}\{g{\,{\subset}\,} G_N(t)\}\right)\notag\\
    &=\sum_{\substack{g_1\in{\tt G}_N({H_i}) \\ g_2\in{\tt G}_N({H_j})}}\operatorname{Cov}\bigl(\mathbb{I}\{g_1{\,{\subset}\,} G_N(s)\}, \mathbb{I}\{g_2{\,{\subset}\,} G_N(t)\}\bigr).\label{X_cov_1}
\end{align}

 Using  \eqref{I_N_def} and the fact that the processes {$a_{N,u,v}(\cdot)$} are independent for different pairs $(u,v)$, we obtain that the dependence of the random variables $\mathbb{I}\{g_1{\,{\subset}\,} G_N(s)\}$ and $\mathbb{I}\{g_2{\,{\subset}\,} G_N(t)\}$ can be expressed in terms of the intersection ${\tt E}(g_1)\cap {\tt E}(g_2) = {\tt E}(g_1\cap g_2)$, i.e.,
\begin{align}
    \operatorname{Cov}\bigl(\mathbb{I}\{g_1{\,{\subset}\,} G_N(t)\}, &\mathbb{I}\{g_2{\,{\subset}\,} G_N(s)\}\bigr) = \operatorname{Cov}\left(\prod_{(u,v)\in{\tt E}(g_1)}a_{N,u,v}(t), \prod_{(u,v)\in{\tt E}(g_2)}a_{N,u,v}(s)\right)\notag\\
    &= \E{\prod_{(u,v)\in{\tt E}(g_1)\setminus {\tt E}(g_2)}a_{N,u,v}(t)}\E{\prod_{(u,v)\in{\tt E}(g_2)\setminus {\tt E}(g_1)}a_{N,u,v}({s})}\notag\\
    &\qquad\times\operatorname{Cov}\left(\prod_{(u,v)\in{\tt E}(g_1)\cap{\tt E}(g_2)}a_{N,u,v}(t), \prod_{(u,v)\in{\tt E}(g_1)\cap {\tt E}(g_2)}a_{N,u,v}(s)\right)\notag\\
    &= \bigl(p_N(t)\bigr)^{\mathcal{E}({H}_i) - \mathcal{E}(g_1\cap g_2)}\bigl(p_N(s)\bigr)^{\mathcal{E}({H}_j) - \mathcal{E}(g_1\cap g_2)}\notag\\
    &\qquad\times\left( \bigl(\kappa_N(t,s) + p_N(t)p_N(s)\bigr)^{\mathcal{E}(g_1\cap g_2)} - \bigl(p_N(t)p_N(s)\bigr)^{\mathcal{E}(g_1\cap g_2)}\right).\label{X_cov_2}
\end{align}
{Hence, each term in the sum \eqref{X_cov_1} depends on the {pattern} of the intersection $g_1\cap g_2$. We show that the dominant part of the asymptotics of \eqref{X_cov_1} is due to such pairs $g_1,g_2$ that are such that their intersection is isomorphic to {the optimal common subgraph {pattern} of $H_i$ and $H_j$}.} 

At this point we need to distinguish between two different scenarios. {Suppose $\lim_{N\to\infty}\varrho_N = 0$, i.e., we consider the case in which the edge probabilities go to zero as $N\to\infty$.} Then,
as $N\to\infty$, by combining \eqref{X_cov_2} with \eqref{p_asymptotics} and \eqref{kappa_asymptotics}, we obtain that
\begin{align}
    \operatorname{Cov}&\bigl(\mathbb{I}\{g_1{\,{\subset}\,} G_N(s)\}, \mathbb{I}\{g_2{\,{\subset}\,} G_N(t)\}\bigr) =\notag\\ &\varrho_N^{\mathcal{E}({H}_i) + \mathcal{E}({H}_j) - \mathcal{E}(g_1\cap g_2)}\bigl(p_{\star}(s)\bigr)^{\mathcal{E}({H}_i) - \mathcal{E}(g_1\cap g_2)}\bigl(p_{\star}(t)\bigr)^{\mathcal{E}({H}_j) - \mathcal{E}(g_1\cap g_2)}
    \bigl(\kappa_{\star}(s,t)\bigr)^{\mathcal{E}(g_1\cap g_2)}\bigl(1+o(1)\bigr).\label{X_cov_3}
\end{align}
Inserting \eqref{X_cov_3} into \eqref{X_cov_1}, we have, as $N\to\infty$,
\begin{align*}
    \frac{\operatorname{Cov}\bigl(X_{N,i}(t), X_{N,j}(s)\bigr)}{\varrho_N^{\mathcal{E}({H}_i)+\mathcal{E}({H}_j)}\bigl(p_{\star}(s)\bigr)^{\mathcal{E}({H}_i)}\bigl(p_{\star}(t)\bigr)^{\mathcal{E}({H}_j)}}
    &=(1+o(1))\sum_{\substack{g_{\cap}\in{\tt CS}({H}_i,{H}_j),\\ \mathcal{V}{(g_{\cap})\geqslant 1}}}\sum_{\substack{g_1\in{\tt G}_N({H}_i) \\ g_2\in{\tt G}_N({H}_j) \\ {\pi(g_1\cap g_2)= {g_{\cap}}}}}\varrho_N^{-\mathcal{E}(g_{\cap})}\left(\frac{\kappa_{\star}(s,t)}{p_{\star}(t)p_{\star}(s)}\right)^{\mathcal{E}(g_{\cap})}\\
    &=(1+o(1))\sum_{\substack{g_{\cap}\in{\tt CS}({H}_i,{H}_j),\\ \mathcal{V}(g_{\cap})\geqslant 1}}\mathcal{C}({H}_i,{H}_j,g_{\cap})N^{\mathcal{V}({H}_i) + \mathcal{V}({H}_j) - \mathcal{V}(g_{\cap})}\varrho_N^{-\mathcal{E}(g_{\cap})}\\
    &\qquad\qquad\qquad\qquad\qquad\times \left(\frac{\kappa_{\star}(s,t)}{p_{\star}(t)p_{\star}(s)}\right)^{\mathcal{E}(g_{\cap})}\\
    &=(1+o(1)){N^{\mathcal{V}({H}_i) + \mathcal{V}({H}_j)}}{\mathcal{F}_{N}^{\rm opt}({H}_i,{H}_j)}\\
    &\qquad\times \sum_{\substack{g_{\cap}\in{\tt CS}({H}_i,{H}_j),\\ \mathcal{V}(g_{\cap})\geqslant 1}}\mathcal{C}({H}_i,{H}_j,g_{\cap})\frac{\mathcal{F}_N(g_{\cap})}{\mathcal{F}_N^{\rm opt}({H}_i,{H}_j)}\left(\frac{\kappa_{\star}(s,t)}{p_{\star}(t)p_{\star}(s)}\right)^{\mathcal{E}(g_{\cap})}\\
    &=(1+o(1)){N^{\mathcal{V}({H}_i) + \mathcal{V}({H}_j)}}{\mathcal{F}_{N}^{\rm opt}({H}_i,{H}_j)}\\
    &\qquad\qquad\qquad\times \sum_{\substack{g_{\cap}\in{\tt OCS}({H}_i,{H}_j),\\ \mathcal{V}(g_{\cap})\geqslant 1}}\mathcal{C}({H}_i,{H}_j,g_{\cap})\left(\frac{\kappa_{\star}(s,t)}{p_{\star}(t)p_{\star}(s)}\right)^{\mathcal{E}(g_{\cap})},
\end{align*}
where the constant $\mathcal{C}({H}_i,{H}_j,g_{\cap})$ is defined in \eqref{C_def}. We conclude that \eqref{Cov_formula} follows in case $\lim_{N\to\infty}\varrho_N = 0$.

We are left with the case that $\varrho_N=1$. In this case one has that for any $i,j\in\sprod{m}$
\begin{align*}
    \mathcal{F}_N({H}_i) = N^{-\mathcal{V}({H}_i)},\qquad \mathcal{F}_N^{\rm{opt}}({H}_i,{H}_j) = N^{-2}, \qquad {\tt OCS}({H}_i,{H}_j) = \{\mathfrak e\},
\end{align*}
where $\mathfrak{e}$ is {the pattern of an edge}. We observe that the covariance representation obtained in \eqref{X_cov_2} depends only on the number of common edges in the graphs $g_1$ and $g_2$. Hence we can classify all the pairs $g_1,g_2$ on the right hand side of \eqref{X_cov_1} by $k=\mathcal{E}(g_1\cap g_2)$. We thus obtain, {as $N\to\infty$},
\begin{align*}
    \operatorname{Cov}\bigl({X_{N,i}(t), X_{N,i}(s)}\bigr)&=\sum_{k=0}^{\min(\mathcal{E}({H}_i),\mathcal{E}({H}_j))}\sum_{\substack{g_1\in{\tt G}_N({H}_i) \\ g_2\in{\tt G}_N({H}_j) \\ \mathcal{E}(g_1\cap g_2)=k}}\bigl(p_N(t)\bigr)^{\mathcal{E}({H}_i) - k}\bigl(p_N(s)\bigr)^{\mathcal{E}({H}_j) - k}\\
    &\qquad\qquad\qquad\qquad\qquad\qquad\qquad\times\left( \bigl(\kappa_N(t,s) + p_N(t)p_N(s)\bigr)^{k} - \bigl(p_N(t)p_N(s)\bigr)^{k}\right)\\
    &=(1+o(1))\sum_{k=1}^{\min(\mathcal{E}({H}_i),\mathcal{E}({H}_j))}\mathcal{C}_{N,k}({H}_i,{H}_j)\bigl(p_{\star}(t)\bigr)^{\mathcal{E}({H}_i) - k}\bigl(p_{\star}(s)\bigr)^{\mathcal{E}({H}_j) - k}\\
    &\qquad\qquad\qquad\qquad\qquad\qquad\qquad\times\left( \bigl(\kappa_{\star}(t,s) + p_{\star}(t)p_{\star}(s)\bigr)^{k} - \bigl(p_{\star}(t)p_{\star}(s)\bigr)^k\right),
\end{align*}
{where}
\begin{align*}
    \mathcal{C}_{N,k}({H}_i,{H}_j) = \sum_{\substack{g_1\in{\tt G}_N({H}_i) \\ g_2\in{\tt G}_N({H}_j) \\ \mathcal{E}(g_1\cap g_2)=k}} 1.
\end{align*}
The crucial idea behind this representation is that the only factor of each term which depends on $N$ is $\mathcal{C}_{N,k}({H}_i,{H}_j)$. We can find its asymptotic behavior (in the regime that $N\to\infty$) via
\begin{align*}
    \mathcal{C}_{N,1}({H}_i,{H}_j) &= (1+o(1)) \,C({H}_i,{H}_j,\mathfrak{e}) N^{\mathcal{V}({H}_i) + \mathcal{V}({H}_j) - 2},\\
    \mathcal{C}_{N,k}({H}_i,{H}_j) &= o(N^{\mathcal{V}({H}_i) + \mathcal{V}({H}_j) - 2})\qquad \text { for } k\geqslant 2.
\end{align*}
Hence, applying this asymptotic expression, we obtain {that, as $N\to\infty$,}
\begin{align}
    &\operatorname{Cov}\bigl(X_N(t), X_N(s)\bigr) \\
    &\qquad= (1+o(1)) C({H}_i,{H}_j,\mathfrak{e}) N^{\mathcal{V}({H}_i) + \mathcal{V}({H}_j) - 2}\kappa_{\star}(t,s)\bigl(p_{\star}(t)\bigr)^{\mathcal{E}({H}_i) - 1}\bigl(p_{\star}(s)\bigr)^{\mathcal{E}({H}_j) - 1}.\label{X_cov_3_alt}
\end{align}
We conclude that claim \eqref{Cov_formula} follows as well in the case that $\varrho_N=1$.

\QED

\section{Proof of \nelem{lem:asymptotic_normality}}\label{section:proof_of_normality}

Without loss of generality we restrict ourselves to the case where $p_{\star}(t_k)>0$ for all $k\in\sprod{n}$.
{To simplify the upcoming calculations, we define the} {\it normalized process} $\bar{\vk X}_N(t)$, which is defined by, for $t\in\{t_1,\ldots,t_n\}$,
\begin{align*}
    \bar{X}_{N,i}(t) = \frac{X^{\star}_{N,i}(t)}{\bigl(p_{\star}(t)\bigr)^{\mathcal{E}({H}_i)}}{.}
\end{align*}
{Due to \neprop{prop:X_vk_cov}, to} obtain the assertion of \nelem{lem:asymptotic_normality}, it is enough to show that the vector $\bigl(\bar{\vk X}_N(t_1),\ldots,\bar{\vk X}_{N}(t_n)\bigr)$ converges in distributions to a centered Gaussian random vector with a covariance block matrix $\bar{\vk \Sigma}$ with the entries
\begin{align}
\left(\Sigma_{k,l}\right)_{i,j} = \mathbb{I}\{{H}_i\cong {H}_j\}\sum_{g\in{\tt OCS}({H}_i,{H}_j)}\mathcal{C}({H}_i,{H}_j,g)\left(\frac{\kappa_{\star}(t_k,t_l)}{p_{\star}(t_k)p_{\star}(t_l)}\right)^{\mathcal{E}(g)}\label{sigma_klij_def}
\end{align}
for any $k,l\in\sprod{n}$ and any $i,j\in\sprod{m}$.
{To do so} we use the {\it Cram\'er--Wold device}: it is enough to show that, for any vectors $\vk q_1,\ldots \vk q_n\in\R^m$, {the random variable}
\begin{align}
    {\xi_N} = \sum_{k=1}^{n}\sprod{\vk q_k, \bar{\vk X}_N(t_k)}\label{xi_init_def}
\end{align}
converges in distribution to a Gaussian random variable (with the prescribed variance). 
Using \cite[Theorem 30.2]{billingsley2017probability}, this entails that we only need to verify that
\begin{align}
    {\lim_{N\to\infty}\E{\xi_{{N}}^z} = \begin{cases} \sigma^z(z-1)!!,\qquad &\text{if $z$ is even},\\
    0,\qquad &\text{if $z$ is odd,}
    \end{cases}}\label{normality_claim}
\end{align}
for any  $z\in\N$ and the constant 
\begin{align}
    \sigma^2 = \sum_{i,j=1}^{m}\sum_{k,l=1}^{n}q_{k,i}q_{l,j}\left(\Sigma_{k,l}\right)_{i,j}.\label{sigma_def}
\end{align}
From \eqref{X_def_with_g} we know that $X_{N,i}^{\star}(t)$, as {defined in \eqref{X_star_def}}, can be represented as a sum of centered indicators over all subgraphs $g\in{\tt G}_N({H}_i)$. We would like to {obtain} a similar {representation} for $\bar{X}_{N,i}(t)$.  To this end, define, for any $i\in\sprod{m}$ and any $g\in{\tt G}_{N}({H}_i)$,
\begin{align}
    \xi_{i,{N}}(g) = \sum_{k=1}^{n}\frac{q_{k,i}\bigl(\mathbb{I}\{g{\,{\subset}\,} G_N(t_k)\} - \pk{g{\,{\subset}\,} G_N(t_k)}\bigr)}{\bigl(p_{\star}(t_k)\bigr)^{\mathcal{E}({H}_i)}}.\label{xi_def}
\end{align} 
This means that, appealing to \eqref{X_def_with_g}, we can represent $\xi_{{N}}$ defined in \eqref{xi_init_def} as a sum over all possible $\xi_{i,{N}}(g)$:  
\begin{align}
    \xi_{{N}} = \sum\limits_{i=1}^{m}\frac{\mathcal{F}({H}_i)}{\sqrt{\mathcal{F}^{\rm opt}_N({H}_i)}}\sum\limits_{g\in{\tt G}_N({H}_i)}\xi_{i,{N}}(g).\label{xi_sum_repr}
\end{align}
Hence, taking the expectation of the sum on the right hand side of \eqref{xi_sum_repr} raised to the power $z$, we obtain the following  expression of $z$-th moment of $\xi_{{N}}$: 
\begin{align}
    \E{\xi_{{N}}^z} &= \sum_{f_1,\ldots,f_{z}\in\sprod{m}}\left(\prod\limits_{i=1}^{z}\frac{\mathcal{F}({H}_{f_i})}{\sqrt{\mathcal{F}^{\rm opt}_N({H}_{f_i})}}\right)\sum_{\substack{g_1\in{\tt G}_N({H}_{f_1}) \\ \ldots \\ g_{z}\in{\tt G}_N({H}_{f_z})}}\E{\prod_{i=1}^{z}\xi_{f_i,{N}}(g_i)}.\label{moment_calculation_1}
\end{align}
To simplify the notation, {for any $\vk f = (f_{1},\ldots,f_z)^{\top}\in\R^z$} introduce the set
\begin{align*}
    {\tt G}_{\vk f,N} = \bigtimes\limits_{i=1}^{z}{\tt G}_{N}({H}_{f_i}).
\end{align*}
Then, the last sum in \eqref{moment_calculation_1} is simply a sum over ${\vk g = (g_1,\ldots,g_z)^{\top}\in{\tt G}_{\vk f,N}}$. {The primary objective of this proof is to demonstrate that \eqref{moment_calculation_1} satisfies the asymptotics given by \eqref{normality_claim}. The first step involves identifying the terms in the inner sum of \eqref{moment_calculation_1} that are non-zero. Next, we isolate the terms within that sum which make the dominant contribution to the asymptotics. Finally, we compute the asymptotics of the identified dominant terms.} We start the {our first step} by investigating for which collections $\vk g$ the expectation at the end of $\eqref{moment_calculation_1}$ is different from zero. Using \eqref{I_N_def} we observe that if $\mathcal{E}(g^{\star}\cap g^{\star\star}) = 0$ for some subgraphs $g^{\star},g^{\star\star}{\,{\subset}\,} K_N$, then, for any $i,i^{\star}\in\sprod{m}$ the random variables $\xi_{i,{N}}(g^{\star})$ and $\xi_{i^{\star},{N}}(g^{\star\star})$ are independent. Hence, if for some $j\in\sprod{z}$ and all $i\in\sprod{k}$, $i\not= j$ we have that $\mathcal{E}(g_{i}\cap g_{j}) = 0$, then
\begin{align}
    \E{\prod_{i=1}^{z}\xi_{f_i,{N}}(g_i)} = \E{\prod_{\substack{i=1 \\ i\not= j}}^{z}\xi_{f_i,{N}}(g_i)}\E{\xi_{f_{j},{N}}(g_{j})} = 0.\label{exp_decomposition}
\end{align}
Thus, we can consider only such collections $\vk g\in{\tt G}_{N,\vk f}$ which belongs to the following set
\begin{align}
    {\tt G}^{\star}_{\vk f, N} = \bigl\{\vk g\in{\tt G}_{\vk f,N}\colon\forall j\in\sprod{m} \,\exists i\in\sprod{m}\setminus j\colon \mathcal{E}(g_i\cap g_j)\geqslant 1\bigr\}.\label{g_condition}
\end{align}
Applying this observation,  \eqref{moment_calculation_1} can be expressed as
\begin{align}
     \E{\xi_{{N}}^z} = \sum_{f_1,\ldots,f_{z}\in\sprod{m}}\sum_{\vk g\in{\tt G}^{\star}_{\vk f,N}}\E{\prod_{i=1}^{z}\frac{\mathcal{F}_{{N}}({H}_{f_i})\xi_{f_i,{N}}(g_i)}{\sqrt{\mathcal{F}^{\rm opt}_N({H}_{f_i})}}}.\label{moment_calculation_2}
\end{align}
{We now move to the next step of our proof: identifying the terms in the inner sum of \eqref{moment_calculation_2} that determine the asymptotics of this sum (i.e., the remaining terms are asymptotically negligible).}
To do so,  classify all the collections $\vk g$ that appear in \eqref{moment_calculation_2} by the number of vertices in their union, i.e., the value $\mathcal{V}(\vk g) := \mathcal{V}\left(\bigcup_{i=1}^{z}g_i\right)$. In this case we can first choose a complete graph $K_{\mathcal{V}(\vk g)}{\,{\subset}\,} K_N$ where the graph $\bigcup_{i=1}^{z}g_i$ is located (where it is noted that the number of such choices is dependent on $N$), after which we choose $g_i{\,{\subset}\,} K_{\mathcal{V}(\vk g)}$ (where the number of such choices no longer depends neither on $N$, nor on the particular choice of the graph $K_{\mathcal{V}(\vk g)}$). Let{
\begin{align}
    V_{\max} = \max\limits_{\vk g\in{\tt G}^{\star}_{\vk f,N}}\mathcal{V}(\vk g).\label{V_max_def}
\end{align}}
Note that constant $V_{\max}$ does not depend on $N$, based on the above considerations. Then
\begin{align}
     \sum_{\vk g\in{\tt G}^{\star}_{\vk f,N}}\E{\prod_{i=1}^{z}\frac{\mathcal{F}_{{N}}({H}_{f_i})\xi_{f_i,{N}}(g_i)}{\sqrt{\mathcal{F}^{\rm opt}_N({H}_{f_i})}}} &= \sum_{{v=2}}^{V_{\max}}\sum_{\substack{\vk g\in{\tt G}^{\star}_{\vk f,N}\\ \mathcal{V}\left(\vk g\right) = v }}\E{\prod_{i=1}^{z}\frac{\mathcal{F}_{{N}}({H}_{f_i})\xi_{f_i,{N}}(g_i)}{\sqrt{\mathcal{F}^{\rm opt}_N({H}_{f_i})}}}\notag\\
     &= \sum_{{v=2}}^{V_{\max}}\sum_{K_v{\,{\subset}\,} K_N}\sum_{\substack{\vk g\in{\tt G}^{\star}_{\vk f,v}\\
     \mathcal{V}\left(\vk g\right) = v }}\E{\prod_{i=1}^{z}\frac{\mathcal{F}_{{N}}({H}_{f_i})\xi_{f_i,{N}}(g_i)}{\sqrt{\mathcal{F}^{\rm opt}_N({H}_{f_i})}}}\notag\\
     &= \sum_{{v=2}}^{V_{\max}}\frac{N!}{(N-v)!v!}\sum_{\substack{\vk g\in{\tt G}^{\star}_{\vk f,v}\\
     \mathcal{V}\left(\vk g\right) = v }}\E{\prod_{i=1}^{z}\frac{\mathcal{F}_{{N}}({H}_{f_i})\xi_{f_i,{N}}(g_i)}{\sqrt{\mathcal{F}^{\rm opt}_N({H}_{f_i})}}}\notag\\
     &= \sum_{{v=2}}^{V_{\max}}\frac{N!}{N^{v}(N-v)!v!}\sum_{\substack{\vk g\in{\tt G}^{\star}_{\vk f,v}\\
     \mathcal{V}\left(\vk g\right) = v }}\E{N^{\mathcal{V}(\vk{g})}\prod_{i=1}^{z}\frac{\mathcal{F}_{{N}}({H}_{f_i})\xi_{f_i,{N}}(g_i)}{\sqrt{\mathcal{F}^{\rm opt}_N({H}_{f_i})}}},\label{moment_calculation_3}
\end{align}
 where $V_{\max}$ is the maximal possible value of $\mathcal{V}\left(\vk g\right)$ over $\vk g\in{\tt G}^{\star}_{\vk f,N}$. 
{To analyze the expectations on the right-hand side of \eqref{moment_calculation_3}, it is useful to apply a factorization in terms of expectations for independent groups of random variables $\xi_{f_i,N}(g_i)$. From the definition \eqref{xi_def}, it follows that the random variables $\xi_{f_i,N}(g_i)$ and $\xi_{f_j,N}(g_j)$ are independent if and only if the graphs $g_i$ and $g_j$ do not share any common edges. Therefore, we need to partition the graphs $g_i$ for $i \in \sprod{z}$ into groups such that any two graphs $g_i$ and $g_j$ from different groups do not have any common edges.}

{To formalize this idea,} {define for $\vk g\in\tt G^{\star}_{\vk f,N}$}  a graph $\mathcal{G}(\vk g)$ with vertices ${\tt V}\bigl(\mathcal{G}(\vk g)\bigr)=\sprod{z}$, and the edge between vertices $i,j\in\sprod{z}$ exists if and only if $\mathcal{E}(g_i\cap g_j)\geqslant 1$. We define by ${\tt CC}\bigl(\mathcal{G}(\vk g)\bigr)$ the set of connected components of the graph $\mathcal{G}(\vk g)$. Then, {if $i$ and $j$ belong to different connected components of the graph $\mathcal{G}(\vk g)$, the  random variables $\xi_{f_i, N}(g_i)$ and $\xi_{f_j, N}(g_j)$ are independent.} {Hence,} similarly to how \eqref{exp_decomposition} was obtained, we see that
\begin{align}
    \E{N^{\mathcal{V}(\vk{g})}\prod_{i=1}^{z}\frac{\mathcal{F}_{{N}}({H}_{f_i})\xi_{f_i,{N}}(g_i)}{\sqrt{\mathcal{F}^{\rm opt}_N({H}_{f_i})}}} = \prod_{{\tt G}\in {\tt CC}(\mathcal{G}(\vk g))} \E{N^{\mathcal{V}({\tt G})}\prod_{i\in{\tt V}({\tt G})}\frac{\mathcal{F}_{{N}}({H}_{f_i})\xi_{f_i,{N}}(g_i)}{\sqrt{\mathcal{F}^{\rm opt}_N({H}_{f_i})}}}\label{E_approx_prod},
\end{align}
where $\mathcal{V}({\tt G}) = \mathcal{V}(\cup_{i\in{\tt V}({\tt G})}g_i)$. Note that if at least for one connected component of $\mathcal{G}(\vk g)$ the expectation right hand side of \eqref{E_approx_prod} tends to zero as $N\to\infty$, and for the rest connected components the expectations are bounded, we immediately obtain that for such $\vk g\in{\tt G}^{\star}_{\vk f,v}$
\begin{align*}
    \lim_{N\to\infty}\E{N^{\mathcal{V}(\vk{g})}\prod_{i=1}^{z}\frac{\mathcal{F}_{{N}}({H}_{f_i})\xi_{f_i,{N}}(g_i)}{\sqrt{\mathcal{F}^{\rm opt}_N({H}_{f_i})}}} = 0,
\end{align*}
and we can disregard them in the inner sum of \eqref{moment_calculation_3}. The following lemma, proven in the appendix, tells us how to find such a $\vk g$.

\begin{lem}\label{lem:connected_components}
    Let $\vk g\in{\tt G}^{\star}_{\vk f,v}$ and ${\tt G}\in{\tt CC}(\mathcal{G}(\vk g))$. Then
    \begin{align*}
        \lim_{N\to\infty}\E{N^{\mathcal{V}({\tt G})}\prod_{i\in{\tt V}({\tt G})}\frac{\mathcal{F}_{{N}}({H}_{f_i})\xi_{f_i,{N}}(g_i)}{\sqrt{\mathcal{F}^{\rm opt}_N({H}_{f_i})}}}\in[0,\infty).
    \end{align*}
    If also $\mathcal{V}({\tt G})\geqslant 3$ then
    \begin{align*}
        \lim_{N\to\infty}\E{N^{\mathcal{V}({\tt G})}\prod_{i\in{\tt V}({\tt G})}\frac{\mathcal{F}_{{N}}({H}_{f_i})\xi_{f_i,{N}}(g_i)}{\sqrt{\mathcal{F}^{\rm opt}_N({H}_{f_i})}}}=0.
    \end{align*}
\end{lem}

\noindent
{\it Continuation of proof of \nelem{lem:asymptotic_normality}}. Recall that from \eqref{g_condition} we know that, for any $\vk g\in{{\tt G}^{\star}_{\vk f,v}}$ each connected component of $\mathcal{G}(\vk g)$ has at least two vertices. Thus, {if the graph $\mathcal{G}(\vk g)$ contains} at least one connected component with at least three vertices, {then} due to \nelem{lem:connected_components}, in combination with \eqref{moment_calculation_2} and \eqref{moment_calculation_3}, we obtain that
\begin{align*}
    \lim_{N\to\infty} \E{\xi_{{N}}^z} &= \lim_{N\to\infty}\sum_{v=0}^{V_{\max}}\frac{N!}{N^{v}(N-v)!v!}\sum_{\substack{\vk g\in{\tt G}^{\star}_{\vk f,v}\\
     \mathcal{V}\left(\vk g\right) = v }}\E{N^{\mathcal{V}(\vk{g})}\prod_{i=1}^{z}\frac{\mathcal{F}_{{N}}({H}_{f_i})\xi_{f_i,{N}}(g_i)}{\sqrt{\mathcal{F}^{\rm opt}_N({H}_{f_i})}}}\\
     &= \sum_{v=0}^{V_{\max}}\frac{1}{v!}\sum_{\substack{\vk g\in{\tt G}^{\star}_{\vk f,v}\\
     \mathcal{V}\left(\vk g\right) = v }}\lim_{N\to\infty}\E{N^{\mathcal{V}(\vk{g})}\prod_{i=1}^{z}\frac{\mathcal{F}_{{N}}({H}_{f_i})\xi_{f_i,{N}}(g_i)}{\sqrt{\mathcal{F}^{\rm opt}_N({H}_{f_i})}}}=0.
\end{align*}
{Hence, for odd $z$, \eqref{normality_claim} follows from \eqref{moment_calculation_2}, because the graph $\mathcal{G}(\vk g)$ has an odd number of vertices and no isolated vertices, which implies that $\mathcal{G}(\vk g)$ must have at least one connected component with at least three vertices. 

In addition, for even $z$, we restrict our attention to the ${\vk g} \in {\tt G}_{\vk f,v}$ for which every connected component of $\mathcal{G}(\vk g)$ consists of exactly two vertices, meaning that it forms a {`perfect matching'} (again, because $\mathcal{G}(\vk g)$ has no isolated vertices). 

At this point, we complete the second step of our proof: we show that the dominant contribution to the asymptotics of \eqref{moment_calculation_2} is provided by those $\vk g$ for which $\mathcal{G}(\vk g)$ is a {perfect matching}. Let us denote by ${\tt M}_z$ the set of all perfect matchings of $K_z$. This implies, by \eqref{moment_calculation_3}, that, as $N \to \infty$,} 
\begin{align}
    \sum_{\vk g\in{\tt G}^{\star}_{\vk f,N}}\E{\prod_{i=1}^{z}\frac{\mathcal{F}_{{N}}({H}_{f_i})\xi_{f_i,{N}}(g_i)}{\sqrt{\mathcal{F}^{\rm opt}_N({H}_{f_i})}}}={\bigl(1+o(1)\bigr)}\sum_{v=0}^{V_{\max}}\frac{1}{v!}\sum_{{H}\in{\tt M}_z}\sum_{\substack{\vk g\in{\tt G}^{\star}_{\vk f,v}\\
     \mathcal{V}\left(\vk g\right) = v \\ \mathcal{G}(\vk g)={H}}}\E{N^{\mathcal{V}(\vk{g})}\prod_{i=1}^{z}\frac{\mathcal{F}_{{N}}({H}_{f_i})\xi_{f_i,{N}}(g_i)}{\sqrt{\mathcal{F}^{\rm opt}_N({H}_{f_i})}}}\notag\\
     = {\bigl(1+o(1)\bigr)}\sum_{{H}\in{\tt M}_z}\sum_{v=0}^{V_{\max}}\frac{1}{v!}\sum_{\substack{\vk g\in{\tt G}^{\star}_{\vk f,v}\\
     \mathcal{V}\left(\vk g\right) = v \\ \mathcal{G}(\vk g)={H}}}\prod_{(i,j)\in{\tt E}({H})}\E{N^{\mathcal{V}(g_i\cup g_j)}\frac{\mathcal{F}_{{N}}({H}_{f_i})\xi_{f_i}(g_i)}{\sqrt{\mathcal{F}^{\rm opt}_N({H}_{f_i})}}\frac{\mathcal{F}_{{N}}({H}_{f_j})\xi_{f_j,{N}}(g_j)}{\sqrt{\mathcal{F}^{\rm opt}_N({H}_{f_j})}}}.\label{moment_calculation_4}
\end{align}
{We now proceed to the final step of our proof: calculating the asymptotics of \eqref{moment_calculation_4}. We start by deriving the asymptotics of each inner expectation in \eqref{moment_calculation_4}.} From \eqref{xi_def} we conclude that
\begin{align}
    \E{N^{\mathcal{V}({g_i}\cup g_j)}\frac{\mathcal{F}_{{N}}({H}_{f_i})\xi_{f_i,{N}}(g_i)}{\sqrt{\mathcal{F}^{\rm opt}_N({H}_{f_i})}}\frac{\mathcal{F}_{{N}}({H}_{f_j})\xi_{f_j,{N}}(g_j)}{\sqrt{\mathcal{F}^{\rm opt}_N({H}_{f_j})}}} &= \frac{N^{\mathcal{V}(g_i\cup g_j)}\mathcal{F}_{{N}}({H}_{f_i})\mathcal{F}_{{N}}({H}_{f_j})}{\sqrt{\mathcal{F}^{\rm opt}_N({H}_{f_i})}\sqrt{\mathcal{F}^{\rm opt}_N({H}_{f_j})}}\notag\\&\hspace{-10mm}
\times\sum_{k,l=1}^{n}\frac{q_{k,i}q_{l,j}\operatorname{Cov}\bigl(\mathbb{I}\{g_i{\,{\subset}\,} G_N(t_{k})\},\mathbb{I}\{g_j{\,{\subset}\,}G_N(t_{l})\}\bigr)}{\bigl(p_{\star}(t_{k})\bigr)^{\mathcal{E}({H}_i)}\bigl(p_{\star}(t_{l})\bigr)^{\mathcal{E}({H}_j)}}.\label{moment_calculation_5}
\end{align}
From \eqref{X_cov_3} we know that, as $N\to\infty$,
\begin{align*}
    \operatorname{Cov}\bigl(\mathbb{I}\{g_i\subset G_N(t_{k})\},\mathbb{I}\{g_j\subset G_N(t_{l})\}\bigr)&={\bigr(1+o(1)\bigl)} \varrho_N^{\mathcal{E}({H}_i) + \mathcal{E}({H}_j) - \mathcal{E}(g_i\cap g_j)}\bigl(p_{\star}({t_k})\bigr)^{\mathcal{E}({H}_i) - \mathcal{E}(g_i\cap g_j)}\notag\\
    &\qquad\times\bigl(p_{\star}({t_l})\bigr)^{\mathcal{E}({H}_j) - \mathcal{E}(g_i\cap g_j)}\bigl(\kappa_{\star}({t_k,t_l})\bigr)^{\mathcal{E}(g_i\cap g_j)}.
\end{align*}
Hence, using that
\begin{align*}
    N^{\mathcal{V}(g_i\cup g_j)}\mathcal{F}_{{N}}({H}_{f_i})\mathcal{F}_{{N}}({H}_{f_j})\varrho_N^{\mathcal{E}({H}_i) + \mathcal{E}({H}_j) - \mathcal{E}({g_i\cap g_j})} &= N^{\mathcal{V}(g_i\cup g_j)}N^{- \mathcal{V}(g_i)}N^{-\mathcal{V}(g_j)}\\
    &\qquad\times\varrho_N^{-\mathcal{E}(g_i)}\varrho_N^{-\mathcal{E}(g_j)}\varrho_N^{\mathcal{E}({H}_i) + \mathcal{E}({H}_j) - \mathcal{E}(g_1\cap g_2)}\\
    &=N^{-\mathcal{V}(g_i\cap g_j)}\varrho_N^{-\mathcal{E}(g_i\cap g_j)}=\mathcal{F}_{{N}}(g_i\cap g_j),
\end{align*}
we obtain from \eqref{moment_calculation_5} that, as $N\to\infty$,
\begin{align}
    \E{N^{\mathcal{V}({g_i}\cup g_j)}\frac{\mathcal{F}_{{N}}({H}_{f_i})\,\xi_{f_i,{N}}(g_i)}{\sqrt{\mathcal{F}^{\rm opt}_N({H}_{f_i})}}\frac{\mathcal{F}_{{N}}({H}_{f_j})\,\xi_{f_j,{N}}(g_j)}{\sqrt{\mathcal{F}^{\rm opt}_N({H}_{f_j})}}}\notag\\&\hspace{-6.55cm}={\bigl(1+o(1)\bigr)} \sqrt{\frac{\mathcal{F}_N(g_i\cap g_j)}{\mathcal{F}^{\rm{opt}}_N({H}_{f_i})}}\sqrt{\frac{\mathcal{F}_N(g_i\cap g_j)}{\mathcal{F}^{\rm{opt}}_N({H}_{f_j})}}\label{moment_calculation_5.1}
    \times\sum_{k,l=1}^{n}\frac{q_{k,i}q_{l,j}\bigl(\kappa_{\star}({t_k,t_l})\bigr)^{\mathcal{E}(g_i\cap g_j)}}{\bigl(p_{\star}({t_k})\,p_{\star}({t_l})\bigr)^{\mathcal{E}(g_i\cap g_j)}}.
\end{align}
{Recall that from \eqref{F_cal_opt_def}, as $g_i\cap g_j {\,{\subset}\,} g_i\sim {H}_{f_i}$, we know that
\begin{align*}
    \lim_{N\to\infty}\sqrt{\frac{\mathcal{F}_N(g_i\cap g_j)}{\mathcal{F}^{\rm{opt}}_N({H}_{f_i})}} = \mathbb{I}\{{\pi(g_i\cap g_j) \in\,}{\tt OCS}({H}_{f_i})\}.
\end{align*}
By a similar reasoning 
\begin{align*}
    \lim_{N\to\infty}\sqrt{\frac{\mathcal{F}_N(g_i\cap g_j)}{\mathcal{F}^{\rm{opt}}_N({H}_{f_j})}} = \mathbb{I}\{{\pi(g_i\cap g_j)\in\,}{\tt OCS}({H}_{f_j})\}.
\end{align*}
}{Thus, combining \eqref{GCS_cong} and \eqref{GCS_not_cong}, from \eqref{moment_calculation_5.1} we obtain that} 
\begin{align}
    \lim_{N\to\infty}\E{N^{\mathcal{V}(g_1\cup g_j)}\frac{\mathcal{F}_{{N}}({H}_{f_i})\xi_{f_i,{N}}(g_i)}{\sqrt{\mathcal{F}^{\rm opt}_N({H}_{f_i})}}\frac{\mathcal{F}_{{N}}({H}_{f_j})\xi_{f_j,{N}}(g_j)}{\sqrt{\mathcal{F}^{\rm opt}_N({H}_{f_j})}}} \notag\\& \hspace{-8cm}=\mathbb{I}\{{H}_{f_i}\cong {H}_{f_j}\}\mathbb{I}\{{\pi(g_i\cap g_j)\in\,}{\tt OCS}({H}_{f_i},{H}_{f_j})\}\label{moment_calculation_6}\times\sum_{k,l=1}^{n}\frac{q_{k,i}q_{l,j}\bigl(\kappa_{\star}({t_k,t_l})\bigr)^{\mathcal{E}(g_i\cap g_j)}}{\bigl(p_{\star}({t_k})p_{\star}({t_l})\bigr)^{\mathcal{E}(g_i\cap g_j)}}.
\end{align}
(As a side remark we notice that even though \eqref{X_cov_3} holds only in case $\varrho_N\to 0$ as $N\to\infty$, one can still obtain the asymptotic relation \eqref{moment_calculation_6} even in case $\varrho_N=1$).
Plugging \eqref{moment_calculation_6} into \eqref{moment_calculation_4}, {we conclude} 
\begin{align}
    {\lim_{N\to\infty}}\sum_{\vk g\in{\tt G}^{\star}_{\vk f,N}}\E{\prod_{i=1}^{z}\frac{\mathcal{F}_{{N}}({H}_{f_i})\xi_{f_i,{N}}(g_i)}{\sqrt{\mathcal{F}^{\rm opt}_N({H}_{f_i})}}}
    &=\sum_{{H}\in{\tt M}_z}\sum_{v=0}^{V_{\max}}\frac{1}{v!}\sum_{\substack{\vk g\in{\tt G}^{\star}_{\vk f,v}\\
     \mathcal{V}\left(\vk g\right) = v \\ \mathcal{G}(\vk g)={H}}}{\mathbb J}({\vk f},{\vk g},H)\notag\\
     &{=\sum_{{H}\in{\tt M}_z}\sum_{v=0}^{V_{\max}}\sum_{\substack{\vk g\in{\tt G}^{\star}_{\vk f}\\
     \mathcal{V}\left(\vk g\right) = v \\ \mathcal{G}(\vk g)={H}}}\frac{1}{\mathcal{V}(\vk g)!}{\mathbb J}({\vk f},{\vk g},H)},\label{moment_calculation_6.5}
\end{align}
where we denote 
\begin{align*}
    {\tt G}_{\vk f}^{\star} &= \bigcup_{v=0}^{V_{\max}}\{\vk g\in{\tt G}_{\vk f, v}^{\star}\colon \mathcal{V}(\vk g) = v\},\\
    {\mathbb J}({\vk f},{\vk g},H)&=\prod_{(i,j)\in{\tt E}({H})}\mathbb{I}\{{H}_{f_i}\cong {H}_{f_j}\}\Bigg(\mathbb{I}\{{\pi(g_i\cap g_j)\in\,}{\tt OCS}({H}_{f_i},{H}_{f_j})\}\sum_{k,l=1}^{n}\frac{q_{k,i}q_{l,j}\bigl(\kappa_{\star}(t_k,t_l)\bigr)^{\mathcal{E}(g_i\cap g_j)}}{\bigl(p_{\star}(t_k)\,p_{\star}(t_l)\bigr)^{\mathcal{E}(g_i\cap g_j)}}\Bigg).
\end{align*}
{The last equality in \eqref{moment_calculation_6.5} holds because
\begin{align*}
    \bigl\{\vk g\colon \vk g\in {\tt G}_{\vk f}^{\star},\, \mathcal{V}(\vk g) = v\bigr\} = \bigl\{\vk g\colon \vk g\in {\tt G}_{\vk f, v}^{\star},\, \mathcal{V}(\vk g) = v\bigr\}.
\end{align*}

\noindent 
Note that each \( \vk g \in {\tt G}_{\vk f}^{\star} \) from the third sum in \eqref{moment_calculation_6.5} corresponds uniquely to a specific \( v \in \sprod{V_{\max}} \) in the second sum (where \( v = \mathcal{V}(\vk g) \)). { Consequently, recalling the definition of $V_{\max}$ given in \eqref{V_max_def},
\begin{align*}
    \bigcup_{v=1}^{V_{\max}}\left\{\vk g\in{\tt G}^{\star}_{\vk f}\colon \mathcal{V}(\vk g) = v,\,\mathcal{G}(\vk g) = H\right\} &= \left\{\vk g\in{\tt G}_{\vk f}^{\star}\colon \mathcal{G}(\vk g) = H,\, \mathcal{V}(\vk g)\in\sprod{V_{\max}}\right\}\\
    &= \left\{\vk g\in{\tt G}_{\vk f}^{\star}\colon \mathcal{G}(\vk g) = H\right\},
\end{align*}
and for any $v_1,v_2\in\sprod{V_{\max}}$, $v_1\not= v_2$
\begin{align*}
    \left\{\vk g\in{\tt G}^{\star}_{\vk f}\colon \mathcal{V}(\vk g) = v_1,\,\mathcal{G}(\vk g) = H\right\}\cap \left\{\vk g\in{\tt G}^{\star}_{\vk f}\colon \mathcal{V}(\vk g) = v_2,\,\mathcal{G}(\vk g) = H\right\} = \varnothing,
\end{align*}
implying that
\begin{align*}
    \sum_{v=0}^{V_{\max}} \sum_{\substack{\vk g \in {\tt G}^{\star}_{\vk f} \\ 
    \mathcal{V}\left(\vk g\right) = v \\ \mathcal{G}(\vk g) = {H}}} 
    \frac{1}{\mathcal{V}(\vk g)!}{\mathbb J}({\vk f},{\vk g},H) =  \sum_{\substack{\vk g \in {\tt G}^{\star}_{\vk f} \\ \mathcal{G}(\vk g) = {H}}}\frac{1}{\mathcal{V}(\vk g)!}{\mathbb J}({\vk f},{\vk g},H).
\end{align*}}
Hence, from \eqref{moment_calculation_6.5} we can obtain that}
\begin{align}
    {\lim_{N\to\infty}}\sum_{\vk g\in{\tt G}^{\star}_{\vk f,N}}\E{\prod_{i=1}^{z}\frac{\mathcal{F}_{{N}}({H}_{f_i})\xi_{f_i,{N}}(g_i)}{\sqrt{\mathcal{F}^{\rm opt}_N({H}_{f_i})}}}
     &{=\sum_{{H}\in{\tt M}_z}\sum_{\substack{\vk g\in
     {\tt G}^{\star}_{\vk f} \\ \mathcal{G}(\vk g)={H}}}\frac{1}{\mathcal{V}(\vk g)!}{\mathbb J}({\vk f}{,\,}{\vk g},H)},\label{moment_calculation_7}
\end{align}
{Denote the edges of graph ${H}\in{\tt M}_z$ as follows: ${\tt E}({H}) = \{(i_1,j_1),\ldots,(i_{z/2},j_{z/2})\}$. Note that the right-hand side of \eqref{moment_calculation_7} does not depend on the full structure of vector $\vk g$, but only on the {patterns} of graphs $g_{i_b}\cap g_{j_b}$, for $b\in\sprod{z/2}$. Hence we decompose this sum by considering the various {patterns $g_b^{\cap}$} of the graphs {$g_{i_b}\cap g_{j_b}$}, {noting that due to the indicators appearing in the right-hand side of \eqref{moment_calculation_7} we can restrict ourselves to ${g_b^{\cap} \in} {\tt OCS}(H_{f_{i_b}}, H_{f_{j_b}})$ only}, so as to obtain 
\begin{align}
    &\lim_{N\to\infty}\sum_{\vk g\in{\tt G}^{\star}_{\vk f,N}}\E{\prod_{i=1}^{z}\frac{\mathcal{F}_{{N}}({H}_{f_i})\xi_{f_i,{N}}(g_i)}{\sqrt{\mathcal{F}^{\rm opt}_N({H}_{f_i})}}}
    \notag\\
     &=\sum_{{H}\in{\tt M}_z}\sum_{\substack{\vk g\in
     {\tt G}^{\star}_{\vk f} \\ \mathcal{G}(\vk g)={H}}}\frac{1}{\mathcal{V}(\vk g)!}\prod_{b=1}^{z/2}\mathbb{I}\{{H}_{f_{i_b}}\cong {H}_{f_{j_b}}\}\Bigg(\mathbb{I}\{g_b^{\cap}\in{\tt OCS}({H}_{f_i},{H}_{f_j})\}\sum_{k,l=1}^{n}\frac{q_{k,i}q_{l,j}\bigl(\kappa_{\star}(k,l)\bigr)^{\mathcal{E}(g_b^{\cap})}}{\bigl(p_{\star}(k)\,p_{\star}(l)\bigr)^{\mathcal{E}(g_b^{\cap})}}\Bigg)\notag\\
     &=\sum_{{H}\in{\tt M}_z}\sum_{\vk g^{\cap}\in{\tt OCS}_{\vk f}({H})}\sum_{\vk g\in
     {\tt G}^{\star\star}_{\vk f}({H},\vk g^{\cap}) }\frac{1}{\mathcal{V}(\vk g)!}\prod_{b=1}^{z/2}\mathbb{I}\{{H}_{f_{i_b}}\cong {H}_{f_{j_b}}\}\sum_{k,l=1}^{n}\frac{q_{k,i}q_{l,j}\bigl(\kappa_{\star}(k,l)\bigr)^{\mathcal{E}(g_b^{\cap})}}{\bigl(p_{\star}(k)\,p_{\star}(l)\bigr)^{\mathcal{E}(g_b^{\cap})}},\label{moment_calculation_7.1}
\end{align}
where
\begin{align*}
    {\tt G}^{\star\star}_{\vk f}({H},\vk g^{\cap}) := \left\{\vk g\in{\tt G}_{\vk f}^{\star}\colon \mathcal{G}(\vk g) = {H},\, \forall b\in\sprod{z/2}:{\pi(g_{i_b}\cap g_{j_b}) =} g^{\cap}_{b}\right\},
\end{align*}
and 
\begin{align*}
    {\tt OCS}_{\vk f}({H}) = \prod_{b=1}^{z/2}{\tt OCS}({H}_{f_{i_b}},{H}_{f_{j_b}}).
\end{align*}
}Notice that for any $\vk g\in {\tt G}^{\star\star}_{\vk f}({H},\vk g^{\cap})$
\begin{align*}
    \mathcal{V}(\vk g) = \sum_{b=1}^{z}\mathcal{V}({H}_{f_b}) - \sum_{b=1}^{z/2}\mathcal{V}(g_b^{\cap}).
\end{align*}
Hence, {the terms in the inner sum of \eqref{moment_calculation_7.1} do not depend on $\vk g$, allowing us to rewrite this sum as} 
\begin{align}
    \lim_{N\to\infty}\sum_{\vk g\in{\tt G}^{\star}_{\vk f,N}}\E{\prod_{i=1}^{z}\frac{\mathcal{F}_{{N}}({H}_{f_i})\xi_{f_i,{N}}(g_i)}{\sqrt{\mathcal{F}^{\rm opt}_N({H}_{f_i})}}}&= \sum_{{H}\in{\tt M}_z}\sum_{\vk g^{\cap}\in{\tt OCS}_{\vk f}({H})}\sum_{\vk g\in{\tt G}^{\star\star}_{\vk f}({H},\vk g^{\cap})}\frac{\prod_{b=1}^{z/2}\mathbb{I}\{{H}_{f_{i_b}}\cong {H}_{f_{j_b}}\}}{\bigl(\sum_{b=1}^{z}\mathcal{V}({H}_{f_b}) - \sum_{b=1}^{z/2}\mathcal{V}(g_b^{\cap})\bigr)!}\notag\\
    &\qquad\qquad\qquad\times\prod_{b=1}^{z/2}\sum_{k,l=1}^{n}\frac{q_{k,i}q_{l,j}\bigl(\kappa_{\star}({t_k,t_l})\bigr)^{\mathcal{E}(g_b^{\cap})}}{\bigl(p_{\star}({t_k})\,p_{\star}({t_l})\bigr)^{\mathcal{E}(g_b^{\cap})}}\notag\\
    &= \sum_{{H}\in{\tt M}_z}\sum_{\vk g^{\cap}\in{\tt OCS}_{\vk f}({H})}\frac{\prod_{b=1}^{z/2}\mathbb{I}\{{H}_{f_{i_b}}\cong {H}_{f_{j_b}}\}\abs{{\tt G}^{\star\star}_{\vk f}({H},\vk g^{\cap})}}{\bigl(\sum_{b=1}^{z}\mathcal{V}({H}_{f_b}) - \sum_{b=1}^{z/2}\mathcal{V}(g_b^{\cap})\bigr)!}\label{moment_calculation_8}\\
    &\qquad\qquad\qquad\times\prod_{b=1}^{z/2}\sum_{k,l=1}^{n}\frac{q_{k,i}q_{l,j}\bigl(\kappa_{\star}({t_k,t_l})\bigr)^{\mathcal{E}(g_b^{\cap})}}{\bigl(p_{\star}({t_k})\,p_{\star}({t_l})\bigr)^{\mathcal{E}(g_b^{\cap})}}.\notag
\end{align}

 {The goal of the last step} is to calculate the cardinality $|{\tt G}^{\star\star}_{\vk f}({H},\vk g^{\cap})|$, {which appears in \eqref{moment_calculation_8}. We state this result in the following lemma, proven in the appendix.}

{\begin{lem}\label{lem:G_star_star_cardinality} For any ${H}\in{\tt M}_z$ and any $\vk g^{\cap}\in{\tt OCS}_{\vk f}({H})$,
    \begin{align*}
        \abs{{\tt G}^{\star\star}_{\vk f}({H},\vk g^{\cap})} = \left(\sum_{b=1}^{z}\mathcal{V}({H}_{f_b}) - \sum_{b=1}^{z/2}\mathcal{V}(g_b^{\cap})\right)!\prod_{b=1}^{z/2}\mathcal{C}({H}_{f_{i_b}},{H}_{f_{j_b}},g^{\cap}_{b}),
    \end{align*}
    where the constants $\mathcal{C}({H}_{f_{i_b}},{H}_{f_{j_b}},g^{\cap}_{b})$ are as defined in \eqref{C_def}.
\end{lem}}

\noindent
{{\it Continuation of proof of \nelem{lem:asymptotic_normality}}.} Combining \eqref{moment_calculation_8} with {the result obtained in} \nelem{lem:G_star_star_cardinality}, we {can now derive the asymptotics of the inner sum of \eqref{moment_calculation_2}:} 
\begin{align}
   & {\lim_{N\to\infty}}\sum_{\vk g\in{\tt G}^{\star}_{\vk f,N}}\E{\prod_{i=1}^{z}\frac{\mathcal{F}_{{N}}({H}_{f_i})\xi_{f_i,{N}}(g_i)}{\sqrt{\mathcal{F}^{\rm opt}_N({H}_{f_i})}}}\notag\\
    &=\sum_{{H}\in{\tt M}_z} \left(\prod_{b=1}^{z/2}\mathbb{I}\{{H}_{f_{i_b}}\cong {H}_{f_{j_b}})\}\right)
    \times\sum_{\vk g^{\cap}\in{\tt OCS}_{\vk f}({H})}\prod_{b=1}^{z/2}\mathcal{C}({H}_{f_{i_b}},{H}_{f_{j_b}},g^{\cap}_{b})\sum_{k,l=1}^{n}\frac{q_{k,i}q_{l,j}\bigl(\kappa_{\star}(k,l)\bigr)^{\mathcal{E}(g_b^{\cap})}}{\bigl(p_{\star}(k)p_{\star}(l)\bigr)^{\mathcal{E}(g_b^{\cap})}}\notag\\
    &= \sum_{{H}\in{\tt M}_z} \left(\prod_{b=1}^{z/2}\mathbb{I}\{{H}_{f_{i_b}}\cong {H}_{f_{j_b}})\}\right) \\
    &\qquad \times\prod_{b=1}^{z/2}\left(\sum_{k,l=1}^{n}q_{k,i}q_{l.j}\sum_{g^{\cap}\in{\tt OCS}({H}_{f_{i_b}},{H}_{f_{j_b}})}\frac{\mathcal{C}({H}_{f_{i_b}},{H}_{f_{j_b}},g^{\cap})\bigl(\kappa_{\star}(k,l)\bigr)^{\mathcal{E}(g^{\cap})}}{\bigl(p_{\star}(k)p_{\star}(l)\bigr)^{\mathcal{E}(g^{\cap})}}\right)\notag\\
    &=\sum_{{H}\in{\tt M}_z}\prod_{b=1}^{z/2}\left(\sum_{k,l=1}^{n}q_{k,i}q_{l.j}\left(\Sigma_{k,l}\right)_{f_{i_b},f_{j_b}}\right),\label{moment_calculation_9}
\end{align}
where the block matrix $\vk{\Sigma}$ is defined in \eqref{sigma_klij_def}.
{By combining \eqref{moment_calculation_9} with \eqref{moment_calculation_2}, we obtain the asymptotics of the moments $\E{\xi_{{N}}^z}$: with the random variable $\xi_{{N}}$ as defined in \eqref{xi_init_def},}
\begin{align*}
    \lim_{N\to\infty}\E{\xi_{{N}}^z} &= \sum_{f_1,\ldots,f_{z}\in\sprod{m}}\sum_{{H}\in{\tt M}_z}\prod_{b=1}^{z/2}\left(\sum_{k,l=1}^{n}q_{k,i}q_{l.j}\left(\Sigma_{k,l}\right)_{f_{i_b},f_{j_b}}\right)\\
    &= \sum_{{H}\in{\tt M}_z}\sum_{f_1,\ldots,f_{z}\in\sprod{m}}\prod_{b=1}^{z/2}\left(\sum_{k,l=1}^{n}q_{k,i}q_{l.j}\left(\Sigma_{k,l}\right)_{f_{i_b},f_{j_b}}\right)\\
    &=\sum_{{H}\in{\tt M}_z}\left(\sum_{i,j=1}^{m}\sum_{k,l=1}^{n}q_{k,i}q_{l,j}\left(\Sigma_{k,l}\right)_{i,j}\right)^{z/2}\\
    &=\abs{{\tt M}_z}(\sigma^2)^{z/2}=(z-1)!!(\sigma^2)^{z/2},
\end{align*}
where $\sigma^2$ is defined in \eqref{sigma_def}. Hence, the claim \eqref{normality_claim} follows for even $z$ as well, which finishes the proof.
\QED

\section{Proof of \nelem{lem:tightness}}\label{section:proof_of_tightness}
The objective of this section is to establish {inequality \eqref{tightness_claim}, which due to \cite[Theorem 13.5]{billingsley2013convergence} leads to tightness.}  Using that
\begin{align*}
    \abs{\vk X^{\star}_N(t)-\vk X^{\star}_N(s)}^2 = \sum_{i=1}^{m}\abs{X^{\star}_{N,i}(t)-X^{\star}_{N,i}(s)}^2,
\end{align*}
{to establish \eqref{tightness_claim}, it suffices to prove that, for all $i,j\in\sprod{m}$,}
\begin{align}
    \Delta_{N,i,j}(r,s,t) \leqslant \abs{F_{i,j}(t)-F_{i,j}(r)}^2,\label{tightness_separate_claim}
\end{align}
where $F_{i,j}(t) =t\sqrt{51F_{i,j}^{\star}C_{i,j}^{\star}}$
and {for $0\leqslant r\leqslant s\leqslant t\leqslant T$,}
\begin{align}
    \Delta_{N,i,j}(r,s,t) = \E{\abs{X^{\star}_{N,i}(t)- X^{\star}_{N,i}(s)}^2\abs{X^{\star}_{N,j}(s)-X^{\star}_{N,j}(r)}^2}.\label{Delta_def}
\end{align}

{We begin by expressing the function $\Delta_{N,i,j}(r,s,t)$ from \eqref{Delta_def} in a form similar to \eqref{X_def_with_g}. By applying \eqref{X_def_with_g} and \eqref{X_star_def}, we can decompose the differences $X^{\star}_{N,i}(t) - X^{\star}_{N,i}(s)$ into a sum of indicator processes over all possible inclusions $g \in {\tt G}_N(H_i)$:}
\begin{align*}
    X^{\star}_{N,i}(t)-X^{\star}_{N,i}(s) = \frac{\mathcal{F}_N({H}_i)}{\sqrt{\mathcal{F}^{\rm opt}_N({H}_i)}}\sum_{g\in{\tt G}_N({H}_i)}\eta_{i,g,N}(s,t),
\end{align*}
{where}
\begin{align}
    \eta_{i,g,N}(s,t) &= \mathbb{I}\{g{\,{\subset}\,} G_N(t)\} - p^{\mathcal{E}({H}_i)}(t) - \mathbb{I}\{g{\,{\subset}\,} G_N(s)\} + p^{\mathcal{E}({H}_i)}(s).\label{nu_def}
\end{align}
{Hence, \eqref{Delta_def} has the following form}
\begin{align}
    \Delta_{N,i,j}(r,s,t) &= \frac{\mathcal{F}_N^2({H}_i)\mathcal{F}_N^2({H}_j)}{\mathcal{F}^{\rm opt}_N({H}_i)\mathcal{F}^{\rm opt}_N({H}_j)}\sum_{\substack{g_1,g_2\in {\tt G}_N({H}_i)\\ g_3,g_4\in {\tt G}_N({H}_j)}} \E{\eta_{i,g_1,N}(s,t)\eta_{i,g_2,N}(s,t)\eta_{j,g_3,N}(r,s)\eta_{j,g_4,N}(r,s)}\notag\\
    &= \frac{\mathcal{F}_N^2({H}_i)\mathcal{F}_N^2({H}_j)}{\mathcal{F}^{\rm opt}_N({H}_i)\mathcal{F}^{\rm opt}_N({H}_j)}\sum_{\substack{g_1,g_2\in {\tt G}_N({H}_i)\\ g_3,g_4\in {\tt G}_N({H}_j)}} \E{\prod_{k=1}^{4}\eta_{i_k,g_k,N}(x_k,y_k)},\label{tightness_sum}
\end{align}
where
\begin{align*}
    x_k = \begin{cases} s,\qquad k\in\{1,2\},\\
        r,\qquad k\in\{3,4\},
    \end{cases}\qquad 
    y_k = \begin{cases} t,\qquad k\in\{1,2\},\\
        s,\qquad k\in\{3,4\}{,}
    \end{cases}\qquad
    i_k = \begin{cases} i,\qquad k\in\{1,2\},\\
        j,\qquad k\in\{3,4\}.
    \end{cases}
\end{align*}
{We now proceed to prove \eqref{tightness_separate_claim} using the representation \eqref{tightness_sum}. The first step is to distinguish the zero terms of \eqref{tightness_sum} from the non-zero ones. Next, we derive an upper bound for each non-zero term in \eqref{tightness_sum}, and finally we show that the sum of the derived bounds does not exceed the right-hand side of \eqref{tightness_separate_claim}.}

{We identify the non-zero terms in the sum appearing in \eqref{tightness_sum}.} From the definition \eqref{nu_def} we conclude that the random variables $\eta_{i,g,N}(s,t)$ are centered, meaning that for any $i\in\sprod{m}$, $g\in{\tt G}_N({H}_i)$ and any $s,t\in[0,T]$
\begin{align*}
    \E{\eta_{i,g,N}(s,t)} &= 0,
\end{align*}
and $\eta_{i,g_1,N}(s_1,t_1)$, $\eta_{j,g_2,N}(s_2,t_2)$ are independent for any $s_1,t_1,s_2,t_2\in[0,T]$ if $\mathcal{E}(g_1\cap g_2) =0$.
Notice that if there exist $l\in\{1,2,3,4\}$ such that for every $k\in\{1,2,3,4\}$, ${l}\not= k$,
\begin{align*}
    \mathcal{E}(g_{{l}}\cap g_k) = 0,
\end{align*}
and in this case
\begin{align*}
    \E{\prod_{k=1}^{4}\eta_{i_k,g_k,N}(x_k,y_k)} = \E{\eta_{g_l}(x_l,y_l)}\E{\prod_{\substack{k=1 \\ k\not= l}}^{4}\eta_{i_k,g_k,N}(x_k,y_k)} = 0.
\end{align*}
So, in the sum \eqref{tightness_sum} we can consider only those sets $\vk g = (g_1,g_2,g_3,g_4)$ which belong to 
\begin{align}
    {\tt G}_{N}(i,j) = 
   \bigl\{\vk g\in {{\tt G}_N(H_i)\times{\tt G}_N(H_i)\times{\tt G}_N(H_j)\times{\tt G}_N(H_j)}\colon \forall k\in\sprod{4}\, \exists l\in\sprod{4}\setminus k\colon \mathcal{E}({g_k\cap g_l})\geqslant 1\bigr\}.
    \label{tightness_g_condition}
\end{align}
If $\vk g\not\in{\tt G}_{N}(i,j)$, the corresponding term of the sum \eqref{tightness_sum} equals zero.
We can then rewrite \eqref{tightness_sum} as follows
\begin{align}
    \Delta_{N,i,j}(r,s,t) = \frac{\mathcal{F}_N^2({H}_i)\mathcal{F}_N^2({H}_j)}{\mathcal{F}^{\rm opt}_N({H}_i)\mathcal{F}^{\rm opt}_N({H}_j)}\sum_{\vk g\in{{\tt G}_{N}(i,j)}} \E{\prod_{k=1}^{4}\eta_{i_k,g_k,N}(x_k,y_k)}{.}\label{tightness_new_sum}
\end{align}
{Next, we proceed to the second step, where we establish an upper bound on all non-zero terms in \eqref{tightness_new_sum} (equivalent to \eqref{tightness_sum}). We partition the probability space into disjoint events: for some of these events, we limit the value of the inner random variable, while for others we bound the probability of the event itself.} Define for $k\in\{1,\ldots,4\}$ the events
\begin{align}
    \Omega_k &= \bigl\{\mathbb{I}\{g_{k}{\,{\subset}\,} G_N(y_k)\} \not= \mathbb{I}\{g_{k}{\,{\subset}\,} G_N(x_k)\}\bigr\}.\label{Omega_def}
\end{align}
Then, each expectation in \eqref{tightness_new_sum} can be represented as a sum of expectations over all possible intersections of events $\Omega_k$ and $\Omega_k^{\rm c}$, for $k\in\{1,\ldots,4\}$.  
On each event $\Omega_k^{\rm c}$, immediately from \eqref{p_star_bound} and \eqref{p_inc_bound}, we have the upper bound
\begin{align}
    \abs{\eta_{i_k,g_k,N}(x_k,y_k)} &= \abs{p_N^{\mathcal{E}({H}_{i_k})}(y_k) - p_N^{\mathcal{E}({H}_{i_k})}(x_k)}\notag\\
    &=\abs{p_N(y_k) - p_N(x_k)}\left(\sum_{l=0}^{\mathcal{E}({H}_{i_k}) - 1} p_N^{l}(x_k)p_N^{\mathcal{E}({H}_{i_k}) - 1 - l}(y_k) \right)\notag\\
    &\leqslant \varrho^{\mathcal{E}({H}_{i_k})}\mathcal{E}({H}_{i_k})p_+^{\mathcal{E}({H}_{i_k})-1}{\mathfrak{C}}\,(y_k-x_k)\leqslant \varrho_N^{\mathcal{E}({H}_{i_k})}\mathcal{E}({H}_{i_k})p_+^{\mathcal{E}({H}_{i_k})-1}{\mathfrak{C}}\,(t-r)\label{O_c_bound}
\end{align}
{for} each event $\Omega_k$ we bound the corresponding probability, using the following lemma, which we prove in the appendix.

\begin{lem}\label{lem:omega_bound}
     For the events $\Omega_k$ defined in \eqref{Omega_def}, and for any permutation $\{k_1,k_2,k_3,k_4\} = \sprod{4}$,
    \begin{align*}
        \pk{\Omega_{k_1}}&\leqslant C_0\varrho_N^{\mathcal{E}(g_{k_1})}\abs{t-r},\\
        \pk{\Omega_{k_1}\Omega_{k_2}}&\leqslant C_0\varrho_N^{\mathcal{E}(g_{k_1}\cup g_{k_2})}\abs{t-r},\\
        \pk{\Omega_{k_1}\Omega_{k_2}\Omega_{k_3}}&\leqslant C_0\varrho_N^{\mathcal{E}(g_{k_1}\cup g_{k_2}\cup g_{k_3})}\abs{t-r}^2,\\
        \pk{\Omega_{k_1}\Omega_{k_2}\Omega_{k_3}\Omega_{k_4}}&\leqslant C_0\varrho_N^{\mathcal{E}(g_{1}\cup g_{2}\cup g_{3}\cup g_{4})}\abs{t-r}^2
    \end{align*}
    for the constant $C_0$ defined in \eqref{C_0_def}.
\end{lem}

\noindent {\it Continuation of proof of \nelem{lem:tightness}.}
Combining the bounds presented in \nelem{lem:omega_bound} with \eqref{O_c_bound} and the evident bound $\eta_{i_k,g_k,N}(x_k,y_k)\leqslant 2$,  we obtain that
\begin{align}
    &\hspace{-1.2cm}\E{\prod_{k=1}^{4}\eta_{i_k,g_k,N}(x_k,y_k)}\notag\\\leqslant &\:\pk{\Omega_1^{c}\Omega_2^{c}\Omega_3^{c}\Omega_4^{c}}\prod_{k=1}^{4}\left(\varrho_N^{\mathcal{E}({H}_{i_k})}\mathcal{E}({H}_{i_k})p_+^{\mathcal{E}({H}_{i_k})-1}{\mathfrak{C}}\abs{t-r}\right)\:+\notag\\
    &\frac{2}{3!}\sum_{\{k_1,k_2,k_3,k_4\} = \sprod{4}} \pk{\Omega_{k_1}\Omega_{k_2}^{c}\Omega_{k_3}^{c}\Omega_{k_4}^{c}}\prod_{q=2}^{4}\left(\varrho_N^{\mathcal{E}({H}_{i_{k_q}})}\mathcal{E}({H}_{i_{k_q}})p_+^{\mathcal{E}({H}_{i_{k_q}})-1}{\mathfrak{C}}\abs{t-r}\right)\:+\notag\\
    &\frac{4}{2!2!}\sum_{\{k_1,k_2,k_3,k_4\} = \sprod{4}} \pk{\Omega_{k_1}\Omega_{k_2}\Omega_{k_3}^{c}\Omega_{k_4}^{c}}\prod_{q=3}^{4}\left(\varrho_N^{\mathcal{E}({H}_{i_{k_q}})}\mathcal{E}({H}_{i_{k_q}})p_+^{\mathcal{E}({H}_{i_{k_q}})-1}{\mathfrak{C}}\abs{t-r}\right)\:+\notag\\
    &\frac{8}{3!}\sum_{\{k_1,k_2,k_3,k_4\} = \sprod{4}} \pk{\Omega_{k_1}\Omega_{k_2}\Omega_{k_3}\Omega_{k_4}^{c}}\varrho_N^{\mathcal{E}({H}_{i_{k_4}})}\mathcal{E}({H}_{i_{k_4}})p_+^{\mathcal{E}({H}_{i_{k_4}})-1}{\mathfrak{C}}\abs{t-r}\:+\notag\\
    &16\,\pk{\Omega_1\Omega_2\Omega_3\Omega_4}\notag\\
    \leqslant&\: \prod_{k=1}^{4}\left(\varrho_N^{\mathcal{E}({H}_{i_k})}\mathcal{E}({H}_{i_k})p_+^{\mathcal{E}({H}_{i_k})-1}{\mathfrak{C}}\abs{t-r}\right)\:+\label{tightness_expectation_bound_0}\\
    &\frac{2}{3!}\sum_{\{k_1,k_2,k_3,k_4\} = \sprod{4}} C_0\varrho_N^{\mathcal{E}(g_{k_1})}\abs{t-r}\prod_{q=2}^{4}\left(\varrho_N^{\mathcal{E}({H}_{i_{k_q}})}\mathcal{E}({H}_{i_{k_q}})p_+^{\mathcal{E}({H}_{i_{k_q}})-1}{\mathfrak{C}}\abs{t-r}\right)\:+\notag\\
    &\frac{4}{2!2!}\sum_{\{k_1,k_2,k_3,k_4\} = \sprod{4}} C_0\varrho_N^{\mathcal{E}(g_{k_1}\cup g_{k_2})}\abs{t-r}\prod_{q=3}^{4}\left(\varrho_N^{\mathcal{E}({H}_{i_{k_q}})}\mathcal{E}({H}_{i_{k_q}})p_+^{\mathcal{E}({H}_{i_{k_q}})-1}{\mathfrak{C}}\abs{t-r}\right)\:+\notag\\
    &\frac{8}{3!}\sum_{\{k_1,k_2,k_3,k_4\} = \sprod{4}} C_0\varrho_N^{\mathcal{E}(g_{k_1}\cup g_{k_2}\cup g_{k_3})}\abs{t-r}^2\varrho_N^{\mathcal{E}({H}_{i_{k_4}})}\mathcal{E}({H}_{i_{k_4}})p_+^{\mathcal{E}({H}_{i_{k_4}})-1}{\mathfrak{C}}\abs{t-r}\:+\notag\\
    &16\,C_0\varrho_N^{\mathcal{E}(g_{1}\cup g_{2}\cup g_{3}\cup g_{4})}\abs{t-r}^2{.}\notag
\end{align}
{We bound  each term {in the} sum \eqref{tightness_expectation_bound_0} separately.

\noindent {\it --~First term:}
\begin{align}
    \prod_{k=1}^{4}\left(\varrho_N^{\mathcal{E}({H}_{i_k})}\mathcal{E}({H}_{i_k})p_+^{\mathcal{E}({H}_{i_k})-1}{\mathfrak{C}}\abs{t-r}\right)\leqslant C_1F_1(N)\abs{t-r}^2{,}\label{CF_bound_1}
\end{align}
where
\begin{align*}
    C_1 &={\mathfrak{C}^4}\,T^2\,\mathcal{E}^2({H}_i)\,\mathcal{E}^{2}({H}_j)\,p_+^{2\mathcal{E}({H}_{i})+2\mathcal{E}({H}_j)-4},\\
    F_1(N) &=\varrho_N^{2\mathcal{E}({H_i})+2\mathcal{E}({H_j})}.
\end{align*}
\noindent {\it --~Second term:}
\begin{align}
    &\frac{2}{3!}\sum_{\{k_1,k_2,k_3,k_4\} = \sprod{4}} C_0\varrho_N^{\mathcal{E}(g_{k_1})}\abs{t-r}\prod_{q=2}^{4}\left(\varrho_N^{\mathcal{E}({H}_{i_{k_q}})}\mathcal{E}({H}_{i_{k_q}})p_+^{\mathcal{E}({H}_{i_{k_q}})-1}{\mathfrak{C}}\abs{t-r}\right) \notag\\
    &\qquad\leqslant C_2F_2(N)\abs{t-r}^2{,}\label{CF_bound_2}
\end{align}
where
\begin{align*}
    C_2 &=\frac{1}{3}\,C_0\,{\mathfrak{C}^3}\,T^2\max_{\{k_1,k_2,k_3,k_4\} = \sprod{4}} \prod_{q=2}^{4}\left(\mathcal{E}({H}_{i_{k_q}})p_+^{\mathcal{E}({H}_{i_{k_q}})-1}\right),\\
    F_2(N) &=\varrho_N^{2\mathcal{E}({H_i})+2\mathcal{E}({H_j})}.
\end{align*}
\noindent {\it --~Third term:}
\begin{align}
    &\frac{4}{2!2!}\sum_{\{k_1,k_2,k_3,k_4\}= \sprod{4}} C_0\varrho_N^{\mathcal{E}(g_{k_1}\cup g_{k_2})}\abs{t-r}\prod_{q=3}^{4}\left(\varrho_N^{\mathcal{E}({H}_{i_{k_q}})}\mathcal{E}({H}_{i_{k_q}})p_+^{\mathcal{E}({H}_{i_{k_q}})-1}{\mathfrak{C}}\abs{t-r}\right)\notag\\
    &\qquad \leqslant 4!\,C_3F_3(N)\abs{t-r}^2{,}\label{CF_bound_3}
\end{align}
where
\begin{align*}
    C_3 &= C_0\,{\mathfrak{C}^2}\,T\max_{\{k_1,k_2,k_3,k_4\} = \sprod{4}} \prod_{q=3}^{4}\left(\mathcal{E}({H}_{i_{k_q}})p_+^{\mathcal{E}({H}_{i_{k_q}})-1}\right),\\
    F_3(N) &=\frac{1}{4!}\sum_{\{k_1,k_2,k_3,k_4\}=\sprod{4}}\varrho_N^{\mathcal{E}(g_{k_1}\cup g_{k_2}) + \mathcal{E}(g_{k_3})+\mathcal{E}(g_{k_4})}.
\end{align*}
\noindent {\it --~Fourth term:}
\begin{align}
    &\frac{8}{3!}\sum_{\{k_1,k_2,k_3,k_4\} = \sprod{4}} C_0\varrho_N^{\mathcal{E}(g_{k_1}\cup g_{k_2}\cup g_{k_3})}\abs{t-r}^2\varrho_N^{\mathcal{E}({H}_{i_{k_4}})}\mathcal{E}({H}_{i_{k_4}})p_+^{\mathcal{E}({H}_{i_{k_4}})-1}{\mathfrak{C}}\abs{t-r}\notag\\
    &\qquad\leqslant 4!\,C_4F_4(N)\abs{t-r}^2{,}\label{CF_bound_4}
\end{align}
where
\begin{align*}
     C_4 &= \frac{4}{3}C_0\,{\mathfrak{C}}\,T\max_{\{k_1,k_2,k_3,k_4\} = \sprod{4}} \mathcal{E}({H}_{i_{k_4}})\,p_+^{\mathcal{E}({H}_{i_{k_4}})-1} ,\\
     F_4(N) &=\frac{1}{4!}\sum_{\{k_1,k_2,k_3,k_4\}=\sprod{4}}\varrho_N^{\mathcal{E}(g_{k_1}\cup g_{k_2}\cup g_{k_3}) + \mathcal{E}(g_{k_4})}.
\end{align*}
\noindent {\it --~Fifth term:}
\begin{align}
    16\,C_0\varrho_N^{\mathcal{E}(g_{1}\cup g_{2}\cup g_{3}\cup g_{4})}\abs{t-r}^2\leqslant C_5F_5(N)\abs{t-r}^2\label{CF_bound_5}
\end{align}
where
\begin{align*}
    C_5 &=16\,C_0,\\
     F_5(N) &= \varrho_N^{\mathcal{E}(g_1\cup g_2\cup g_3\cup g_4)}.
\end{align*}
Combining the bounds \eqref{CF_bound_1}--\eqref{CF_bound_5} with \eqref{tightness_expectation_bound_0}, we obtain that, with $\phi_1=\phi_2=\phi_5=1$ and $\phi_3=\phi_4=4!$,
\begin{align}
    \E{\prod_{k=1}^{4}\eta_{i_k,g_k,N}(x_k,y_k)}\leqslant&\: \abs{t-r}^2 \sum_{i=1}^5\phi_i \,C_i\,F_i(N).\label{tightness_expectation_bound}
\end{align}}

\COM{where the constants $C_1,\ldots,C_5$ are defined as
\begin{align*}
    C_1 &={\mathfrak{C}^4}\,T^2\,\mathcal{E}^2({H}_i)\,\mathcal{E}^{2}({H}_j)\,p_+^{2\mathcal{E}({H}_{i})+2\mathcal{E}({H}_j)-4},\\
    C_2 &=\frac{1}{3}\,C_0\,{\mathfrak{C}^3}\,T^2\max_{\{k_1,k_2,k_3,k_4\} = \sprod{4}} \prod_{q=2}^{4}\left(\mathcal{E}({H}_{i_{k_q}})p_+^{\mathcal{E}({H}_{i_{k_q}})-1}\right),\\
    C_3 &= C_0\,{\mathfrak{C}^2}\,T\max_{\{k_1,k_2,k_3,k_4\} = \sprod{4}} \prod_{q=3}^{4}\left(\mathcal{E}({H}_{i_{k_q}})p_+^{\mathcal{E}({H}_{i_{k_q}})-1}\right),\\
    C_4 &= \frac{4}{3}C_0\,{\mathfrak{C}}\,T\max_{\{k_1,k_2,k_3,k_4\} = \sprod{4}} \mathcal{E}({H}_{i_{k_4}})\,p_+^{\mathcal{E}({H}_{i_{k_4}})-1} ,\\
    C_5 &=16\,C_0,
\end{align*}
and the functions $F_1(\cdot),\ldots F_5(\cdot)$ as
\begin{align*}
    F_1(N) &=\varrho_N^{2\mathcal{E}({H_i})+2\mathcal{E}({H_j})},\\
    F_2(N) &=\varrho_N^{2\mathcal{E}({H_i})+2\mathcal{E}({H_j})},\\
    F_3(N) &=\frac{1}{4!}\sum_{\{k_1,k_2,k_3,k_4\}=\sprod{4}}\varrho_N^{\mathcal{E}(g_{k_1}\cup g_{k_2}) + \mathcal{E}(g_{k_3})+\mathcal{E}(g_{k_4})},\\
    F_4(N) &=\frac{1}{4!}\sum_{\{k_1,k_2,k_3,k_4\}=\sprod{4}}\varrho_N^{\mathcal{E}(g_{k_1}\cup g_{k_2}\cup g_{k_3}) + \mathcal{E}(g_{k_4})},\\
    F_5(N) &= \varrho_N^{\mathcal{E}(g_1\cup g_2\cup g_3\cup g_4)}.
\end{align*}}

{At this point we have upper bounds \eqref{tightness_expectation_bound} for each term of \eqref{tightness_new_sum}.}
This means that, in order to verify \eqref{tightness_separate_claim}, {and to conclude the proof of \nelem{lem:tightness},} it is enough to show that for every $k\in\sprod{5}$
\begin{align}
    \frac{\mathcal{F}_N^2({H}_i)\mathcal{F}_N^2({H}_j)}{\mathcal{F}^{\rm opt}_N({H}_i)\mathcal{F}^{\rm opt}_N({H}_j)}\sum_{\vk g\in{\tt G}_{N}(i,j)} F_k &< F^*_{i,j},\label{cal_F_bound}\\
    C_k&\leqslant C^{\star}_{i,j},\label{cal_C_bound}
\end{align}
where $F^*_{i,j}$ and $C^*_{i,j}$ are defined in \eqref{frak_f_def} and $\eqref{frak_c_def}$, respectively. The bound \eqref{cal_C_bound} is trivial, {as $\mathfrak{C}\geqslant 1$, $p_{+}\geqslant 1$, $C_0\geqslant 1$ and $\mathcal{E}(H_{i_k})\geqslant 1$ for any $k\in\sprod{4}$ implying
\begin{align*}
    C_1 &= {\mathfrak{C}^4}\,T^2\,\mathcal{E}^2({H}_i)\,\mathcal{E}^{2}({H}_j)\,p_+^{2\mathcal{E}({H}_{i})+2\mathcal{E}({H}_j)-4} \leqslant C^{\star}_{i,j}, \\
    C_2 &= \frac{1}{3}\,C_0\,{\mathfrak{C}^3}\,T^2\max_{\{k_1,k_2,k_3,k_4\} = \sprod{4}} \prod_{q=2}^{4}\left(\mathcal{E}({H}_{i_{k_q}})p_+^{\mathcal{E}({H}_{i_{k_q}})-1}\right)\\
    &\leqslant \frac{1}{3}\,C_0\,{\mathfrak{C}^3}\,T^2\max_{\{k_1,k_2,k_3,k_4\} = \sprod{4}} \prod_{q=1}^{4}\left(\mathcal{E}({H}_{i_{k_q}})p_+^{\mathcal{E}({H}_{i_{k_q}})-1}\right)\\
    &= \frac{1}{3}\,C_0\,{\mathfrak{C}^3}\,T^2\mathcal{E}^2({H}_i)\,\mathcal{E}^{2}({H}_j)\,p_+^{2\mathcal{E}({H}_{i})+2\mathcal{E}({H}_j)-4}\leqslant C_{i,j}^{\star},\\
    C_3 &= C_0\,{\mathfrak{C}^2}\,T\max_{\{k_1,k_2,k_3,k_4\} = \sprod{4}} \prod_{q=3}^{4}\left(\mathcal{E}({H}_{i_{k_q}})p_+^{\mathcal{E}({H}_{i_{k_q}})-1}\right)\\
    &\leqslant C_0\,{\mathfrak{C}^2}\,T\max_{\{k_1,k_2,k_3,k_4\} = \sprod{4}} \prod_{q=1}^{4}\left(\mathcal{E}({H}_{i_{k_q}})p_+^{\mathcal{E}({H}_{i_{k_q}})-1}\right)\\
    &= C_0\,{\mathfrak{C}^2}\,T\mathcal{E}^2({H}_i)\,\mathcal{E}^{2}({H}_j)\,p_+^{2\mathcal{E}({H}_{i})+2\mathcal{E}({H}_j)-4}\leqslant C_{i,j}^{\star},\\
    C_4 &= \frac{4}{3}C_0\,{\mathfrak{C}}\,T\max_{\{k_1,k_2,k_3,k_4\} = \sprod{4}} \mathcal{E}({H}_{i_{k_4}})\,p_+^{\mathcal{E}({H}_{i_{k_4}})-1}\\
    &\leqslant \frac{4}{3}C_0\,{\mathfrak{C}}\,T\max_{\{k_1,k_2,k_3,k_4\} = \sprod{4}}\prod_{q=1}^{4}\left(\mathcal{E}({H}_{i_{k_q}})p_+^{\mathcal{E}({H}_{i_{k_q}})-1}\right)\\
    &= \frac{4}{3}C_0\,{\mathfrak{C}}\,T\mathcal{E}^2({H}_i)\,\mathcal{E}^{2}({H}_j)\,p_+^{2\mathcal{E}({H}_{i})+2\mathcal{E}({H}_j)-4}\leqslant C_{i,j}^{\star},\\
    C_5 &= 16\,C_0 \leqslant C_{i,j}^{\star}.
\end{align*}
This entails that} we can consider only $\eqref{cal_F_bound}$. For any set $\vk g\in{\tt G}_{N}(i,j)$ and any permutation $\{k_1,k_2,k_3,k_4\}=\sprod{4}$ we {have that}
\begin{align*}
    \mathcal{E}(g_1\cup g_1\cup g_3\cup g_4)&\leqslant \mathcal{E}(g_{k_1}\cup g_{k_2}\cup g_{k_3}) + \mathcal{E}(g_{k_4})\leqslant \mathcal{E}(g_{k_1}\cup g_{k_2}) + \mathcal{E}(g_{k_3})+\mathcal{E}(g_{k_4}) \\&\leqslant 2\mathcal{E}({H_i})+2\mathcal{E}({H_j}){.}
\end{align*}
Hence, for $N$ large enough ensuring that $\varrho_N\leqslant 1$, it follows that
\begin{align*}
    {{F}_5(N)\geqslant {F}_4(N)\geqslant {F}_3(N)\geqslant {F}_2(N)= {F}_1(N),}
\end{align*}
so it is enough to check \eqref{cal_F_bound} only for $k=5$. Factorize the sum in \eqref{cal_F_bound} in terms of the number of vertices in the graph $g_{\cup}:=g_1\cup g_1\cup g_3\cup g_4$. Denoting $V^{\star} = 2\mathcal{V}({H}_i) + 2\mathcal{V}({H_j}) - 2$ and $V_{\star} = \max(\mathcal{V}({H}_i),\mathcal{V}({H}_j))$,
\begin{align}
    \frac{\mathcal{F}_N^2({H}_i)\mathcal{F}_N^2({H}_j)}{\mathcal{F}^{\rm opt}_N({H}_i)\mathcal{F}^{\rm opt}_N({H}_j)}\sum_{\vk g\in{\tt G}_{N}(i,j)} F_5 &= \frac{\mathcal{F}_N^2({H}_i)\mathcal{F}_N^2({H}_j)}{\mathcal{F}^{\rm opt}_N({H}_i)\mathcal{F}^{\rm opt}_N({H}_j)}\sum_{V=V_{\star}}^{V^{\star}}\sum_{\substack{\vk g\in{\tt G}_{N}(i,j) \\ \mathcal{V}(g_{\cup}) = V}} F_5\notag\\
    &\leqslant \frac{\mathcal{F}_N^2({H}_i)\mathcal{F}_N^2({H}_j)}{\mathcal{F}^{\rm opt}_N({H}_i)\mathcal{F}^{\rm opt}_N({H}_j)}\sum_{V=V_{\star}}^{V^{\star}}N^{V}\sum_{\substack{\vk g\in{\tt G}_{V}(i,j) \\ \mathcal{V}(g_{\cup}) = V}} F_5\notag\\
    &= \frac{\mathcal{F}_N^2({H}_i)\mathcal{F}_N^2({H}_j)}{\mathcal{F}^{\rm opt}_N({H}_i)\mathcal{F}^{\rm opt}_N({H}_j)}\sum_{V=V_{\star}}^{V^{\star}}\sum_{\substack{\vk g\in{\tt G}_{V}(i,j) \\ \mathcal{V}(g_{\cup}) = V}} \frac{1}{\mathcal{F}_N(g_{\cup})}{.}\label{tightness_final_sum}
\end{align}
Without loss of generality we may assume from \eqref{tightness_g_condition} that
\begin{align}
    \mathcal{E}(g_1\cap g_2)\geqslant 1,\qquad \mathcal{E}(g_3\cap g_4)\geqslant 1{.}\notag
\end{align}
Let us represent $\mathcal{F}_N(g_{\cup})$ by
\begin{align}
    \mathcal{F}_N(g_{\cup}) &= \mathcal{F}_N\bigl(g_1\cup g_2\cup g_3\bigr)\frac{\mathcal{F}_N(g_4)}{\mathcal{F}_N\bigl(\bigl(g_1\cup g_2\cup g_3\bigr)\cap g_4\bigr)}\notag\\
    &= \mathcal{F}_N\bigl(g_1\cup g_2\bigr)\frac{\mathcal{F}_N(g_3)}{\mathcal{F}_N\bigl(\bigl(g_1\cup g_2\bigr)\cap g_3\bigr)}\frac{\mathcal{F}_N(g_4)}{\mathcal{F}_N\bigl(\bigl(g_1\cup g_2\cup g_3\bigr)\cap g_4\bigr)}\notag\\
    &= \frac{\mathcal{F}_N(g_1)\mathcal{F}_{N}(g_2)}{\mathcal{F}_N(g_1\cap g_2)}\frac{\mathcal{F}_N(g_3)}{\mathcal{F}_N\bigl(\bigl(g_1\cup g_2\bigr)\cap g_3\bigr)}\frac{\mathcal{F}_N(g_4)}{\mathcal{F}_N\bigl(\bigl(g_1\cup g_2\cup g_3\bigr)\cap g_4\bigr)}{.}\label{F_cal_decomposition}
\end{align}
From \eqref{F_cal_zero_limit} we conclude that, for large enough $N$,
\begin{align*}
    \frac{1}{\mathcal{F}_N\bigl(\bigl(g_1\cup g_2\bigr)\cap g_3\bigr)}\geqslant 1.
\end{align*}
Hence, from \eqref{F_cal_decomposition} we see that, for large enough $N$, denoting {$g_{\star} = \bigl(g_1\cup g_2\cup g_3\bigr)\cap g_4$}
\begin{align}
    \mathcal{F}_N(g_{\cup})\geqslant  \frac{\prod_{k=1}^{4}\mathcal{F}_N(g_k)}{\mathcal{F}_N(g_1\cap g_2)\mathcal{F}_N(g_{\star})}{.}\label{F_cal_union_bound}
\end{align}
Applying the bound \eqref{F_cal_union_bound} to each term of \eqref{tightness_final_sum}, we obtain 
\begin{align}
     \frac{\mathcal{F}_N^2({H}_i)\mathcal{F}_N^2({H}_j)}{\mathcal{F}^{\rm opt}_N({H}_i)\mathcal{F}^{\rm opt}_N({H}_j)}\sum_{\vk g\in{\tt G}_{N}(i,j)} F_5&\leqslant \frac{\mathcal{F}_N^2({H}_i)\mathcal{F}_N^2({H}_j)}{\mathcal{F}^{\rm opt}_N({H}_i)\mathcal{F}^{\rm opt}_N({H}_j)}\sum_{V=V_{\star}}^{V^{\star}}\sum_{\substack{\vk g\in{\tt G}_{V}(i,j) \\ \mathcal{V}(g_{\cup}) = V}}\frac{\mathcal{F}_N(g_1\cap g_2)\mathcal{F}_N(g_{\star})}{\prod_{k=1}^{4}\mathcal{F}_N(g_k)}\\
     &= \sum_{V=V_{\star}}^{V^{\star}}\sum_{\substack{\vk g\in{\tt G}_{V}(i,j) \\ \mathcal{V}(g_{\cup}) = V}}\frac{\mathcal{F}_N(g_1\cap g_2)}{\mathcal{F}^{\rm opt}_N({H}_i)}\frac{\mathcal{F}_N(g_{\star})}{\mathcal{F}^{\rm opt}_N({H}_j)}{.}\label{F_cal_decomposition_2}
\end{align}
{Using \eqref{F_cal_opt_def} for {$g=\pi(g_1\cap g_2)\in {\tt CS}(H_i)$} we obtain for large enough $N$ that
\begin{align}
    \frac{\mathcal{F}_N(g_1\cap g_2)}{\mathcal{F}^{\rm opt}_N(H_i)} \leqslant 2,\label{2_bound_1}
\end{align}
and using \eqref{F_cal_opt_def} for {$g = \pi(g_{\star})\in {\tt CS}(H_j)$} we obtain for large enough $N$ that
\begin{align}
    \frac{\mathcal{F}_N(g_{\star})}{\mathcal{F}^{\rm opt}_N(H_j)} \leqslant 2.\label{2_bound_2}
\end{align}
Applying \eqref{2_bound_1} and \eqref{2_bound_2} for any of the individual terms of \eqref{F_cal_decomposition_2}, we derive that 
\begin{align*}
     \frac{\mathcal{F}_N^2(H_i)\mathcal{F}_N^2(H_j)}{\mathcal{F}^{\rm opt}_N(H_i)\mathcal{F}^{\rm opt}_N(H_j)}\sum_{\vk g\in{\tt G}_{N}(i,j)} F_5&\leqslant \sum_{V=V_{\star}}^{V^{\star}}\sum_{\substack{\vk g\in{\tt G}_{V}(i,j) \\ \mathcal{V}(g_{\cup}) = V}}2\cdot 2\\
     &\leqslant \sum_{V=V_{\star}}^{V^{\star}}\bigg(\sum_{\substack{g_1,g_2\in {\tt G}_V(H_i)\\ g_3,g_4\in {\tt G}_V(H_j)}}4\bigg)\,=\, \sum_{V=V_{\star}}^{V^{\star}}4\abs{{\tt G}_V(H_i)}^2\abs{{\tt G}_V(H_j)}^2\\
     &= 4\sum_{V=V_{\star}}^{V^{\star}}\left(\frac{V!}{(V-\mathcal{V}(H_i))!\mathcal{A}(H_i)}\right)^2\left(\frac{V!}{(V-\mathcal{V}(H_j))!\mathcal{A}(H_j)}\right)^2\\
     &=F^{\star}_{i,j},
\end{align*}
}where $F^{\star}_{i,j}$ is as defined in \eqref{frak_f_def}. Hence \eqref{cal_F_bound} follows, thus establishing \eqref{tightness_separate_claim} {and finalizing the proof}.
\QED

\appendix

\section{{Auxiliary proofs}}\label{section:appendix}

\textit{Proof of \neprop{example}.} 
We first show that in the setting considered the expected number of switching points of the alternating renewal process is finite. Denote for any given $N\in\N$ and starting state $x\in\{0,1\}$ the $i$-th (for $i\in\N$) switching point of the alternating renewal process {$a_N(\cdot)$} by
\begin{align*}
    \tau_{N,x,i} = \sum_{j=0}^{i}\xi_{N,x,j}. 
\end{align*}
Then, for any $x\in\{0,1\}$ and $N\in\N$,
\begin{align}
    \E{\max\{n\colon \tau_{N,x,n}\leqslant T\}} & = \sum_{i=1}^{\infty}\pk{\tau_{N,x,i}\leqslant T}\notag\\ & = \sum_{i=1}^{\infty}\pk{\sum_{j=1}^{i}\xi_{N,x,j} \leqslant T} = \sum_{i=1}^{\infty}\int_{0}^{T}\rho_{N, x} ^{\ast i}(t)\td t,\label{example_exp_bound}
\end{align}
where (with $\ast$ denoting the convolution of two functions)
$\rho_{N, x}^{\ast i} = \rho_{N, x,1}\ast\rho_{N, x,2}\ast\ldots\ast\rho_{N, x,n}.$
 Consider first $x=1$. By induction in $k$, it immediately follows that, for any $k\in\N$ and $t\in[0,T]$,
\begin{align}
    \rho_{N, 1}^{\ast k}(t)\leqslant \frac{P^{k}t^{k-1}}{(k-1)!}.\label{convolution_bound}
\end{align}
{Indeed, for $k=1$ it is given, while for $k>1$ we have}
\begin{align*}
    \rho_{N, 1}^{\ast k}(t)=\int_{0}^{t}\rho_{N,1}^{\ast k-1}(y)\rho_{N,1,k}(t-y)\td y\leqslant P\int_{0}^{t} \frac{P^{k-1}y^{k-2}}{(k-2)!}\td y= \frac{P^{k}t^{k-1}}{(k-1)!}.
\end{align*}
Inserting \eqref{convolution_bound} into \eqref{example_exp_bound} yields the upper bound
\begin{align}
    \E{\max\{n\colon \tau_{N,1,n}\leqslant T\}} \leqslant \sum_{i=1}^{\infty}\int_{0}^{T}\frac{P^{i}t^{i-1}}{(i-1)!}\td t = e^{PT}-1{.}\label{example_exp_bound2}
\end{align}
Similarly, for $x=0$ we obtain that
\begin{align}
    \E{\max\{n\colon \tau_{N,0,n}\leqslant T\}} \leqslant \varrho_N(e^{PT}-1){.} \label{example_exp_bound3}
\end{align}
We proceed by deriving the properties \eqref{claim_1} and \eqref{claim_2}. Fix some $r,s,t\in[0,T]$ such that $r\leqslant s\leqslant t$. {Define, for any random variable $\chi$, by $F_{\chi}(\cdot)$ the cumulative distribution function of this random variable.} For \eqref{claim_1} we observe that
\begin{align*}
    \pk{a_N(s)\not= a_N(r)} &\leqslant \pk{\exists n\in\N,\, x\in\{0,1\}\colon a_N(0)=x,\,\tau_{N,x,n}\in[r,s]}\\
    &\leqslant \sum_{{x\in\{0,1\}}}\pk{a_N(0)=x}\sum_{i=1}^{\infty} \pk{\tau_{N,x,i-1}\leqslant r,\, \tau_{N,x,i}\in[r,s]}\\
    &= \sum_{{x\in\{0,1\}}}\pk{a_N(0)=x}\sum_{i=1}^{\infty} \int_{0}^{r}\pk{\tau_{N_x,i}\in[r,s] \mid\tau_{N,x,i-1}=t_0}\td F_{\tau_{N,x,i-1}}(t_0)\\
    &= \sum_{{x\in\{0,1\}}}\pk{a_N(0)=x}\sum_{{i=0}}^{\infty} \int_{0}^{r}\pk{\xi_{N,x,{i+1}}\in[r-t_0,\, s-t_0]}\td F_{\tau_{N,x,{i}}}(t_0){.}
\end{align*}
{Now applying the second and third assumption of \neprop{example}, we conclude that}
\begin{align*}
     \pk{a_N(s)\not= a_N(r)}&\leqslant \pk{a_N(0)=1}\sum_{{i=0}}^{\infty} \int_{0}^{r}P\abs{s-r}\td F_{\tau_{N,1,{i}}}(t_0)\\
    &\qquad + \pk{a_N(0)=0}\left(P\varrho_N\abs{s-r} + \sum_{i=1}^{\infty} \int_{0}^{r}P\abs{s-r}\td F_{\tau_{N,0,i}}(t_0)\right).
\end{align*}
Furthermore, using the first assumption of \neprop{example} for the probability $\pk{a_N(0)=1}$, and combining \eqref{example_exp_bound2} and \eqref{example_exp_bound3}, we obtain for large enough $N$ that
\begin{align*}
    \pk{a_N(s)\not= a_N(r)}&\leqslant \varrho_N(p_{\star}(0)+1)\sum_{{i=0}}^{\infty}P\abs{s-r}\pk{\tau_{N,1,{i}}<T}\\
    &\qquad + \left(P\varrho_N\abs{s-r} + \sum_{i=1}^{\infty}P\abs{s-r}\pk{\tau_{N,0,i}<T}\right)\\
    &= \varrho_N(p_{\star}(0)+1)P\abs{s-r}\bigl(\E{\max\{n\colon \tau_{N,1,n}\leqslant T\}}+1\bigr)\\
    &\qquad + \left(P\varrho_N\abs{s-r} + P\abs{s-r}\E{\max\{n\colon \tau_{N,0,n}\leqslant T\}}\right)\\
    &\leqslant (p_{\star}(0)+2)Pe^{PT}\varrho_N\abs{s-r}.
\end{align*}
{We use a similar approach for \eqref{claim_2}. As $\xi_{N,x,i+1}$ is independent of $\tau_{N,x,i}$,
\begin{align*}
    &\hspace{-6mm}\pk{a_N(s)\not= a_N(r),\, a_N(t)\not= a_N(s)} \\&\leqslant \sum_{{x\in\{0,1\}}}\pk{a_N(0)=x}\pk{\exists n\in\N,\colon \tau_{N,x,n},\tau_{N,x,n+1}\in{[r,t]}}\\
    &\leqslant \sum_{{x\in\{0,1\}}}\pk{a_N(0)=x}\pk{\exists n\in\N:\tau_{N,x,n-1}\leqslant r,\,\tau_{N,x,n}\in [r,t],\, \xi_{N,x,n+1}<\abs{t-r}}\\
    &= \sum_{{x\in\{0,1\}}}\pk{a_N(0)=x}\pk{\exists n\in\N:\tau_{N,x,n-1}\leqslant r,\,\tau_{N,x,n}\in[r,t]}\pk{\xi_{N,x,n+1}<\abs{t-r}}.
\end{align*}
    Then applying the second assumption of \neprop{example} to the probability $\pk{\xi_{N,x,n+1}<\abs{t-r}}$, we have
\begin{align}
    &\pk{a_N(s)\not= a_N(r),\, a_N(t)\not= a_N(s)}\nonumber\\ 
    &\leqslant P\abs{t-r}\sum_{{x\in\{0,1\}}}\pk{a_N(0)=x}\pk{\exists n\in\N\colon\tau_{N,x,n-1}\leqslant r,\,\tau_{N,x,n}\in[t,r]}.\label{example_2.6_1}
\end{align}
Following the same calculations as the ones for \eqref{claim_1}, one observes that
\begin{align}
    \sum_{{x\in\{0,1\}}}\pk{a_N(0)=x}\pk{\exists n\in\N\colon\tau_{N,x,n-1}\leqslant r,\,\tau_{N,x,n}\in[t,r]}\leqslant (p_{\star}(0)+2)Pe^{PT}\varrho_N\abs{t-r}.\label{example_2.6_2}  
\end{align}
Combining \eqref{example_2.6_1} with \eqref{example_2.6_2} we obtain that
\begin{align*}
     \pk{a_N(s)\not= a_N(r),\, a_N(t)\not= a_N(s)}&\leqslant (p_{\star}(0)+2)P^2e^{PT}\varrho_N\abs{t-r}^2.
\end{align*}
}
Hence, the claim follows.
\QED

\medskip

\noindent {\it Proof of Lemma \ref{lem:constant_formula}.}
{In order to keep our notation compact, in this proof we locally abbreviate $\mathcal{N}(H,H^{\star},g)$ to $\mathcal{N}$.} So as to calculate $\mathcal{C}({H},{H^{\star}},g)$, we first choose $\mathcal{V}(g)$ vertices which belong to both graphs $g_1$ and $g_2$, and then divide the remaining $\mathcal{N} - \mathcal{V}(g)$ vertices between $g_1$ and $g_2$. To do this, we introduce some new notation. We define, for any ${\tt S}\subset \sprod{\mathcal{N}}$, $V\leqslant \abs{{\tt S}}$,
\begin{align*}
    {\tt P}_{{\tt S}}(V) = \left\{\vk v\in {\tt S}^{V}: v_1<\ldots < v_{V}\right\}.
\end{align*}
Let also, for any set $\vk v = \{v_1,\ldots, v_n\}\subset \sprod{\mathcal{N}}$, the graph $K_{\vk v}\in K_{\mathcal{N}}$ denote the complete graph on the vertices $v_1,\ldots, v_n$.
Hence,
\begin{align}
   \mathcal{C}({H},{H^{\star}},g) &= \frac{1}{\mathcal{N}!}\sum_{\substack{g_1\in{\tt G}_{\mathcal{N}}({H}) \notag\\ g_2\in{\tt G}_{\mathcal{N}}({H^{\star}})}}\mathbb{I}\{{\pi(g_1\cap g_2)= g}\}\\
   &= \frac{1}{\mathcal{N}!}\sum_{\vk v\in {\tt P}_{\sprod{\mathcal{N}}}(\mathcal{V}(g))}\sum_{\substack{g^{\star}{\,{\subset}\,} K_{\vk v} \\ {\pi(g^{\star})= g}}}\sum_{\substack{g_1\in{\tt G}_{\mathcal{N}}({H}) \\ g_2\in{\tt G}_{\mathcal{N}}({H^{\star}})}}\mathbb{I}\{g_1\cap g_2= g^{\star}\}\notag\\
    &= \frac{1}{\mathcal{N}!}\sum_{\vk v\in {\tt P}_{\sprod{\mathcal{N}}}(\mathcal{V}(g))}\sum_{\vk v_1\in{\tt P}_{\sprod{N}\setminus\vk v}(\mathcal{V}({H})-\mathcal{V}(g))}\sum_{\substack{g^{\star}{\,{\subset}\,} K_{\vk v} \\ {\pi(g^{\star})= g}}}\sum_{\substack{g_1{\,{\subset}\,} K_{\vk v\cup \vk v_1} \\ g_1\sim {H}}}\sum_{\substack{g_2{\,{\subset}\,} K_{\sprod{\mathcal{N}}\setminus \vk v_1} \\ g_2\sim {H^{\star}}}}\mathbb{I}\{g_1\cap g_2= g^{\star}\}{.}\label{C_representation_1}
\end{align}
Consider first the inner indicator in \eqref{C_representation_1}:
\begin{align}
    \mathbb{I}\{g_1\cap g_2 = g^{\star}\} &= \mathbb{I}\{g^{\star}{\,{\subset}\,} g_1\}\mathbb{I}\{g^{\star}{\,{\subset}\,} g_2\}\mathbb{I}\{g_1\cap g_2 = g^{\star}\}\notag\\
    &= \mathbb{I}\{g^{\star}{\,{\subset}\,} g_1\cap K_{\vk v}\}\mathbb{I}\{g^{\star}{\,{\subset}\,} g_2\cap K_{\vk v}\}\mathbb{I}\{g_1\cap g_2 = g^{\star}\}.\label{indicator_1}
\end{align}
We would like to show that {for any {$g^{\star}\subset K_{\vk v}$ such that $\pi(g^{\star})\in{\tt OCS}(H,H^{\star})$}}
\begin{align}
    \mathbb{I}\{g^{\star}{\,{\subset}\,} g_1\cap K_{\vk v}\} = \mathbb{I}\{g^{\star}= g_1\cap K_{\vk v}\},\qquad \mathbb{I}\{g^{\star}{\,{\subset}\,} g_2\cap K_{\vk v}\} = \mathbb{I}\{g^{\star}= g_2\cap K_{\vk v}\}.\label{indicator_equiv}
\end{align}
{In other words, we show that if \eqref{indicator_equiv} is not satisfied, {the pattern of the graph $g^{\star}$ is not an optimal common subgraph pattern of $H$ and $H^{\star}$.} In case $\varrho_N=1$, we see from combining \eqref{OCS_def} and \eqref{F_claim_0} that the only optimal common subgraph possible is an edge, meaning that $g^{\star}$ is isomorphic to an edge. Hence, as $\abs{\vk v} = \mathcal{V}(g^{\star})=2$, graph $K_{\vk v}$ as well is isomorphic to an edge, thus \eqref{indicator_equiv} follows.} So, we need to consider only the case $\varrho_N\to 0$ as $N\to\infty$. As ${H}\cong {H^{\star}}$ and $g\in{\tt OCS}({H},{H^{\star}})$, applying \eqref{GCS_cong} we have that
\begin{align*}
    g\in{\tt OCS}({H}),\qquad  g\in{\tt OCS}({H^{\star}}).
\end{align*}
Assume that \eqref{indicator_equiv} does not hold for some $g_1$, $\vk v$ and $g^{\star}$, meaning that
\begin{align}
    g^{\star}{\,{\subset}\,} g_1\cap K_{\vk v},\qquad g^{\star}\not= g_1\cap K_{\vk v}.\label{wrong_ex}
\end{align}
Define $g^{\star\star} = g_1\cap K_{\vk v}$. It is clear that ${\pi(g^{\star\star})\in\,} {\tt CS}({H})$, $\mathcal{V}(g^{\star\star}) = \mathcal{V}(g^{\star}) = \mathcal{V}(g)$, but due to \eqref{wrong_ex} $\mathcal{E}(g^{\star\star}) > {\mathcal{E}(g^{\star}) = \mathcal{E}(g)}$. Hence,
\begin{align*}
    \lim_{N\to\infty}\frac{\mathcal{F}_N(g^{\star\star})}{\mathcal{F}_N(g)} = \lim_{N\to\infty}\frac{N^{-\mathcal{V}(g^{\star\star})}\varrho_N^{-\mathcal{E}(g^{\star\star})}}{N^{-\mathcal{V}(g)}\varrho_N^{-\mathcal{E}(g)}} = \varrho_N^{\mathcal{E}(g) - \mathcal{E}(g^{\star\star})} = \infty,
\end{align*}
which contradicts \eqref{opt_property} as ${\pi(g^{\star\star})}\in{\tt CS}({H})$ and $g\in{\tt OCS}({H})$.
This implies that \eqref{indicator_equiv} is met. Applying it to \eqref{indicator_1} we obtain that 
\begin{align}
    \mathbb{I}\{g_1\cap g_2 = g^{\star}\} &= \mathbb{I}\{g^{\star}= g_1\cap K_{\vk v}\}\mathbb{I}\{g^{\star}{\,{\subset}\,} g_2= K_{\vk v}\}\mathbb{I}\{g_1\cap g_2 = g^{\star}\}\notag\\
    &=\mathbb{I}\{g_1\cap K_{\vk v} = g^{\star}\} \mathbb{I}\{g_2\cap K_{\vk v} = g^{\star}\}.\label{indicator_2}
\end{align}
Using \eqref{indicator_2}, we observe by \eqref{C_representation_1} that 
\begin{align}
    \mathcal{C}({H},{H^{\star}},g) \:= \frac{1}{\mathcal{N}!}\sum_{\vk v\in {\tt P}_{\sprod{\mathcal{N}}}(\mathcal{V}(g))}\sum_{\vk v_1\in{\tt P}_{\sprod{N}\setminus\vk v}(\mathcal{V}({H})-\mathcal{V}(g))}\sum_{\substack{g^{\star}{\,{\subset}\,} K_{\vk v} \\ {\pi(g^{\star})= g}}}&\left(\sum_{\substack{g_1{\,{\subset}\,} K_{\vk v\cup \vk v_1} \\ g_1\sim {H}}}\mathbb{I}\{g_1\cap K_{\vk v} = g^{\star}\}\right)\notag\\
    \hspace{6.3cm}\times\,&\left(\sum_{\substack{g_2{\,{\subset}\,} K_{\sprod{\mathcal{N}}\setminus \vk v_1} \\ g_2\sim {H^{\star}}}}\mathbb{I}\{g_2\cap K_{\vk v} = g^{\star}\}\right){.}\label{C_representation_2}
\end{align}
For any particular choice of $\vk v$ and $\vk v_1$, by relabeling the vertices, we see that
\begin{align*}
    \sum_{\substack{g_1{\,{\subset}\,} K_{\vk v\cup \vk v_1} \\ g_1\sim {H}}}\mathbb{I}\{g_1\cap K_{\vk v} = g^{\star}\} = \frac{\mathcal{K}({H},g)}{\abs{{\tt G}_{\mathcal{V}(g)}(g)}}=\frac{\mathcal{K}({H},g)\mathcal{A}(g)}{\mathcal{V}(g)!},
\end{align*}
where
\begin{align*}
        \mathcal{K}({H},g) = \sum_{g_1\in{\tt G}_{\mathcal{V}({H})}({H})}\mathbb{I}\{{\pi(g_1\cap K_{\sprod{\mathcal{V}(g)}})= g}\}.
    \end{align*}
Similarly, we have that
\begin{align*}
    \sum_{\substack{g_2{\,{\subset}\,} K_{\sprod{\mathcal{N}}\setminus \vk v_1} \\ g_2\sim {H^{\star}}}}\mathbb{I}\{g_2\cap K_{\vk v} = g^{\star}\} = \frac{\mathcal{K}({H^{\star}},g)\mathcal{A}(g)}{\mathcal{V}(g)!}.
\end{align*}
and thus, from \eqref{C_representation_2} we obtain, {also using \eqref{G_cardinality}},
\begin{align}
    \mathcal{C}({H},{H^{\star}},g) &= \frac{\mathcal{K}({H},g)\mathcal{K}({H^{\star}},g)\mathcal{A}^2(g)}{\bigl(\mathcal{V}(g)!\bigr)^2\mathcal{N}!}\sum_{\vk v\in {\tt P}_{\sprod{\mathcal{N}}}(\mathcal{V}(g))}\sum_{\vk v_1\in{\tt P}_{\sprod{N}\setminus\vk v}(\mathcal{V}({H})-\mathcal{V}(g))}\sum_{\substack{g^{\star}{\,{\subset}\,} K_{\vk v} \\ {\pi(g^{\star})= g}}}1\notag\\
    &= \frac{\mathcal{K}({H},g)\mathcal{K}({H^{\star}},g)\mathcal{A}^2(g)}{\bigl(\mathcal{V}(g)!\bigr)^2\mathcal{N}!}\sum_{\vk v\in {\tt P}_{\sprod{\mathcal{N}}}(\mathcal{V}(g))}\sum_{\vk v_1\in{\tt P}_{\sprod{N}\setminus\vk v}(\mathcal{V}({H})-\mathcal{V}(g))}\abs{{\tt G}_{\mathcal{V}(g)}(g)}\notag\\
    &= \frac{\mathcal{K}({H},g)\mathcal{K}({H^{\star}},g)\mathcal{A}(g)}{\mathcal{V}(g)!\mathcal{N}!}\sum_{\vk v\in {\tt P}_{\sprod{\mathcal{N}}}(\mathcal{V}(g))}\sum_{\vk v_1\in{\tt P}_{\sprod{N}\setminus\vk v}(\mathcal{V}({H})-\mathcal{V}(g))}1\notag\\
    &= \frac{\mathcal{K}({H},g)\mathcal{K}({H^{\star}},g)\mathcal{A}(g)}{\mathcal{V}(g)!\mathcal{N}!}\frac{\mathcal{N}!}{\mathcal{V}(g)!\bigl(\mathcal{V}({H})-\mathcal{V}(g)\bigr)!\bigl(\mathcal{V}({H^{\star}})-\mathcal{V}(g)\bigr)!}\notag\\
    &=\left(\frac{\sqrt{\mathcal{A}(g)}\mathcal{K}({H},g)}{\mathcal{V}(g)!\bigr(\mathcal{V}({H})-\mathcal{V}(g)\bigl)!}\right)\left(\frac{\sqrt{\mathcal{A}(g)}\mathcal{K}({H^{\star}},g)}{\mathcal{V}(g)!\bigr(\mathcal{V}({H^{\star}})-\mathcal{V}(g)\bigl)!}\right).\label{C_in_K}
\end{align}

Finally, consider the constants $\mathcal{K}({H},g)$. Let us for all $g_i\in{\tt G}_{\mathcal{V}({H})}({H})$ calculate the total number of subgraphs which {have pattern} $g$. On one hand, this number can be represented using the constant $\mathcal{K}({H},g)$:
\begin{align*}
    \sum_{\vk v\in{\tt P}_{\sprod{\mathcal{V}({H})}}(\mathcal{V}(g))}\sum_{g_1\in{\tt G}_{\mathcal{V}({H})}({H})}\mathbb{I}\{{\pi(g_1\cap K_{\vk v})= g}\} = \abs{{\tt P}_{\sprod{\mathcal{V}({H})}}(\mathcal{V}(g))}\mathcal{K}({H},g).
\end{align*}
On the other hand, the same number can be represented using the constant $\mathcal{S}({H}_i,g)$, as defined in \eqref{S_cal_def}:
\begin{align*}
    \sum_{g_1\in{\tt G}_{\mathcal{V}({H})}({H})}\sum_{\vk v\in{\tt P}_{\sprod{\mathcal{V}({H})}}(\mathcal{V}(g))}\mathbb{I}\{{\pi(g_1\cap K_{\vk v})= g}\}&=\sum_{g_1\in{\tt G}_{\mathcal{V}({H})}({H})}\mathcal{S}({H},g)=\abs{{\tt G}_{\mathcal{V}({H})}({H})}\mathcal{S}({H},g).
\end{align*}
Combining these two representations, 
\begin{align}
    \mathcal{K}({H},g)&=\frac{\abs{{\tt G}_{\mathcal{V}({H})}({H})}}{\abs{{\tt P}_{\sprod{\mathcal{V}({H})}}(\mathcal{V}(g))}}\mathcal{S}({H},g)\notag\\
    &=\frac{\mathcal{V}({H})!}{\mathcal{A}({H})}\frac{\mathcal{V}(g)!\bigl(\mathcal{V}({H})-\mathcal{V}(g)\bigr)!}{\mathcal{V}({H})!}\mathcal{S}({H},g)=\frac{\mathcal{V}(g)!\bigl(\mathcal{V}({H})-\mathcal{V}(g)\bigr)!}{\mathcal{A}({H})}\mathcal{S}({H},g).\label{K_formula}
\end{align}
Hence, the claim follows combining \eqref{C_in_K} with \eqref{K_formula}.
\QED

\medskip

\noindent {\it Proof of Proposition \ref{prop:X_vk_cov}.}
Using the definition of process {$\vk X(\cdot)$} given in \netheo{main},
\begin{align*}
    \operatorname{Cov}(X_i(t),X_j(s)) &= \operatorname{Cov}\left(\sum_{g\in{\tt OCS}({H}_i)}\frac{\sqrt{\mathcal{A}(g)}\,\mathcal{S}({H}_i,g)}{\mathcal{A}({H}_i)}\bigl(p_{\star}(t)\bigr)^{\mathcal{E}({H}_i) - \mathcal{E}(g)}X^{\star}_{g}(t),\right.\\
    &\qquad\qquad\qquad\left.\sum_{g\in{\tt OCS}({H}_j)}\frac{\sqrt{\mathcal{A}(g)}\,\mathcal{S}({H}_j,g)}{\mathcal{A}({H}_j)}\bigl(p_{\star}(s)\bigr)^{\mathcal{E}({H}_j) - \mathcal{E}(g)}X^{\star}_{g}(t)\right)\\
    &=\sum_{\substack{g_1\in {\tt OCS}({H}_i)\\ g_2\in {\tt OCS}({H}_j)}}\frac{\sqrt{\mathcal{A}(g_1)\,\mathcal{A}(g_2)}}{\mathcal{A}({H}_i)\,\mathcal{A}({H}_j)}\frac{\bigl(p_{\star}(t)\bigr)^{\mathcal{E}({H}_i)}\bigl(p_{\star}(s)\bigr)^{\mathcal{E}({H}_j)}}{\bigl(p_{\star}(t)\bigr)^{\mathcal{E}(g_1)}\bigl(p_{\star}(s)\bigr)^{\mathcal{E}(g_2)}}\\
    &\qquad\qquad\qquad\qquad\qquad\times\mathcal{S}({H}_i,g_1)\,\mathcal{S}({H}_j,g_2)\operatorname{Cov}(X_{g_1}^{\star}(t),X^{\star}_{g_2}(s))\\
    &=\sum_{\substack{g\in {\tt OCS}({H}_i)\cap{\tt OCS}({H}_j)}}{\frac{\mathcal{A}(g)}{\mathcal{A}({H}_i)\mathcal{A}({H}_j)}}\frac{\bigl(p_{\star}(t)\bigr)^{\mathcal{E}({H}_i)}\bigl(p_{\star}(s)\bigr)^{\mathcal{E}({H}_j)}}{\bigl(p_{\star}(t)\,p_{\star}(s)\bigr)^{\mathcal{E}(g)}}\\
    &\qquad\qquad\qquad\qquad\qquad\times\mathcal{S}({H}_i,{g})\,\mathcal{S}({H}_j,{g})\bigl(\kappa_{\star}(t,s)\bigr)^{\mathcal{E}(g)}.
\end{align*}
where in the last line we use the independence of the individual processes {$X^{\star}_{g}(\cdot)$}.
Hence, the claim follows from applying \eqref{GCS_cong}, \eqref{GCS_not_cong} and \nelem{lem:constant_formula}.
\QED

\medskip

\noindent {\it Proof of Lemma \ref{lem:omega_bound}.}
For every edge $(u,v)\in{\tt E}(K_N)$ and time $t^{\star}\in[0,T]$, define the event that the edge $(u,v)$ exists at time $t^{\star}$:
\begin{align*}
    \Omega_{(u,v)}^{+}(t^{\star}) = \{a_{N,u,v}(t^{\star})=1\}.
\end{align*}
Likewise, for every pair of times $t^{\star},t^{\star\star}\in[0,T]$, we define the event that the edge $(u,v)$ switches between times $t^{\star}$ and $t^{\star\star}$:
\begin{align*}
    \Omega_{(u,v)}^{\not =}({t^{\star},t^{\star\star}})= \{a_{N,u,v}(t^{\star})\not= a_{N,u,v}(t^{\star\star})\}.
\end{align*}
By the definition of the processes {$a_{N,u,v}(\cdot)$}, all the events $\Omega_{(u,v)}^{+}(t^{\star})$ and $\Omega_{(u,v)}^{\not =}({t^{\star},t^{\star\star}})$ are independent for different edges $(u,v)$. 
One readily verifies that for every $k\in\sprod{4}$ we can represent the event $\Omega_k$ in terms of events $\Omega^{+}$ and $\Omega^{\neq}$ as follows:
\begin{align}
    \Omega_k = \left(\left(\bigcap_{(u,v)\in{\tt E}(g_k)}\Omega^{+}_{{(u,v)}}(x_k)\right)\cup \left(\bigcap_{(u,v)\in{\tt E}(g_k)}\Omega^{+}_{{(u,v)}}(y_k)\right)\right){\bigcap}\left(\bigcup_{(u,v)\in{\tt E}(g_k)}\Omega^{\not =}_{(u,v)}(x_k,y_k)\right){.}\label{omega_in_omega_0}
\end{align}
As $x_k,y_k\in\{r,s,t\}$ it is convenient to also define the following event:
\begin{align*}
    \Omega^{+}_{(u,v)} = \Omega_{(u,v)}(r)\cup \Omega_{(u,v)}(s)\cup\Omega_{(u,v)}(t).
\end{align*}
Then we can slightly enlarge the event on the right hand side of \eqref{omega_in_omega_0}, obtaining
\begin{align}
    \Omega_k &\subset \left(\bigcap_{(u,v)\in{\tt E}(g_k)}\Omega^{+}_{{(u,v)}}\right){\bigcap}\left(\bigcup_{(u,v)\in{\tt E}(g_k)}\Omega^{\not =}_{(u,v)}(x_k,y_k)\right){.}\label{omega_in_omega}
\end{align}
Based on the inclusions \eqref{omega_in_omega}, we can obtain similar inclusions for the `double event': for any $k_1,k_2\in\sprod{4}$,
\begin{align}
    \Omega_{k_1}\Omega_{k_2}  &\subset \left(\bigcap_{(u,v)\in{\tt E}(g_{k_1}\cup g_{k_2})}\Omega^{+}_{u,v}\right)\left(\bigcup_{(u,v)\in{\tt E}(g_{k_1})}\Omega^{\not =}_{(u,v)}(x_{k_1},y_{k_1})\right)\label{omega_double_bound}.
\end{align}
Along the same lines, for triple event for different $k_1,k_2,k_3\in\sprod{4}$ (without loss of generality we may assume that $k_1\in\{1,2\}$ and $k_2\in\{3,4\}$),
\begin{align}
    \Omega_{k_1}\Omega_{k_2}\Omega_{k_3} &\subset \left(\bigcap_{(u,v)\in{\tt E}(g_{k_1}\cup g_{k_2}\cup g_{k_3})}\Omega^{+}_{{(u,v)}}\right){\bigcap}\left(\bigcup_{\substack{(u_1,v_1)\in{\tt E}(g_{k_1}) \\ (u_2,v_2)\in{\tt E}(g_{k_2})}}\Omega^{\not =}_{(u_1,v_1)}(r,s)\Omega^{\not =}_{(u_2,v_2)}(s,t)\right)\label{omega_triple_bound},
\end{align}
and for the intersection of all four events $k_1,k_2,k_3,k_4\in\sprod{4}$ (again, without loss of generality we may again assume that $k_1\in\{1,2\}$ and $k_2\in\{3,4\}$),
\begin{align}
    \Omega_{k_1}\Omega_{k_2}\Omega_{k_3}\Omega_{k_4} &\subset \left(\bigcap_{(u,v)\in{\tt E}(g_{k_1}\cup g_{k_2}\cup g_{k_3}\cup g_{k_4})}\Omega^{+}_{{(u,v)}}\right){\bigcap}\left(\bigcup_{\substack{(u_1,v_1)\in{\tt E}(g_{k_1}) \\ (u_2,v_2)\in{\tt E}(g_{k_2})}}\Omega^{\not =}_{(u_1,v_1)}(r,s)\Omega^{\not =}_{(u_2,v_2)}(s,t)\right)\label{omega_quad_bound}.
\end{align}
Using \eqref{p_star_bound}, \eqref{claim_1} and \eqref{claim_2}, we thus find that, for any $(u,v)\in{\tt E}(K_N)$ and any $x_k,y_k$,
\begin{align}
    \pk{\Omega^{+}_{(u,v)}}&\leqslant {p_+}\,\varrho_N{,}\label{omega_plus_bound}\\
    \pk{ \Omega_{(u,v)}^{\not =}(x_k,y_k)}&\leqslant {\mathfrak{C}}\,\varrho_N\abs{t-r}\label{omega_neq_bound},\\
    \pk{ \Omega_{(u,v)}^{\not =}(r,s)\Omega_{(u,v)}^{\not =}(s,r)}&\leqslant {\mathfrak{C}}\,\varrho_N\abs{t-r}^2.\label{omega_two_bound}
\end{align}
Hence, applying the bounds \eqref{omega_plus_bound} and \eqref{omega_neq_bound} to the probability of the event \eqref{omega_in_omega}, we find that
\begin{align*}
    \pk{\Omega_k}&\leqslant \sum_{{(u^{\star},v^{\star})\in{\tt E}(g_k)}}\pk{\Omega^{\not =}_{(u^{\star},v^{\star})}(x_k,y_k)\left(\bigcap_{(u,v)\in{\tt E}(g_k)}\Omega^{+}_{{(u,v)}}\right)}\\
    &\leqslant \sum_{{(u^{\star},v^{\star})\in{\tt E}(g_k)}}\pk{\Omega^{\not =}_{(u^{\star},v^{\star})}(x_k,y_k)}\prod_{\substack{(u,v)\in{\tt E}(g_k) \\ (u,v)\not= (u^{\star},v^{\star})}}\pk{\Omega^{+}_{(u,v)}}\\
    &\leqslant {\mathfrak{C}}\,\mathcal{E}(g_k)\,{p_+^{\mathcal{E}(g_k)-1}}\varrho_N^{\mathcal{E}(g_k)}\abs{t-r}\leqslant C_0\varrho_N^{\mathcal{E}(g_k)}\abs{t-r};
\end{align*}
recall that constant $C_0$ was defined in \eqref{C_0_def}. Similarly, applying the same two bounds to the probability of the event \eqref{omega_double_bound}, we have that
\begin{align*}
    \pk{\Omega_{k_1}\Omega_{k_2}}&\leqslant {\mathfrak{C}}\,\mathcal{E}(g_{k_1})\,{p_+^{\mathcal{E}(g_{k_2}\cup g_{k_2})-1}}\varrho_N^{\mathcal{E}(g_{k_1}\cup g_{k_2})}\abs{t-r}\leqslant C_0\varrho_N^{\mathcal{E}(g_{k_1}\cup g_{k_2})}\abs{t-r}.
\end{align*}
Applying now all three bounds $\eqref{omega_plus_bound}$, \eqref{omega_neq_bound} and \eqref{omega_two_bound} to the probability of the event \eqref{omega_triple_bound}, we obtain 
\begin{align*}
    \pk{\Omega_{k_1}\Omega_{k_2}\Omega_{k_3}}&\leqslant \sum_{\substack{(u_1^{\star},v_1^{\star})\in{\tt E}(g_{k_1}) \\ (u_2^{\star},v_2^{\star})\in{\tt E}(g_{k_2})} }\pk{\Omega^{\not =}_{(u_1^{\star},v_1^{\star})}(r,s)\Omega^{\not =}_{(u_2^{\star},v_2^{\star})}(s,t)\left(\bigcap_{(u,v)\in{\tt E}(g_{k_1}\cup g_{k_2}\cup g_3)}\Omega^{+}_{u,v}\right)}\\
    &\leqslant \sum_{\substack{(u_1^{\star},v_1^{\star})\in{\tt E}(g_{k_1}) \\ (u_2^{\star},v_2^{\star})\in{\tt E}(g_{k_2})} }\pk{\Omega^{\not =}_{(u_1^{\star},v_1^{\star})}(r,s)}\pk{\Omega^{\not =}_{(u_2^{\star},v_2^{\star})}(s,t)}\prod_{\substack{(u,v)\in{\tt E}(g_{k_1}\cup g_{k_2}\cup g_{k_3}) \\ (u,v)\not= (u_1^{\star},v_1^{\star}) \\ (u,v)\not= (u_2^{\star},v_2^{\star})}}\pk{\Omega^{+}_{u,v}}\\
    &\qquad + \sum_{(u^{\star},v^{\star})\in {\tt E}({g_{k_1}})\cap {\tt E}({g_{k_2}})}\pk{\Omega^{\not =}_{(u^{\star},v^{\star})}(r,s)\Omega^{\not =}_{(u^{\star},v^{\star})}(s,t)}\prod_{\substack{(u,v)\in{\tt E}(g_{k_1}\cup g_{k_2}\cup g_{k_3}) \\ (u,v)\not= (u^{\star},v^{\star})}}\pk{\Omega^{+}_{u,v}}\\
    &\leqslant {\mathfrak{C}^2}\,\mathcal{E}(g_{k_1})\mathcal{E}(g_{k_2})\,{p_+^{\mathcal{E}(g_{k_1}\cup g_{k_2}\cup g_{k_3})-2}}\,\varrho_N^{\mathcal{E}(g_{k_1}\cup g_{k_2}\cup g_{k_3})}\abs{t-r}^2\\
    &\qquad + {\mathfrak{C}}\,\mathcal{E}(g_{k_1}\cap g_{k_2})\,{p_+^{\mathcal{E}(g_{k_1}\cup g_{k_2}\cup g_{k_3})-1}}\varrho_N^{\mathcal{E}(g_{k_1}\cup g_{k_2}\cup g_{k_3})}\abs{t-r}^2\\
    &\leqslant C_0\varrho_N^{\mathcal{E}(g_{k_1}\cup g_{k_2}\cup g_{k_3})}\abs{t-r}^2{.}
\end{align*}
Along the same lines we can show, relying on \eqref{omega_quad_bound}, that
\begin{align*}
    \pk{\Omega_{k_1}\Omega_{k_2}\Omega_{k_3}\Omega_{k_4}}&\leqslant {\mathfrak{C}^2}\,\mathcal{E}(g_{k_1})\mathcal{E}(g_{k_2})\,{p_+^{\mathcal{E}(g_{k_1}\cup g_{k_2}\cup g_{k_3}\cup g_{k_4})-2}}\varrho_N^{\mathcal{E}(g_{k_1}\cup g_{k_2}\cup g_{k_3}\cup g_{k_4})}\abs{t-r}^2\\
    &\qquad + \,{\mathfrak{C}}\,\mathcal{E}(g_{k_1}\cap g_{k_2})\,{p_+^{\mathcal{E}(g_{k_1}\cup g_{k_2}\cup g_{k_3}\cup g_{k_4})-1}}\varrho_N^{\mathcal{E}(g_{k_1}\cup g_{k_2}\cup g_{k_3}\cup g_{k_4})}\abs{t-r}^2\\
    &\leqslant C_0\,\varrho_N^{\mathcal{E}(g_{k_1}\cup g_{k_2}\cup g_{k_3}\cup g_{k_4})}\abs{t-r}^2{.}
\end{align*}

Hence, the claim of \nelem{lem:omega_bound} follows.
\QED

\medskip

\noindent {\it Proof of Lemma \ref{lem:connected_components}.}
We establish the statement of this lemma for any connected graph ${\tt G}{\,{\subset}\,} K_v$ with at least two vertices; recall that from \eqref{g_condition} we know that $\mathcal{V}({\tt G})\geqslant 2$ for any ${\tt G}\in{\tt CC}(\mathcal{G}(\vk g))$.
As $N\to\infty$,
\begin{align} 
    \E{N^{\mathcal{V}({\tt G})}\prod_{i\in{\tt V}({\tt G})}\frac{\mathcal{F}({H}_{f_i})\xi_{f_i, {N}}(g_i)}{\sqrt{\mathcal{F}^{\rm opt}_N({H}_{f_i})}}}&{=N^{\mathcal{V({\tt g})}}\E{\prod_{i\in{\tt V}({\tt g})}\xi_{f_i,N}(g_i)}\prod_{i\in{\tt V}({\tt G})}\frac{\mathcal{F}_N({H}_{f_i})}{\sqrt{\mathcal{F}_N^{\rm{opt}}({H}_{f_i})}}=:\mathcal{E}_{N}({\tt G})}{.}\label{E_to_F_cal_1}
\end{align}
From \eqref{xi_def} we derive that
\begin{align}
    \prod_{i\in{\tt V}({\tt g})}\xi_{f_i,N}(g_i) = \sum_{k_1,\ldots, k_{\mathcal{V}({\tt g})}=1}^{n}\prod_{i\in{\tt V}({\tt g})}\frac{q_{k,f_i}\bigl(\mathbb{I}\{g_i{\,{\subset}\,} G_N(t_{k_i})\} - \pk{g_i{\,{\subset}\,} G_N(t_{k_i})}\bigr)}{\bigl(p_{\star}(t_{k_i})\bigr)^{\mathcal{E}(g_i)}}{.}\label{xi_prod_decomp}
\end{align} As $N\to\infty$, using Assumption \ref{ass:ass1} we obtain that
\begin{align}
\E{\prod_{i\in{\tt V}({\tt g})}\bigl(\mathbb{I}\{g_i{\,{\subset}\,} G_N(t_{k_i})\} - \pk{g_i{\,{\subset}\,} G_N(t_{k_i})}\bigr)} \leqslant (p_+\varrho_N)^{\mathcal{E}(\cup_{i\in{\tt V}({\tt G})} g_i)}.\label{exp_prod_assympt}
\end{align}
Hence, from \eqref{xi_prod_decomp} and \eqref{exp_prod_assympt} we can see that, as $N\to\infty$,
\begin{align}
    \E{\prod_{i\in{\tt V}({\tt g})}\xi_{f_i,N}(g_i)} = \bigl(C_{{\tt G}}+o(1)\bigr)\varrho_N^{\mathcal{E}(\cup_{i\in{\tt V}({\tt G})} g_i)}\label{xi_prod_assymptotics}
\end{align}
for some constant $C_{{\tt G}}\in[0,\infty)$. Combining \eqref{E_to_F_cal_1} with \eqref{xi_prod_assymptotics} we thus derive that, as $N\to\infty$,
\begin{align} 
    \E{N^{\mathcal{V}({\tt G})}\prod_{i\in{\tt V}({\tt G})}\frac{\mathcal{F}({H}_{f_i})\xi_{f_i, {N}}(g_i)}{\sqrt{\mathcal{F}^{\rm opt}_N({H}_{f_i})}}}&=\bigl(C_{{\tt G}}+o(1)\bigr)\frac{1}{\mathcal{F}_N\left(\cup_{i\in{\tt V}({\tt G})} g_i\right)}\prod_{i\in{\tt V}({\tt G})}\frac{\mathcal{F}_N({H}_{f_i})}{\sqrt{\mathcal{F}_N^{\rm{opt}}({H}_{f_i})}}\label{E_to_F_cal}\\
    &=:\bigl(C_{{\tt G}}+o(1)\bigr)\mathcal{E}_{N}({\tt G}){.}\notag
\end{align}
So, to verify the statement of \nelem{lem:connected_components} it is enough to show that
\begin{align}
    \lim_{N\to\infty}\mathcal{E}_N({\tt G})\in[0,\infty),\label{V(g)=2}
\end{align}
and moreover, if $\mathcal{V}({\tt G})\geqslant 3$ then
\begin{align}
    \lim_{N\to\infty}\mathcal{E}_N({\tt G})=0.\label{V(g)>2}
\end{align}
As we could restrict ourselves to  $\mathcal{V}({\tt G})\geqslant 2$, we are done if we have proved that claim for $\mathcal{V}({\tt G})= 2$ and $\mathcal{V}({\tt G})\geqslant 3$. Consider first the case $\mathcal{V}({\tt G})= 2$. Let ${\tt V}({\tt G}) = \{i,j\}$. Then,
\begin{align*}
    \mathcal{E}_N({\tt G}) &= \frac{\mathcal{F}_N({H}_{f_{i}})\mathcal{F}_N({H}_{f_{j}})}{\mathcal{F}_N(g_{i}\cup g_i)}\frac{1}{\sqrt{\mathcal{F}^{\rm{opt}}_N({H}_{f_i})}\sqrt{\mathcal{F}^{\rm{opt}}_N({H}_{f_j})}}\\
    &=\frac{\mathcal{F}_N(g_i\cap g_j)}{\sqrt{\mathcal{F}^{\rm{opt}}_N({H}_{f_i})}\sqrt{\mathcal{F}^{\rm{opt}}_N({H}_{f_j})}}=\sqrt{\frac{\mathcal{F}_N(g_i\cap g_j)}{\mathcal{F}^{\rm{opt}}_N({H}_{f_i})}}\sqrt{\frac{\mathcal{F}_N(g_i\cap g_j)}{\mathcal{F}^{\rm{opt}}_N({H}_{f_j})}}{.}
\end{align*}
By virtue of \eqref{F_cal_opt_def} (as $(g_i\cap g_j) {\,{\subset}\,} g_i\sim {H}_{f_i}$ and $(g_i\cap g_j) {\,{\subset}\,} g_j\sim {H}_{f_j}$), we obtain that
\begin{align*}
    \lim_{N\to\infty}\mathcal{E}_N({\tt G}) \in\{0,{1}\}\subset[0,\infty),
\end{align*}
and \eqref{V(g)=2} follows. We now proceed to the case $\mathcal{V}({\tt G})\geqslant 3$. We are going to show \eqref{V(g)>2} by induction to $\mathcal{V}({\tt G})$. Our induction step is the following. Assume for any graph ${\tt G}{\,{\subset}\,} K_v$ such that $\mathcal{V}({\tt G})=k-1$ that the statement \eqref{V(g)=2} {applies}. Then for any ${\tt G}{\,{\subset}\,} K_v$ such that 
$\mathcal{V}({\tt G})=k$ we find \eqref{V(g)>2}. 

For any connected graph ${\tt G}{\,{\subset}\,} K_v$ with $k$ vertices, we can find a vertex $j\in{\tt V}({\tt G})$ such that ${\tt G}^{\star} = {\tt G}\setminus\{j\}$ is still a connected graph. Then, 
\begin{align*}
    \mathcal{E}_N({\tt G}) &= \frac{1}{\mathcal{F}_N\left(\cup_{i\in{\tt V}({\tt G})} g_i\right)}\prod_{i\in{\tt V}({\tt G})}\frac{\mathcal{F}_N({H}_{f_i})}{\sqrt{\mathcal{F}_N^{\rm{opt}}({H}_{f_i})}}\\
    &=\mathcal{E}_N({\tt G}^{\star})\frac{\mathcal{F}_N\left(\cup_{i\in{\tt V}({\tt G}^{\star})} g_i\right)}{\mathcal{F}_N\left(\cup_{i\in{\tt V}({\tt G})} g_i\right)}\frac{\mathcal{F}_N({H}_{f_j})}{\sqrt{\mathcal{F}_N^{\rm{opt}}({H}_{f_j})}}\\
    &=\mathcal{E}_N({\tt G}^{\star})\frac{\mathcal{F}_N\biggl(\bigl(\cup_{i\in{\tt V}({\tt G}^{\star})} g_i\bigr) \cap g_j\biggr)}{\mathcal{F}_N\left(g_j\right)}\frac{\mathcal{F}_N({H}_{f_j})}{\sqrt{\mathcal{F}_N^{\rm{opt}}({H}_{f_j})}}\\
    &=\mathcal{E}_N({\tt G}^{\star})\sqrt{\frac{\mathcal{F}_N\biggl(\bigl(\cup_{i\in{\tt V}({\tt G}^{\star})} g_i\bigr) \cap g_j\biggr)}{\mathcal{F}_N^{\rm{opt}}({H}_{f_j})}}\sqrt{\mathcal{F}_N\biggl(\bigl(\cup_{i\in{\tt V}({\tt G}^{\star})} g_i\bigr) \cap g_j\biggr)}{.}
\end{align*}
{Using the induction hypothesis \eqref{V(g)=2} for $\mathcal{E}_N({\tt g}^{\star})$}, from \eqref{F_cal_opt_def} (as $\bigl(\cup_{i\in{\tt V}({\tt G}^{\star})} g_i\bigr) \cap g_j{\,{\subset}\,} g_j\sim {H}_{f_j}$), we obtain
\begin{align*}
    \lim_{N\to\infty}\mathcal{E}_N({\tt G}^{\star})\sqrt{\frac{\mathcal{F}_N\biggl(\bigl(\cup_{i\in{\tt V}({\tt G}^{\star})} g_i\bigr) \cap g_j\biggr)}{\mathcal{F}_N^{\rm{opt}}({H}_{f_j})}}\in[0,\infty).
\end{align*}
Also from \eqref{F_cal_zero_limit} we have
\begin{align*}
    \lim_{N\to\infty}\mathcal{F}_N\biggl(\bigl(\cup_{i\in{\tt V}({\tt G}^{\star})} g_i\bigr) \cap g_j\biggr) = 0.
\end{align*}
Hence, \eqref{V(g)>2} follows for the graph ${\tt G}$, which finishes the proof.
\QED

\medskip

{\noindent {\it Proof of Lemma \ref{lem:G_star_star_cardinality}.} Let us fix ${H}{\,{\subset}\,}{\tt M}_z$ and $\vk g^{\cap} = {(g^{\cap}_1,\ldots, g^{\cap}_{z/2})} \in{\tt OCS}_{\vk f}({H})$. To calculate the cardinality $|\,{\tt G}_{\vk f}^{\star\star}({H},\vk g^{\cap})\,|$, we separate \[\sum_{b=1}^{z}\mathcal{V}({H}_{f_b}) - \sum_{b=1}^{z/2}\mathcal{V}(g_b^{\cap})\] vertices into $z/2$ groups, each of cardinality $\mathcal{V}({H}_{f_{i_b}})+\mathcal{V}({H}_{f_{i_b}})-\mathcal{V}(g^{\cap}_{b})$. It is directly seen that this can be done in
\begin{align}
    \frac{\bigl(\sum_{b=1}^{z}\mathcal{V}({H}_{f_b}) - \sum_{b=1}^{z/2}\mathcal{V}(g_b^{\cap})\bigr)!}{\prod_{b=1}^{z/2}\bigl(\mathcal{V}({H}_{f_{i_b}})+\mathcal{V}({H}_{f_{i_b}})-\mathcal{V}(g^{\cap}_{b})\bigr)!}\label{G_star_star_cardinality_1}
\end{align}
ways, and then for each group separately we should place the graphs $g_{i_b}$ and $g_{j_b}$ in such a way that they cover all $\mathcal{V}({H}_{f_{i_b}})+\mathcal{V}({H}_{f_{i_b}})-\mathcal{V}(g^{\cap}_{b})$ vertices, and $g_{i_b}\cap g_{j_b} = g^{\cap}_b$. From \eqref{C_def}, this can be done in 
\begin{align}
    \bigl(\mathcal{V}({H}_{f_{i_b}})+\mathcal{V}({H}_{f_{i_b}})-\mathcal{V}(g^{\cap}_{b})\bigr)!\,\mathcal{C}({H}_{f_{i_b}},{H}_{f_{j_b}},g^{\cap}_{b})\label{G_star_star_cardinality_2}
\end{align}
ways.
Upon combining \eqref{G_star_star_cardinality_1} and \eqref{G_star_star_cardinality_2}, the claim of \nelem{lem:G_star_star_cardinality} follows{.}
\QED
}

%%%%%%%%%%%%%%%%%%%%%%%%%%%%%%%%%%%%%%%%%%%%%%%%%%%%%%%%%%%%%%%%%%%
%%                                                               %%
%% Use the two commands below for producing your bibliography    %%
%% with bibtex, then comment again the commands and include the  %%
%% content of the .bbl file in this file below the commands.     %%
%%                                                               %%
%%%%%%%%%%%%%%%%%%%%%%%%%%%%%%%%%%%%%%%%%%%%%%%%%%%%%%%%%%%%%%%%%%%

%\bibliographystyle{amsplain}
%\bibliography{reference.bib}

% add below the content of your .bbl file produced by bibtex.

\providecommand{\bysame}{\leavevmode\hbox to3em{\hrulefill}\thinspace}
\providecommand{\MR}{\relax\ifhmode\unskip\space\fi MR }
\providecommand{\MRhref}[2]{%
  \href{http://www.ams.org/mathscinet-getitem?mr=#1}{#2}
}
\providecommand{\href}[2]{#2}

%%%%%%%%%%%%%%%%%%%%%%%%%%%%%%%%%%%%%%%%%%%%%%%%%%%%%%%%%%%%%%%%%%%
%%                                                               %%
%% You may add acknowledgments (optional).                       %%
%%                                                               %%
%%%%%%%%%%%%%%%%%%%%%%%%%%%%%%%%%%%%%%%%%%%%%%%%%%%%%%%%%%%%%%%%%%%
\begin{acks}
{We thank two anonymous referees for their helpful comments.} This research was supported by the European Union’s Horizon 2020 research and innovation programme under the Marie Sklodowska-Curie grant agreement no.\ 945045, and by the NWO Gravitation project NETWORKS under grant agreement no.\ 024.002.003. \includegraphics[height=1em]{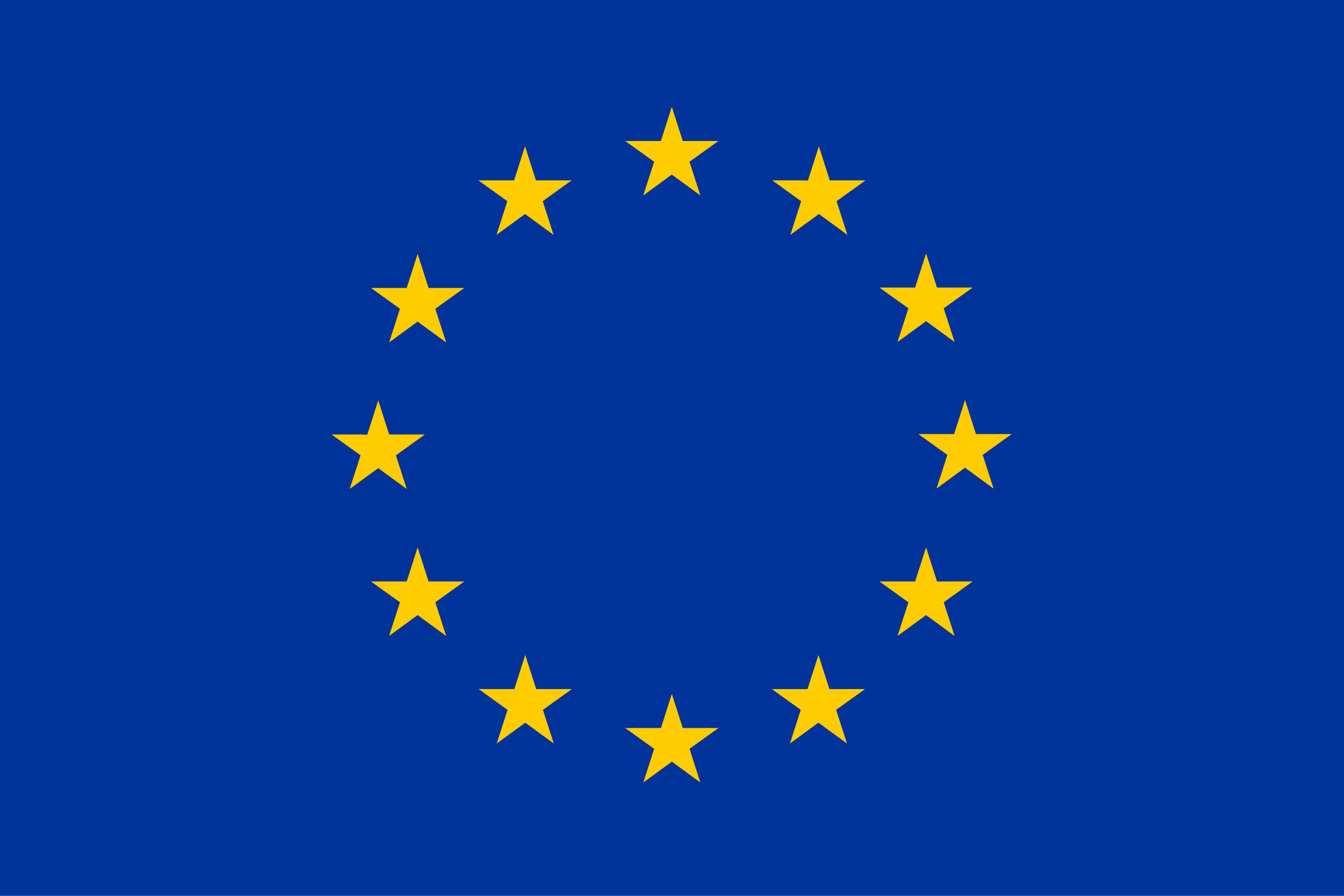} 

\end{acks}

%%%%%%%%%%%%%%%%%%%%%%%%%%%%%%%%%%%%%%%%%%%%%%%%%%%%%%%%%%%%%%%%%%%
%%                                                               %%
%% You have reached the end of your document.                    %%
%%                                                               %%
%%%%%%%%%%%%%%%%%%%%%%%%%%%%%%%%%%%%%%%%%%%%%%%%%%%%%%%%%%%%%%%%%%%

\end{document}